%% file: thesis.tex
\documentclass[12pt]{report}	

\usepackage{utdiss2}
\usepackage{amssymb,amsmath,mathrsfs}
\usepackage{tcolorbox}
\usepackage{amsthm}
\usepackage{enumitem}
\usepackage[all]{xy}
\usepackage[capitalise]{cleveref}
\usepackage{cases}



\usepackage{amsmath,amsthm,amsfonts,amscd} 
\usepackage{eucal} 	 	
\usepackage{verbatim}      	
\usepackage{makeidx}       	
\usepackage{url}		

\author{Ethan Jacob Leeman}  	

\address{2907 West Ave.\\ Austin, Texas 78705}  

\title{Andre-Quillen (co)homology and Equivariant Stable Homotopy Theory}

\supervisor
	{Andrew Blumberg}

\committeemembers
	[David Ben-Zvi]
	[Daniel Allcock]
	{Michael Hill}

\graduationmonth{August}      
\graduationyear{2019}   
\oneandonehalfspacing

\oneandonehalfspacequote

\DeclareMathOperator*{\colim}{colim}
\DeclareMathOperator{\Tor}{Tor}
\DeclareMathOperator{\CoInd}{CoInd}

\DeclareMathOperator{\Ext}{Ext}

\DeclareMathOperator{\coker}{coker}
\DeclareMathOperator{\Exalcomm}{Exalcomm}
\newtheorem{defn}[equation]{Definition}
\newtheorem{theorem}{Theorem}[section]

\newtheorem{definition}{Definition}[section]
\newtheorem{prop}[theorem]{Proposition}
\newtheorem{conj}[theorem]{Conjecture}

\newtheorem{corollary}[theorem]{Corollary}
\newtheorem{cor}[theorem]{Corollary}
\newtheorem{lemma}[theorem]{Lemma}

\newcommand{\term}[2][]{{\bfseries #2}}

\newcommand{\Z}{\mathbb{Z}}

\newcommand{\cat}{\mathsf}

\DeclareMathOperator{\Hom}{Hom}

\def\XXint#1#2#3{{\setbox0=\hbox{$#1{#2#3}{\int}$ }
\vcenter{\hbox{$#2#3$ }}\kern-.6\wd0}}

\newcommand{\Id}{\operatorname{Id}}
\newcommand{\id}{\operatorname{Id}}
\newcommand{\Der}{\operatorname{Der}}

\newcommand{\latexe}{{\LaTeX\kern.125em2%
                      \lower.5ex\hbox{$\varepsilon$}}}

\chardef\bslash=`\\	

\makeatletter		
\def\square{\RIfM@\bgroup\else$\bgroup\aftergroup$\fi
  \vcenter{\hrule\hbox{\vrule\@height.6em\kern.6em\vrule}%
                                              \hrule}\egroup}
\makeatother		

\makeindex    

\begin{document}

\copyrightpage          

%
%
%
\commcertpage           

\titlepage              

%

\begin{acknowledgments}		
\index{Acknowledgments@\emph{Acknowledgments}}%
\noindent From formal to personal:

I would like to thank the committee members for reviewing this dissertation and attending my oral examination. From the committee members I would like to single out Mike Hill, whose work was the launching point of this dissertation's topic and who supported me getting a grasp on this area of mathematics. I would like to acknowledge both past and present Graduate Advisors (Dan Knopf, Thomas Chen, Kui Ren, Tim Perutz) and Graduate Program Administrators (Sandra Catlett, Elisa Bass, Jenny Kondo) for guiding me through the more formal procedures of this program.

The amount of credit due to my advisor Andrew Blumberg cannot be exaggerated. It perhaps goes without saying a thesis would not exist without its supervisor, and this one is certainly no exception. A big thank you for taking me on as a student, guiding me through this project, and imparting your knowledge, both mathematical and not.

I owe a great deal of gratitude to my previous mathematics teachers who indirectly but undoubtably made this work possible. A special thank you to Erica Voolich, David Srebnick, Debbie Seidell and Alex Martsinkovsky. Without their mentorship I would not be here as a graduate student.

A warm thank you to the numerous friends I have made in my time as a graduate student. To name only a few: Eric Korman, Mike Lock, Matt Novack, Cornelia Mihaila, Chris Kennedy, Nicky Reyes, Marjorie Drake, Luis Duque, Val Zakharevich, and Zachary Scherr.

Lastly, I'd like give the largest thanks to my fianc\'ee Misola and my family, whose continuous love supported me in these years. This page is too narrow to state the appreciation they deserve.
\end{acknowledgments}

%
\utabstract
\index{Abstract}%
\indent
Andr\'e and Quillen introduced a (co)homology theory for augmented commutative rings. Strickland \cite{StricklandTambara} initially proposed some issues with the analogue of the abelianization functor in the equivariant setting. These were resolved by Hill \cite{HillAndre} who further gave the notion of a genuine derivation and a module of K\"ahler differentials. We build on this endeavor by expanding to incomplete Tambara functors, introducing the cotangent complex and its various properties, and producing an analogue of the fundamental spectral sequence.
\tableofcontents   


%
%
\include{Intro}

\include{HomChapter}

\include{Modules}

\include{ExtTor}
\include{AQ}
\include{kahlerdiff2}
\include{diffprop}
\include{KahlerOfFree2}

\include{CotangentComplex}

\include{CotangentComplexProp}

\include{SpectralSeq}

\include{Convergence}

%
%

%
%
%

\nocite{*}      
\bibliographystyle{plain}  
\bibliography{diss2}        
\index{Bibliography@\emph{Bibliography}}





\end{document}

%% file: Intro.tex
\chapter{Introduction and Beginning Terminology}
In order to introduce this subject, we will draw parallels between classical homotopy theory and equivariant homotopy theory. Equivariant stable homotopy theory was first created by Segal \cite{segal1970equivariant} and has evolved by the contributions of many mathematicians since. One good source of notes for context is \cite{arun}. Let $G$ be a finite group (in many instances one can take a compact Lie group, but for simplicity we restrict to a finite group). Suppose instead of the category of topological spaces, our interest is on the category $\cat{GTop}$ of topological spaces along with a symmetry by $G$, also known as $G$-spaces. Then $\cat{GTop}$ has as objects the compactly generated weak Hausdorff topological spaces with a $G$-action. The homomorphism sets are the continuous maps that are equivariant, that is $f(g \cdot x) = g \cdot f(x)$ for all $g \in G$ and $x \in X$. There are appropriate notions of a homotopy of maps and homotopy equivalences. Additionally, there is a based version of this category $\cat{GTop}_*.$
The theory of equivariant stable homotopy theory is the application of stable homotopy theory in this realm, which we summarize briefly. Upon doing so, we realize the correct analogues of abelian groups and commutative rings are Mackey functors and Tambara functors. One can then ask what commutative algebra looks like on these new objects. Andr\'e-Quillen cohomology is one such commutative algebra construction, which is the purpose of this thesis.

A natural first question for $\cat{GTop}$ is: what are the appropriate generalizations of CW-complexes? CW-complexes are ubiquitous as test objects and in constructions. The correct answer is to allow attaching cells of the form $(G/H \times D^{n+1},G/H \times S^n).$ That is, a $G$-CW complex is a sequential colimit of spaces $X_n$ where $X_{n+1}$ is the pushout of the following diagram:
\[\xymatrixrowsep{.4in}
\xymatrixcolsep{.5in}\xymatrix{
 \coprod G/H \times S^n \ar[d] \ar[r]&  X_n \ar[d]\\
 \coprod G/H \times D^{n+1} \ar[r]& X_{n+1}\\
}\] 
Then a natural second question is: what are the appropriate generalizations for homotopy groups of a $G$-space? From the definition of $G$-CW complexes the definition must be $$\pi_n^H(X) := [G/H \times S^n, X] \cong [S^n, X^H].$$ So homotopy groups are indexed by subgroups of $G$, as well as by the natural numbers. The following generalization of the Whitehead theorem gives us good justification that we have the right definitions.
\begin{prop}
Let $f: X \to Y$ be a map between $G$-spaces that induces an isomorphism on $\pi_n^H$ for all $H \subset G, n \in \mathbb{N}.$ Then $f$ is a homotopy equivalence.
\end{prop}
Homotopy groups have more structure. Suppose $f: G/H \to G/K$. The map is determined by $eH \mapsto gK$ for some $g \in G$. This requires that $gHg^{-1} \subset K$. This induces $f^*: X^K \to X^H$ by the map $x \mapsto gx$, which induces a map $\pi_n^K(X) \to \pi_n^H(X).$ Suppose $n \geq 2$. Let $O_G$ the orbit category, the full subcategory of $\cat{GTop}$ with objects the orbits $G/H.$ Then $\pi_n^*: O_G^{op} \to \cat{Ab}$ for $n \geq 2.$ Any contravariant functor $O_G^{op} \to \cat{Ab}$ is called a coefficient system.

Homotopy groups describe the coefficients of cohomology theories. Most obviously, the obstruction classes of obstruction theory take values in homotopy groups of the target. But most directly, for a coefficient system $\underline{M}$ there is an Eilenberg-Mac Lane G-space $K(\underline{M},n)$ which has homotopy groups 0 except at degree $n$ for which the homotopy groups are $\underline{M}$. Then a cohomology theory, the Bredon cohomology, is $\tilde{H}_G^n(X;\underline{M}) := [X, K(\underline{M},n)]_G$ which was first introduced in \cite{bredon2006equivariant}. While resembling the ordinary non-equivariant cohomology theory, this cohomology theory is significantly more difficult to compute.

Naturally, we should be considering the generalization of spectra, as these objects represent generalized cohomology theories. The concept of a genuine $G$-spectra is due to \cite{gaunce2006equivariant}. To summarize the key point, non-equivariant prespectra are given by a sequence of spaces $X_n$ with suspension maps that we write as $S^m \wedge X_n \to X_{m+n}$. But equivariantly, this construction seems a bit naive, as the sphere could include a $G$-action itself. Hence, this would be a naive $G$-prespectrum, and a naive $G$-spectrum would be one satisfying the necessary equivalence involving the loop space operator.

Let $V$ be any $G$-representation. Then we could take $S^V$, the one point compactification of $V$. Then we could define a genuine $G$-prespectrum to be a $G$-space $X(V)$ for every finite-dimensional representation of $V$ with maps $S^W \wedge X(V) \to X (V \oplus W).$ The appropriate cohomology theory we get, instead of being $\mathbb{Z}$-graded will be $RO(G)$-graded. For the purposes of this introduction, we will not delve into the definitions of an $RO(G)$-graded cohomology theory, a correct point-set model (symmetric or orthogonal spectra) for spectra in order for the sphere spectrum to be commutative, or any discussion about universes in order to interpolate between naive and genuine spectra.

However, in the category of spectra, our category of orbits has more maps. Suppose $H/K$ embeds into V. Then by the Pontryagin-Thom construction \cite{thom}, there is a sequence of maps $$ S^V \to T(\nu) \to H/K \wedge S^V $$ which by induction gives a map $$G/H \wedge S^V \to G/K \wedge S^V.$$ There is no map $G/H \to G/K$ in $\cat{GTop}$ but these maps do exist between the suspension spectra. As a result, we must replace the orbit category with another suitable category, which we describe in the next section, the Burnside category. Additive maps from the Burnside category to abelian groups are called Mackey functors, and they generalize abelian groups because the homotopy groups of genuine $G$-spectra are Mackey functors.

Finally, we must consider ring structures on $G$-spectra. Non-equivariantly, we have a notion of $E_\infty$-operads which describes the associative and commutative relations in spectra, and an algebra over these operads are called $E_\infty$-ring spectra. In this setting, all the $E_\infty$-operads are homotopy equivalent, giving equivalent definitions of $E_\infty$-ring spectra. In the equivariant setting, there are $N_\infty$-operads which give rise to $N_\infty$ ring $G$-spectra. $N_\infty$ operads are not all homotopy equivalent, but the homotopy classes are given by indexing systems $\mathcal{O}$. Just as the homotopy groups of $E_\infty$ ring spectra are commutative rings, the homotopy groups of $N_\infty$ ring $G$-spectra are $\mathcal{O}$-Tambara functors. $N_\infty$-opeards were first introduced in \cite{blumberg2015operadic} and followed up in \cite{incomplete}. Tambara functors were first introduced by \cite{tambara1993multiplicative}, generalizing the Even's norm \cite{evens1991cohomology} in group cohomology before seeing use in equivariant stable homotopy theory leading to \cite{hill2016nonexistence} which resolved the Kervaire invariant 1 conjecture.

Now that we have Mackey functors and incomplete Tambara functors, our generalizations of abelian groups and commutative rings, a natural question is how much of commutative algebra can we extend to this equivariant regime? This thesis aims to do one such construction, the Andr\'e-Quillen cohomology of augmented rings.

In the 1960's both Andr\'e and Quillen \cite{andre1974homologie} \cite{QuillenMimeo} introduced a (co)homology theory for maps of commutative rings. These determine the obstruction theory for deformations of rings. Many of the properties of (co)homology theories from topology have analogues in Andr\'e-Quillen (co)homology. 

Since then there have been many generalizations of this (co)homology theory. One of note is that Basterra \cite{basterra1999andre} produced the theory of Topological Andr\'e-Quillen homology on commutative ring spectra and went on to show with Mandell \cite{basterra2005homology} that this homology theory is the only homology theory on commutative ring spectra.

This path was first started out by Strickland \cite{StricklandTambara}, who proposed some issues, particularly that the square zero extension from modules to the abelian group objects of algebras was not essentially surjective. This issue was resolved by Hill \cite{HillAndre} by creating the Mackey functor objects in the appropriate algebra category. We discuss these in a later chapter. Hill additionally created the definition of a derivation and the module of K\"ahler differentials. We continue this train of thought by incorporating incomplete Tambara functors, introducing the cotangent complex and its various properties, and producing an analogue of the fundamental spectral sequence.

\section{Mackey Functors, (Incomplete) Tambara functors}
We discuss the objects of interest in terms of polynomial categories, as in \cite{incomplete}.
\begin{definition}
Let $\cat{C}$ be a locally Cartesian closed category. Call $\cat{P}^{\cat{C}}$ the category with objects the same as $\cat{C}$ and morphisms $\cat{P}^{\cat{C}}(X,Y)$ the isomorphism classes of polynomials  \[\xymatrixrowsep{.5in}
\xymatrixcolsep{.25in}\xymatrix{
X & S \ar[l] \ar[r] & T \ar[r] & Y.
}\] Two polynomials are isomorphic if there are isomorphisms $S \overset{\cong}\to S'$ and $T \overset{\cong}\to T'$ that make the following diagram commute
\[\xymatrixrowsep{.25in}
\xymatrixcolsep{.5in}\xymatrix{
 & S \ar[ld] \ar[r] \ar[dd]^\cong & T \ar[rd] \ar[dd]^\cong \\
 X & & & Y\\
 & S' \ar[lu] \ar[r] & T' \ar[ru] \\
}\] 
\end{definition}
We will typically consider $\cat{C}$ to be $\cat{Set}^G$, the finite $G$-sets with equivariant maps. In that case we suppress the superscript and write $\cat{P}$. To describe the compositions of polynomials, we give various double coset rules (pullback diagrams) and Tambara reciprocity relations (exponential diagrams).

\begin{definition}If $f \in \cat{C}(S,T)$ then we give the following morphisms of $\cat{P}^{\cat{C}}$ names 
$$R_f = [T \overset{f} \leftarrow S \overset{1} \rightarrow S \overset{1} \rightarrow S]$$
$$N_f = [S \overset{1} \leftarrow S \overset{f} \rightarrow T \overset{1} \rightarrow T]$$
$$T_f = [S \overset{1} \leftarrow S \overset{1} \rightarrow S \overset{f} \rightarrow T]$$ which are respectively called the \term{restriction, transfer, and norm of $f$}.
\end{definition}
Any homomorphism in $\cat{P}^{\cat{C}}$ can be written in the form of $T_h \circ N_g \circ R_f$ for some equivariant maps $f,g,h$, which we call \term{TNR-form}. Composition of polynomials can be computed by describing commutation relations that put any composition $(T_h \circ N_g \circ R_f) \circ (T_{h'} \circ N_{g'} \circ R_{f'})$ into TNR-form. \cite{StricklandTambara} is one of many good sources for these relations:
\begin{prop} $$T_h \circ T_h' = T_{h \circ h'} $$ 
$$N_g \circ N_{g'} = T_{g \circ g'} $$
$$R_f \circ R_{f'} = R_{f' \circ f}$$
\end{prop}
\begin{prop}
If we have the following pullback diagram 
\[\xymatrixrowsep{.5in}
\xymatrixcolsep{.5in}\xymatrix{
Y' \underset{Y} \times X \ar[r]^{g'} \ar[d]^{f'} & X \ar[d]^f\\
Y' \ar[r]^g& Y
}\] then 
$$R_f \circ N_g = N_{g'} \circ R_{f'} $$
$$R_f \circ T_g = T_{g'} \circ R_{f'}. $$ These are also called the double coset formulas.
\end{prop}
\begin{prop}
Suppose $A \overset{h} \to X \overset{g} \to Y$ are maps in $\cat{Set}^G$. Then define $$\prod_g A := \left\{(y,s: g^{-1}(y) \to A) \, | \, y \in Y, \, h \circ s = \Id_{g^{-1}(y)} \right\}. $$ This gives rise to the following pullback diagram 
\[\xymatrixrowsep{.5in}
\xymatrixcolsep{.5in}\xymatrix{
X \ar[d]^g& A\ar[l]^h & X \underset{Y} \times \prod_g A \ar[l]^(.55){f'} \ar[d]^{g'} \\
Y && \prod_g A \ar[ll]^{h'}
}\] called an exponential diagram,
where \begin{eqnarray*}
h': (y,s) &\mapsto& y,\\ f':(x,y,s) &\mapsto& s(x),\\  g':(x,y,s) &\mapsto& (y,s).
\end{eqnarray*} 
We can rewrite $$X \underset{Y} \times \prod_g A  \cong \left\{(x,s: g^{-1}(g(y)) \to A) \, | \,x \in X, \,h \circ s = \Id_{g^{-1}(g(x))} \right\}.$$
Then $$N_g \circ T_h = R_{f'} \circ N_{g'} \circ T_{h'}$$ which we call the Tambara relation.
\end{prop}
In general, $\Pi_f$ can be replaced for the right adjoint of $f^*$, the pullback functor $\cat{C}/Y \to \cat{C}/X.$

Now we consider the subcategories of $\cat{P}^{\cat{C}}$ given by restricting the exponents. We use the following proposition of \cite{incomplete}:
\begin{definition}
A a subcategory $\cat{D}$ of $\cat{C}$ is \term{pullback stable} if given the pullback diagram \[\xymatrixrowsep{.5in}
\xymatrixcolsep{.5in}\xymatrix{
Y' \underset{Y} \times X \ar[r] \ar[d]^{f'} & X \ar[d]^f\\
Y' \ar[r]& Y
}\]  and $f \in \cat{D}(X,Y)$, then $f' \in \cat{D}(Y' \underset{Y} \times X,Y').$
\end{definition}
\begin{prop}
Let $\cat{D}$ be a wide, pullback stable, finite coproduct complete subcategory of $\cat{C}$. Then the subgraph of $\cat{P}^{\cat{C}}$ of morphisms $T_h N_g R_f$ where $g \in \cat{D}$ is a subcategory called \term{the category of polynomials with exponents restricted to $\cat{D}$} and is denoted $\cat{P}_{\cat{D}}$.
\end{prop}

As described in \cite{incomplete}, there is a bijection between the wide, pullback stable, coproduct complete categories and indexing systems. We recite the key definitions and results now:
\begin{definition}
Let $\underline{Set}: Orb_G^{op} \to Sym$ be the functor that sends every orbit $G/H$ to the category of finite $H$-sets. The monoidal product is disjoint union of $H$-sets. An indexing system $\mathcal{O}$ is a subfunctor of $\underline{Set}$ that 
\begin{enumerate}
\item contains all trivial sets,
\item is closed under finite limits, and
\item is closed under self-induction: if $H/K \in \mathcal{O}(G/H), \, T \in \mathcal{O}(G/K),$ then $H \underset{K}\times T \in \mathcal{O}(G/H).$
\end{enumerate}
\end{definition}
\begin{theorem}
Given an indexing system $\mathcal{O}$, let $\cat{Set}_{\mathcal{O}}^G$ be the subgraph of $\cat{Set}^G$ where $f:S\to T$ is a morphism in $\cat{Set}_{\mathcal{O}}^G$ if and only if for all $s \in S$, $$G_{f(s)}/G_s \in \mathcal{O}(G_{f(s)}).$$ This is a wide, pullback stable, finite coproduct complete subcategory of $\cat{Set}^G$. This map from the poset of indexing systems and the poset of wide, pullback stable, finite coproduct complete subcategories of $\cat{Set}^G$ is an isomorphism.
\end{theorem}
\begin{definition}
For an indexing system, the category of $\mathcal{O}$-Tambara functors is the category of functors from $\cat{P}_{\cat{D}}$ to abelian groups, where $\cat{D}$ is the category given by $\mathcal{O}$ and the bijection in the preceding theorem.
\end{definition}
We can define Mackey functors using this terminology: let $\cat{Isom}$ be the category of isomorphisms of $\cat{Set}^G$. Then \cite{incomplete} shows that functors $\cat{P}_{\cat{Isom}} \to \cat{Ab}$ are Mackey functors, denoted as the category $\cat{Mac}$. The more standard way of defining Mackey functors is functors $\cat{B}_G \to \cat{Ab}$, where $\cat{B}_G$ is the Burnside category (polynomials where the norm map is the identity).
 It is often easier to define Mackey functors using the orbit categories, as in \cite{Mazur}, which we repeat.
\begin{definition}
A Mackey functor $\underline{M}$ consists of a collection of abelian groups $\underline{M}(G/H)$ with transfer maps $T_K^G: \underline{M}(G/K) \to \underline{M}(G/H)$ and restriction maps $R_K^H: \underline{M}(G/H) \to \underline{M}(G/K)$ for all subgroups $K < H \leq G$ such that 
\begin{enumerate}
\item If $K' < K < H$ then $T_{K'}^H = T_K^H T_{K'}^K$ and $R_{K'}^H = R_{K'}^KR_K^H.$
\item if $K < H \leq G$, there is an action of $N_H(K)/K$ on $\underline{M}(G/K)$ such that if $\gamma \in N_H(K)/K$, then $T_K^H(\gamma \cdot (-)) = T_K^H(-)$ and $\gamma \cdot R_K^H(-) = R_K^H(-).$
\item If $K,K'$ are subgroups of $H$, we have the double coset formula $$R_{K'}^H T_K^H(x) = \sum_{\gamma \in N_H(K')/K'} \gamma \cdot T_{K' \cap K}^{K'}(x).  $$
\end{enumerate}
\end{definition}
We can similarly make a definition of incomplete Tambara functors:
\begin{definition}
An $\mathcal{O}$-Tambara functor $\underline{S}$ consists of a collection of commutative rings $\underline{S}(G/H)$ with transfer maps $T_K^G: \underline{S}(G/K) \to \underline{S}(G/H)$ and restriction maps $R_K^H: \underline{S}(G/H) \to \underline{S}(G/K)$ for all subgroups $K < H \leq G$, as well as norm maps $N_K^H:  \underline{S}(G/K) \to \underline{S}(G/H)$ when $\cat{D}$ contains the map $G/K \to G/H$. The transfer and restriction maps have all the properties that make $\underline{S}$ a Mackey functor. The restriction is a map of commutative monoids. We have the additional properties regarding norms:
\begin{enumerate}
\item If $K,K'$ are subgroups of $H$, we have the double coset formula $$R_{K'}^H N_K^H(x) = \sum_{\gamma \in N_H(K')/K'} \gamma \cdot N_{K' \cap K}^{K'}(x).  $$
\item The norm maps commute with transfers by the Tambara reciprocity relation.
\end{enumerate}
\end{definition}

Mackey functors have a symmetric moniodal product, originally described by Day in \cite{Day}, given by the coend $$\underline{M} \boxtimes \underline{N}(X) = \int^{(Y,Z) \in B_G \times B_G} \underline{M}(Y) \otimes \underline{N}(Z) \otimes B_G(X,Y \times Z)$$ sometimes called the box product or tensor product if the meaning is clear. The Burnside Mackey functor $\underline{A} := \Hom_{B_G^{op}}(-,G/G)$ is the unit. $\underline{A}(G/H)$ is the Grothendieck group of finite $H$-sets with coproduct acting as the symmetric monoidal product. The transfer and restriction are induction and restriction. \cite{StricklandTambara} shows that a commutative monoid in the category of Mackey functors is a Green functor, or an incomplete Tambara functor with the trivial indexing system. In this case $\cat{D}$ contains all the maps $f: U \to V$ such that $f(u)$ has the same isotropy group as $u$. We have the following theorem from \cite{incomplete} categorizing incomplete Tambara functors.
\begin{theorem}
For an indexing system $\mathcal{O}$, an $\mathcal{O}$-Tambara functor is a commutative Green functor $\underline{R}$ together with norm maps of multiplicative monoids $N_H^K: \underline{R}(G/H) \to \underline{R}(G/K)$ for each $G/H \to G/K \in Orb_{\mathcal{O}}$ that satisfy the double coset and Tambara relations. 
\end{theorem}
\cite{StricklandTambara} has the following useful description of elements in $\underline{M} \boxtimes \underline{N}$.
\begin{prop}
Fix $\underline{M},\underline{N}$ Mackey functors for a finite $G$-set $X$. Let $\mathcal{E}$ be the set of quadruples $(U,p,m,n)$ where $p:U \to X$, $(m,n)\in \underline{M}(U) \times \underline{N}(U)$. Let $\sim$ be the smallest equivalence relation on $\mathcal{E}$ such that 
\begin{enumerate}
\item For $U' \overset{r} \to U \overset{q} \to X$, $(m',n) \in \underline{M}(U') \times \underline{N}(U)$ then  $$(U', qr, m',R_r(n)) \sim (V,q,T_r(m'),n)$$
\item For $U' \overset{r} \to U \overset{q} \to X$, $(m,n') \in \underline{M}(U) \times \underline{N}(U')$ then $$(U', qr, R_r(m),n') \sim (V,q,m,T_r(n')).$$
\end{enumerate}
Define the map $\mathcal{E} \to (\underline{M} \boxtimes \underline{N})(X)$ by $(U,p,m,n) \mapsto T_p(m \otimes n)$. Then $\mathcal{E}/\sim \, \overset{\cong}\longrightarrow (\underline{M} \boxtimes \underline{N})(X)$ is a bijection.
\end{prop}
Lastly, the arguments of \cite{StricklandTambara} show that the box product of two incomplete Tambara functors is an incomplete Tambara functor, and the box product is the coproduct in incomplete Tambara functors. The fibered box product $(-) \underset{\underline{A}} \boxtimes (-)$ is the fibered coproduct, which has the same description as the above proposition along with the congruences $(am \otimes n) \sim (m \otimes an)$ for $a \in \underline{A}(U)$, $m \in \underline{M}(U), n \in \underline{N}(U).$

%% file: HomChapter.tex
\chapter{Homological Algebra}\label{HomChapter}
\index{Homological Algebra@\emph{Homological Algebra}}%

In this section, we include some properties that are basic to the field of homological algebra but use the appropriate categories from equivariant stable homotopy theory. We consider Mackey functors, $\mathcal{O}$-Tambara functors, and modules over $\mathcal{O}$-Tambara functors rather than abelian groups, commutative rings, and modules over commutative rings respectively. We show that the category of modules is an abelian category satisfying the necessary properties to perform homological algebra. Additionally, we will discuss free, projective, and flat modules and give some examples.

\section{Closed Symmetric Monoidal Category}
Much of this information is a synthesized version of \cite{Lewis80} and \cite{zeng}.
Given Mackey functors $\underline{M}$ and $\underline{N}$, we have a tensor product $\underline{M} \boxtimes \underline{N}$. To reiterate, if $\underline{M}, \underline{N}$ are functors from $\cat{B_G}$ to $\cat{Ab}$ then 
$$\underline{M} \boxtimes \underline{N}(X) = \int^{(Y,Z) \in \cat{B_G} \times \cat{B_G}} \underline{M}(Y) \otimes \underline{N}(Z) \otimes \cat{B_G}(X,Y \times Z).$$ Day \cite{Day} showed that this symmetric monoidal product makes Mackey functors a closed monoidal category. The unit is the Burnside ring $\underline{A} = \underline{B}_G(-,G/G).$ The internal hom is the right adjoint 
$$\cat{Mac}(\underline{A} \boxtimes \underline{B}, \underline{C}) \cong \cat{Mac}(\underline{A}, [\underline{B},\underline{C}]_{Day}) $$ and is given by the end $$[\underline{M}, \underline{N}]_{Day}(X) = \int_{(Y,Z) \in \cat{B_G} \times \cat{B_G}} \Hom( \cat{B_G}(Z,X \times Y)  \otimes \underline{M}(Y), \underline{N}(Z)).$$ We will mostly denote $[-,-]_{Day}$ by $\underline{\Hom}[-,-].$
Let us review the adjunction as relayed in \cite{zeng} from \cite{Lewis80}:
\begin{eqnarray*}
\cat{Mac}(\underline{A},\underline{\Hom}[\underline{B},\underline{C}]) &\cong& \int_X \Hom[\underline{A}(X), \int_{(Y,Z) \in \cat{B_G} \times \cat{B_G}} \Hom( \cat{B_G}( Z, X \times Y) \otimes \underline{B}(Y), \underline{C}(Z)) ] \\
&\cong& \int_X  \int_{(Y,Z) \in \cat{B_G} \times \cat{B_G}} \Hom[ \cat{B_G}(Z,X \times Y) \otimes \underline{A}(X) \otimes \underline{B}(Y), \underline{C}(Z) ]\\
&\cong& \int_Z \Hom[  \int^{(X,Y) \in \cat{B_G} \times \cat{B_G}}  \cat{B_G}(Z,X \times Y) \otimes \underline{A}(X) \otimes \underline{B}(Y), \underline{C}(Z) ]\\
&\cong& \cat{Mac}(\underline{A} \boxtimes \underline{B}, \underline{C}).
\end{eqnarray*}
Because $$\int_Z \Hom[\cat{B_G}(Z, X \times Y), \underline{N}(Z)] \cong \underline{N}(X \times Y)$$ and dually $$\int^Y \cat{B_G}(Z, X \times Y) \otimes \underline{M}(Y) \cong \underline{M}(X \times Z)$$
we get that $$\underline{\Hom}(\underline{M},\underline{N})(X) \cong \Hom(\underline{M}, \underline{N}_X) \cong \Hom(\underline{M}_X,\underline{N})$$ where $\underline{N}_X$ is defined by 
\begin{eqnarray*}
\underline{N}_X(Z) &:=& \underline{N}(Z \times X).
\end{eqnarray*}
In particular, we also have the isomorphisms \begin{eqnarray*}
\cat{B_G}(X,-) \boxtimes \underline{A}&\cong& \underline{A}(X \times -) \\
&\cong& \underline{\Hom}(\cat{B_G}(X,-),\underline{A})
\end{eqnarray*}
\begin{eqnarray*}
\cat{B_G}(X,-) \boxtimes \cat{B_G}(Y,-) &\cong& \cat{B_G}(X \times Y,-) \\
&\cong& \underline{\Hom}(\cat{B_G}(X,-),\cat{B_G}(Y,-)).
\end{eqnarray*}
One can also compute using the adjunction directly
\begin{eqnarray*}
\underline{\Hom}(\underline{M},\underline{N})(X)&=& \cat{Mac}( \cat{B_G}(X,-), [\underline{M},\underline{N}]_{Day}) \\
&=& \cat{Mac}(  \cat{B_G}(X,-) \boxtimes \underline{M}, \underline{N})
\end{eqnarray*}

The box product gives rise to the following definition of Green functors and $\underline{R}$-modules from \cite{StricklandTambara}:
\begin{definition}
A Green functor is a commutative monoid in $(\cat{Mac},\boxtimes,\underline{A}).$
\end{definition}
As \cite{incomplete} shows, all incomplete $\mathit{O}$-Tambara functors are Green functors.

%% file: Modules.tex

\section{Modules}

\begin{definition}
If $\underline{M}$ is a Mackey functor and $\underline{S}$ is an incomplete Tambara functor, then an $\underline{S}$-module structure on $\underline{M}$ is a map $\nu: \underline{S} \boxtimes \underline{M} \to \underline{M}$ such that the diagrams involving $\eta: \underline{M} \to \underline{S} \boxtimes \underline{M}$ and $\mu: \underline{S} \boxtimes \underline{S} \to \underline{S}$ commute
 \[\xymatrixrowsep{.5in}
\xymatrixcolsep{.5in}\xymatrix{
\underline{S} \boxtimes \underline{S} \boxtimes \underline{M} \ar[d]^{1 \boxtimes \eta} \ar[r]^{\mu \boxtimes 1} & \underline{S} \boxtimes \underline{M} \ar[d]^\eta  &\underline{M} \ar[l]_{\eta \boxtimes 1} \ar[dl]^1 \\
\underline{S} \boxtimes \underline{M} \ar[r]^\eta& \underline{M}
  }\]
\end{definition}
An equivalent definition, also shown in \cite{StricklandTambara}, is giving $\underline{M}(X)$ an $\underline{S}(X)$-module structure such that for all $f: X \to Y,$
\begin{enumerate}
\item for all $s_y \in \underline{S}(Y)$ and $ m_y \in \underline{M}(Y)$, $R_f(s_y \cdot m_y) = R_f(s_y) \cdot R_f(m_y);$ 
\item for all $s_y \in \underline{S}(Y)$ and $ m_x \in \underline{M}(X)$, $s_y \cdot T_f(m_x) = T_f(R_f(s_y) \cdot m_x);$ and
\item for all $s_x \in \underline{S}(X)$ and $ m_y \in \underline{M}(Y)$, $T_f(s_x) \cdot m_y = T_f(s_x \cdot R_f(m_y))$.
\end{enumerate}

Now suppose that $\underline{M}$ is also an $\underline{R}$-module and $\underline{N}$ is an $(\underline{R},\underline{S})$-bimodule, and $\underline{C}$ is an $\underline{S}$-module.
Let $\underline{M} \underset{\underline{R}}\boxtimes \underline{N}$ be the coequalizer of the two maps $\underline{R}$ multiplication maps
 \[\xymatrixrowsep{.5in}
\xymatrixcolsep{.25in}\xymatrix{
\underline{M} \boxtimes \underline{R} \boxtimes \underline{N} \ar@<-.5ex>[r] \ar@<.5ex>[r]& \underline{M} \boxtimes \underline{N} \ar[r]& \underline{M} \underset{\underline{R}}\boxtimes \underline{N}.
}\]
This operation is the tensor product that makes $\cat{\underline{R}Mod}$ a closed symmetric monoidal category. The adjoint is the equalizer
 \[\xymatrixrowsep{.5in}
\xymatrixcolsep{.25in}\xymatrix{
\underline{\Hom}_{\underline{S}}[\underline{B},\underline{C}] \ar[r]& \underline{\Hom}[\underline{B},\underline{C}] \ar@<-.5ex>[r] \ar@<.5ex>[r]& \underline{\Hom}[\underline{B} \boxtimes \underline{S},\underline{C}] \cong \underline{\Hom}[\underline{B},\underline{\Hom}[\underline{S},\underline{C}]].
}\]
 where the two maps are given by the maps $\underline{B} \boxtimes \underline{S} \to \underline{B}$ and $\underline{C} \to \underline{\Hom}(\underline{S},\underline{C})$ the adjoint of the multiplication of $\underline{C}.$ We cite the following properties from \cite{Lewis80}.
 
\begin{prop}
Suppose $\underline{M},\underline{N},\underline{L}$ are $\underline{R}$-modules. Then we have the following natural isomorphisms and properties:
\begin{enumerate}
\item $\underline{M}_X \underset{\underline{R}} \boxtimes \underline{R}_Y \cong \underline{R}_Y \underset{\underline{R}} \boxtimes \underline{M}_X  \cong \underline{M}_{X \times Y}$
\item $\underline{\Hom}_{\underline{R}} (\underline{R}_X, \underline{M}_Y) \cong \underline{M}_{X \times Y}$
\item $ \underline{M}_X \underset{\underline{R}} \boxtimes \underline{N}_Y \cong (\underline{M} \underset{\underline{R}} \boxtimes \underline{N})_{X \times Y}$
\item If $\underline{R}, \underline{S},\underline{T}$ are $\mathcal{O}$-Tambara functors and $\underline{B},\underline{C},\underline{D}$ are respectively $\underline{S}-\underline{R}, \underline{R}-\underline{T},\underline{S}-\underline{T}$ bimodules respectively, then $\underline{\Hom}_{\underline{S}-\underline{T}}(\underline{B} \underset{\underline{R}} \boxtimes \underline{C}, \underline{D}) \cong \underline{\Hom}_{\underline{S}-\underline{R}}(\underline{B},\underline{\Hom}_{\underline{T}}(\underline{C},\underline{D})).$
\item $\underline{M} \underset{\underline{R}} \boxtimes \underline{N}$ and $\underline{\Hom}_{\underline{R}}(\underline{M},\underline{N})$ are $\underline{R}$-modules satisfying the adjunction $\underline{\Hom}_{\underline{R}}(\underline{M} \underset{\underline{R}}\boxtimes \underline{N},\underline{L}) \cong  \underline{\Hom}_{\underline{R}}(\underline{M}, \underline{\Hom}_{\underline{R}}(\underline{N},\underline{L}))$, making $\cat{\underline{R}Mod}$ a symmetric monoidal closed category.
\end{enumerate}
\end{prop}

\begin{prop}
The category of Mackey functors is an abelian category satisfying the AB5 condition. $\underline{R} \boxtimes (-)$ and $\underline{\Hom}(\underline{R},-)$ are respectively left and right adjoints for the forgetful functor from $\cat{\underline{R}Mod}$ to $\cat{Mac}.$ Limits and colimits in $\cat{\underline{R}Mod}$ are obtained from taking limits and colimits in $\cat{Mac}$. $\cat{\underline{R}Mod}$ is an abelian category satisfying AB5.
\end{prop}

\subsection{Free Modules}
We now define what it means to be a basis. The following structure is modeled after the algebra notes of \cite{garrett}.
\begin{defn}
A \term{G indexed set} $\underline{S}$ is a collection of sets $\{S_X\}$ for all $X \in \cat{Set}^G$. There is no other necessary structure. A homomorphism $\underline{S} \to \underline{T}$ of $G$ indexed sets is a collection of maps $\{S_X \to T_X \}$. This makes $G$ indexed sets a category $\cat{GIndexedSet}$. There is a forgetful functor $U: \cat{Mac} \to \cat{GIndexedSet}$ which forgets all the transfer and restriction maps. 
\end{defn}

\begin{defn} 
Let $\underline{S}$ be a $G$-indexed set. A \term{free $\underline{R}$-module} $\underline{M}$ on generators $\underline{S}$ is an $\underline{R}$-module $\underline{M}$ and a $G$ indexed set map $i: \underline{S} \to U\underline{M}$ such that, for any $\underline{R}$-module, any $G$-indexed set map $f: \underline{S} \to U \underline{N}$ there is a unique $\underline{R}$-module homomorphism $\tilde{f}: \underline{M} \to \underline{N}$ such that $U\tilde{f} \circ i = f.$ The elements $i(\underline{S})$ are an $\underline{R}$-\term{basis} for $\underline{M}$
\end{defn}

\begin{defn}
Suppose $E = \{ e_i \in \underline{M}(X_i)\}$ is a collection of elements of $\underline{M}$. Then $\underline{E}$ is a \term{basis} for $\underline{M}$ if every $G$ indexed set map from $\underline{E} \to U\underline{N}$, where $\underline{N}$ is another $\underline{R}$-module, gives a unique $\underline{R}$-module map from $\underline{M} \to \underline{N}$. In other words, there is an isomorphism $$\cat{GIndexedSet}(\underline{E},U(\underline{N})) \cong \cat{\underline{R}Mod}(\underline{M},\underline{N}).$$
\end{defn}

\begin{prop}
If a free $\underline{R}$-module  on generators $\underline{S}$ exists, it is unique up to unique isomorphism.
\end{prop}
\begin{proof}
This is identical to the classical proof. By the defining property, the identity is the unique map $F: \underline{M} \to \underline{M}$ that fixes $\underline{S}$. If there is another free module $\underline{M}'$ with $i': \underline{S} \to \underline{M}'$  then there unique maps between $\underline{M}$ and $\underline{M}'$ that commute with the inclusion of $\underline{S}$. These are mutual inverses so $\underline{M}$ is isomorphic to $\underline{M}'$. There is no other map between them that respects the inclusion so this isomorphism is unique.
\end{proof}

\begin{prop}
A free $\underline{R}$-module $\underline{M}$ is generated by $i(\underline{S})$, in that the smallest submodule including $i(\underline{S})$ is $\underline{M}$.
\end{prop}

\begin{proof}
Let $\underline{N}$ be the submodule of $\underline{M}$ that is the intersection of all submodules including $i(\underline{S})$. Consider $f : \underline{S} \to \underline{M}/\underline{N}$, by $f(s) = 0$ for all $s \in \underline{S}$. Both the quotient and the zero map from $\underline{M}$ to $\underline{M}/\underline{N}$ satisfy the property of commuting with $f$, so they must be equal.
\end{proof}

\begin{cor}
Every element of $\underline{M}$ is of the form $\sum_{i \in I} T_{p_i}( r_i \cdot R_{q_i}(m_i)),$ where $m_i$ is the image of the generators of $\underline{S}$, where $I$ goes over the generators $\underline{S}$, and all but finitely many of these terms must be 0.
\end{cor}
\begin{proof}
This is a submodule: using the Frobenius reciprocity relations, one can see it is closed with respect to multiplication by $\underline{R}$, transfers, and restrictions. It is the smallest submodule including all the $m_i$, so we have this corollary. 
\end{proof}

\begin{defn}
Let $m_i \in \underline{M}(X_i)$ be a finite set of elements of $\underline{M}$ indexed by $I$ in the $G$-set $X_i$. If the relation $$\sum_{i \in I} T_{p_i}( r_i \cdot R_{q_i}(m_i)) = 0$$ with $q_i: U_i \to X_i$, $p_i: U_i \to Y$ and $r_i \in \underline{R}(U_i)$ implies $$T_{p_i}(r_i \cdot R_{q_i}(1)) \in \underline{R}(Y)$$ is 0 for all $i \in I$, then we say the elements $m_i$ are \term{linearly independent}.
\end{defn}

This is the best we can say about the uniqueness of the elements above, because $\underline{R}$ is naturally a $\underline{A}$-module, and multiplication by $\underline{R}$ can include the composition of restrictions and transfers.

\begin{prop}
Let $\underline{M}$ be a free $\underline{R}$-module on generators $i: \underline{S} \to U\underline{M}.$ Then the elements $i(S)$ are linearly independent.
\end{prop}
\begin{proof}
Suppose $$\sum_{i \in I} T_{p_i}( r_i \cdot R_{q_i}(s_i)) = 0.$$ To show that every  $T_{p_i}(r_i \cdot R_{q_i}(1)) = 0$, fix $s_k \in \underline{S}(U_k)$ and map $f: \underline{S} \to \underline{R}$ by the map that sends $s_i$ to 1 and every other element to 0. Let $\tilde{f}$ be the associated map $\underline{R}$-module homomorphism $\underline{M} \to \underline{R}$. Then $$\tilde{f}\left(\sum_{i \in I} T_{p_i}( r_i \cdot R_{q_i}(s_i))\right) =T_{p_k}( r_k \cdot R_{q_k}(1)) =  \tilde{f}(0) = 0.$$ This holds for every $k$.
\end{proof}

In the classical case, there is an isomorphism between free modules generated by sets $S$ and $T$ if and only if $|S| = |T|$. That is, the cardinality of the generating set determines the free module up to isomorphism. The equivariant case is similar but a bit more subtle. For example, an $\underline{R}$-module generated by a single element at the $G$-set $X \amalg Y$ is isomorphic to the $\underline{R}$-module generated by two elements, one at $X$ and one at $Y$. Therefore, every free $\underline{R}$-module can be described by the number of generators at each orbit.

\begin{prop}
Let $\underline{M}$ be a free $\underline{R}$-module generated by $i: \underline{S} \to \underline{M}$. Then $\underline{M}$ is also the free module generated by the inclusion of the elements $res_j(i(\underline{S}))$ where if $i(\underline{S}) \in \underline{M}(\coprod_j G/H_j)$ then $res_j$ restricts to one of the orbits $G/H_j$.
\end{prop}

\begin{proof}
Because $\underline{N}(\coprod_j G/H_j) = \prod_j \underline{N}(G/H_j)$, we have an isomorphism between $$\cat{GIndexedSet}(\underline{S}, U\underline{N}) \cong \cat{GIndexedSet}(res(\underline{S}),U\underline{N}).$$ The image of every element of $S$ is exactly the same data of where every restriction maps.
\end{proof}

\subsection{Construction}

\begin{prop}
Let $\underline{Y}[x_H] = \cat{B_G}(G/H, -)$, that is, the Mackey functor represented by $G/H$. Consider the Mackey functor $\underline{F}[x_H] = \underline{R} \boxtimes \underline{Y}[x_H]$. This is the free $\underline{R}$-module generated by the single element at the orbit $G/H$, with the inclusion mapping to $1 \otimes \id$.
\end{prop}

\begin{proof}
This is trivially an $\underline{R}$-module. There is a map $\underline{R} \boxtimes (\underline{R} \boxtimes \underline{Y}[x_H]) \to \underline{R} \boxtimes \underline{Y}[x_H]$ using the multiplication $\underline{R} \boxtimes \underline{R} \to \underline{R}$, and an inclusion $(\underline{R} \boxtimes \underline{Y}[x_H]) \to \underline{R} \boxtimes (\underline{R} \boxtimes \underline{Y}[x_H])$ given by sending an element $T_p(r_x \otimes y_x)$ to $1 \otimes T_p(r_x \otimes y_x).$ 
The $\underline{R}$-module structure explicitly is $$r_y \cdot T_p(r_x \otimes y_x) = T_p( (R_p(r_y) \cdot r_x) \otimes y_x).$$ 


We now prove that $\underline{F}[x_H]$ is the free module described. Let $\underline{M}$ be an $\underline{R}$-module. Then the image of identity span must be an element of $\underline{M}(G/H)$. The element completely determines the map $\underline{F}[x_H] \to \underline{M}$. An element of $\underline{F}$ is of the form $T_p(r_x \otimes R_h)$ by \cite{StricklandTambara} 3.14, and in order for the map to a be a map of $\underline{R}$-modules, if $(1 \otimes T_{\id}R_{\id}) \mapsto m_{G/H} \in \underline{M}(G/H)$, then $T_p(r_x \otimes R_h)$ must map to $T_p( r_x \cdot R_h(m_{G/H})).$ The only remaining piece is to show that for all $m_{G/H}$, the map $$T_p(r_x \otimes R_h) \mapsto T_p( r_x \cdot R_h(m_{G/H}))$$ defines a map of $\underline{R}$-modules. It is clearly a map of Mackey functors. For the multiplication by $\underline{R}$ we see:
$$r_y \cdot T_p(r_x \otimes R_h) = T_p(R_p(r_y) r_x \otimes R_h) \mapsto T_p( R_p(r_y) r_x \cdot R_h(m_{G/H}))= r_y \cdot T_p(  r_x \cdot R_h(m_{G/H}))$$ showing the map respects the $\underline{R}$-module structure.
\end{proof}

\begin{prop}
Suppose $\underline{F_i}$ is a family of free $\underline{R}$-modules generated by $\underline{S_i}$. Then the direct sum $\oplus_i \underline{F_i}$ is the free $\underline{R}$-module generated by $\coprod_i \underline{S_i}$.
\end{prop}

\begin{proof}
This can be proven directly from the universal properties: 
\begin{eqnarray*}
\cat{\underline{R}mod}(\oplus_i \underline{F_i}, \underline{M}) &\cong& \prod_i (\cat{\underline{R}mod}(\underline{F_i},\underline{M})) \\
&\cong& \prod_i (\cat{GIndexedSet}(\underline{S_i},\underline{M})) \\
&\cong& \cat{GIndexedSet}(\coprod_i\underline{S_i},\underline{M}).
\end{eqnarray*}
\end{proof}

The two propositions above combine to make a free $\underline{R}$-module for any set of generators. 

\subsection{Examples}
Here we present as Lewis diagrams the most basic free modules that do not appear in the classical case. Let $G = C_2$, let $\mathcal{O}_c$ be the complete indexing system, and $\mathcal{O}_{triv}$ be the trivial indexing system. The Burnside ring $\underline{A}^{\mathcal{O}_c}$ for $C_2$ is
\[\xymatrixrowsep{.2in}
\xymatrixcolsep{.1in}\xymatrix{
G/G: &\Z[t]/(t^2-2t)  \ar@/^-3pc/[d]^{R^G_e} \\
G/e: &\Z \ar@/^-5pc/[u]^{T^G_e} \ar@/^-3pc/[u]^{N^G_e}
}\]
We have $$R_e^G(t) = 2, T_e^G(1) = t, N_e^G(a) = a^2$$ and the Burnside ring $\underline{A}^{\mathcal{O}_{triv}}$ for $C_2$ is
\[\xymatrixrowsep{.2in}
\xymatrixcolsep{.1in}\xymatrix{
G/G: &\Z[t]/(t^2-2t)  \ar@/^-3pc/[d]^{R^G_e} \\
G/e: &\Z \ar@/^-5pc/[u]^{T^G_e}
}\] with $$R_e^G(t) = 2, T_e^G(1) = t.$$ We only state some parts of the maps $R_e^G,T_e^G,N_e^G$ since the rest of the map is determined.
In both cases, there are two identical free modules, 
\[\xymatrixrowsep{.2in}
\xymatrixcolsep{.1in}\xymatrix{
&\underline{A}^{\mathit{O_c}}\{x_G\} \cong \underline{A}^{\mathit{O_{triv}}}\{x_G\}  \cong \cat{Mac}(G/G,-)\\
G/G: &\Z[t]/(t^2-2t)\{x\}  \ar@/^-3pc/[d]^{R^G_e}\\
G/e: &\Z \{r\} \ar@/^-5pc/[u]_{T^G_e} \\
&\underline{A}^{\mathit{O_c}}\{x_e\} \cong \underline{A}^{\mathit{O_{triv}}}\{x_e\}  \cong \cat{Mac}(G/e,-)\\
G/G: &\Z\{t_x \}  \ar@/^-3pc/[d]^{R^G_e}\\
G/e: &\Z \{x, \overline{x}\}\ar@/^-5pc/[u]_{T^G_e}\\
}\]
where $\overline{x}$ means the Weyl conjugate. $T_e^G(x) = T_e^G(\overline{x})= t_x$ and $R_e^G(t_x) = x + \overline{x}.$

For some more interesting examples, we will consider modules over some of the possible singly generated free Tambara functors over the Burnside ring:

\[\xymatrixrowsep{.2in}
\xymatrixcolsep{.1in}\xymatrix{
&\underline{A}^{\mathit{O_c}}[x_G]\\
&\Z[t]/(t^2-2t)[x,n]/t(n-x^2)  \ar@/^-3pc/[d]^{R^G_e}\\ 
&\Z[r] \ar@/^-5pc/[u]_{T^G_e} \ar@/^-3pc/[u]^{N_e^G}\\
&\\
&\underline{A}^{\mathit{O_{triv}}}[x_G]\\
&\Z[t]/(t^2-2t)[x]  \ar@/^-3pc/[d]^{R^G_e}\\
&\Z[r] \ar@/^-5pc/[u]_{T^G_e} 
}\]

Here, $R^G_e(x)=r$, $T^G_e(\cdot)$ is multiplication by $t$, $N^G_e(r)=n.$   Each one of these Tambara functors has two free modules of rank 1.

We now write down some of the free modules over these Tambara functors:

\[\xymatrixrowsep{.2in}
\xymatrixcolsep{.1in}\xymatrix{
&{\underline{A}^{\mathit{O_c}}[x_G]}\{y_G\}\\
G/G: &\Z[t]/(t^2-2t)[x,n]/t(n-x^2)\{y\}  \ar@/^-3pc/[d]^{R^G_e}\\
G/e: &\Z[r_x] \{r_y\} \ar@/^-5pc/[u]_{T^G_e} \\
\\
&{\underline{A}^{\mathit{O_{triv}}}[x_G]}\{y_G\}\\
G/G: &\Z[t]/(t^2-2t)[x]\{y\}  \ar@/^-3pc/[d]^{R^G_e}\\
G/e: &\Z[r_x] \{r_y\}\ar@/^-5pc/[u]_{T^G_e}\\ 
}\]
The restriction in these two cases is clear given the underlying Tambara functor. The transfer of $p(x) \cdot r_y$ is $T_e^G(p(x)) \cdot y$ which follows from the module condition.
Note that we only give these examples because when we take the generators at $G/e$, the Weyl action on the generator complicates the formulas. The reader is invited to see the computations in \cite{blumberg2019right}.

\subsection{Projective Modules}

We should note a few warnings before we continue. In the classical case, $\Z[X]$, the free ring generated by a set $X$, when we forget the multiplication and consider $\Z[X]$ as an abelian group, it is a free $\Z$-module with countably many generators $1,x,x^2, \ldots$ This is not true in the equivariant setting as shown by \cite{incomplete}. The free Tambara functor $\underline{A}^{\mathcal{O}}[x_i]$, when we forget the norms to make it an $\underline{A}$-module, or a Mackey functor, is not typically a free $\underline{A}$-module.

Secondly, \cite{lewis1999projective} shows the box product of projective objects is not necessarily projective in some module categories of equivariant stable homotopy theory, but they are for $\cat{\underline{R}Mod}$, which we address. The discussion is, again, very similar to the classical case.

\begin{prop}
Free $\underline{R}$-modules are projective: Suppose $\underline{M}$ is a free $\underline{R}$-module on generators $\underline{S}$, a surjective map of $\underline{R}$-modules from $\underline{N_0} \to \underline{N_1}$ and a map $\underline{M} \to \underline{N_1}$. Then there exists a lifting $\underline{P} \to \underline{M_0}$ that commutes with these maps.
\end{prop}

\begin{proof}
The proof is trivial. The inclusion $i: \underline{S} \to \underline{M}$ gives a map $\underline{S} \to \underline{N_1}$ by composition. Since $\underline{N_0} \to \underline{N_1}$ is surjective, there is a map $\underline{S} \to \underline{N_0}$ such that the composition with $\underline{N_0} \to \underline{N_1}$ gives the composition $\underline{S} \to \underline{M} \to \underline{N_1}$. This map induces a map of $\underline{R}$-modules $\underline{M} \to \underline{N_0}$. Because the diagram commutes on the inclusion of $\underline{S} \to \underline{M}$, the diagram commutes by the universal property of being free.
\end{proof}

\begin{prop}
The following are equivalent for an $\underline{R}$-module $\underline{P}$:
\begin{enumerate}
\item $\underline{P}$ is projective.
\item Any short exact sequence $0 \to \underline{M} \to \underline{N} \to \underline{P} \to 0$ splits.
\item $\underline{P}$ is a direct summand of a free $\underline{R}$-module.
\end{enumerate}
\end{prop}
\begin{proof}
Same proof as you might find in \cite{jacobson2012basic}. For (1) implies (2), $\underline{P}$ being projective gives a lift $\underline{P} \to \underline{N}$. For (2) implies (3), any module is the image of a free $\underline{R}$-module, so then we have $0 \to \underline{P'} \to \underline{F} \to \underline{P} \to 0$. This splits so $\underline{P}$ is a summand of $\underline{F}$. For (3) implies (1), we have a map into a free module $\underline{F}$, which is projective. Composition with the lifting map from $\underline{F}$ makes $\underline{P}$ projective.
\end{proof}

\begin{prop}
If $\underline{M},\underline{N}$ are projective $\underline{R}$-modules, then $\underline{M} \underset{\underline{R}} \boxtimes \underline{N}$ is projective.
\end{prop}

\begin{proof}
Suppose $\underline{M}, \underline{N}$ are summands of free modules $\underline{S},\underline{T}$ respectively. Then $\underline{M} \underset{\underline{R}}\boxtimes \underline{N}$ is a direct summand of $\underline{S} \underset{\underline{R}} \boxtimes \underline{T}$, a free module.
\end{proof}

%


\subsection{Flat Modules}

We will address the alarming discussion of \cite{lewis1999projective}, and that in many contexts in equivariant stable homotopy theory, basic properties of projective modules do not hold. In particular, projective modules are not always flat and box products of projective modules are not always projective.

We explicitly stated our assumption that $G$ is a finite group. We also implicitly state that with regard to the underlying Mackey functor, the $G$-universe is complete (for every projection $G/H \to G/K$, there is a corresponding transfer map). With these assumptions, the category of Mackey functors, and similarly $\cat{\underline{R}Mod},$ has the property that projective modules are flat and products of projective modules are projective.
The same is not true for when $G$ is compact Lie or when the $G$-universe is incomplete. In this context, we can still define the Andr\'e-Quillen (co)homology Mackey functors. But there is no reason to believe that many of properties that rely on the above conditions to be true. In particular, the Jacobi-Zariski sequence or the convergence of the fundamental spectral sequence need not be true.

\begin{prop}
Projective $\underline{R}$-modules are flat.
\end{prop}

\begin{proof}
\cite{Lewis80} shows that $\underline{A}_{G/H}$ is flat for all subgroups $H$ of $G$. Therefore, all free $\underline{R}$-modules are flat. $(\underline{R} \boxtimes \underline{A}_{G/H}) \underset{\underline{R}} \boxtimes (-) \cong  \underline{A}_{G/H} \boxtimes (-),$ so the right being exact means the left is exact. Direct sums of flat modules are flat, so all free $\underline{R}$-modules are flat. Direct summands of flat modules must also be flat, proving the proposition.
\end{proof}

%% file: ExtTor.tex
\section{Ext and Tor Functors}

As shown earlier, for $\underline{A}$ an $\underline{R}$-module, $\underline{B}$ an $(\underline{R},\underline{S})$-bimodule, $\underline{C}$ an $\underline{S}$-module, then $$\cat{\underline{S}Mod}(\underline{A} \underset{\underline{R}} \boxtimes \underline{B}, \underline{C} ) \cong \cat{\underline{R}Mod}(\underline{A}, \underline{\Hom}_{\underline{S}}(\underline{B},\underline{C})). $$
Being adjoints, $(-) \underset{\underline{R}} \boxtimes \underline{B}$ is right exact and $\underline{\Hom}_{\underline{S}}(\underline{B},-)$ is left exact. These give rise to derived functors $$\underline{\Tor}_n^{\underline{R}}(\underline{A},\underline{B}):= L_n (- \underset{R}\boxtimes B)(A)$$ and $$\underline{\Ext}^n_{\underline{S}}(\underline{B},\underline{C}) := R^n (\underline{\Hom}_{\underline{S}}(\underline{B},-)) (\underline{C}).$$

%% file: AQ.tex
\chapter{Andr\'e-Quillen (co)homology}
\section{Square zero extension and augmentation ideal functors are inverse equivalences}
The following is a equivariant retelling of \cite{QuillenMimeo}. Andr\'e-Quillen cohomology is a cohomology theory for augmented commutative rings, that is maps of commutative rings $A \to B$. Quillen introduced a model structure on the category of simplicial commutative rings. He then defined the cotangent complex as the left derived functor of abelianization, and the (co)homology as the (co)homology of this cotangent complex. The cotangent complex can be computed explicitly, via free resolutions in the category of simplicial commutative algebras. The Andr\'e-Quillen cohomology and cotangent complex is closely related to the deformation theory of commutative rings.

Let $G$ be a finite group and $\mathcal{O}$ be an indexing system.
Let $\underline{A}$ be an $\mathcal{O}$-Tambara functor and let $\underline{A}\cat{-CAlg} := \underline{A} / (\mathcal{O}\cat{-Tamb})$ be the under category. This is the natural generalization of a commutative algebra of a ring. A morphism is a map of $\mathcal{O}$-Tambara functors that preserves the inclusion of the base $\mathcal{O}$-Tambara functors $\underline{A}$.

In the non-equivariant setting, the category of abelian objects of the over category $A \cat{-CAlg}/B$ is equivalent to the category of $B$-modules. Let us go through this well known fact to highlight the difference in the equivariant case. There is an augmentation ideal functor from $(A \cat{-Alg}/B)_{ab}$ to $B \cat{-Mod}$ that assigns to $R \to B$ the kernel of that map. Being the kernel of a map, it is naturally an ideal. $R$, being an abelian group object, has a map from $B \to R$, making $R$ a $B$-module. Because all of these are maps of rings, this extends to a $B$-module structure on the kernel from $R$ to $B$.

The inverse is the square-zero extension $B \ltimes (-) : B\cat{-Mod} \to (A \cat{-CAlg}/B)_{ab}.$ As an abelian group, $B \ltimes M$ is isomorphic to $B \oplus M$. By the name, we assert that the products on the second term are zero, so $M^2 = 0$ inside $B \ltimes M$. Additionally, the product of an element of $B$ and an element of $R$ must reflect the $B$-module structure. Together with distributivity, this implies that $(b,m)\cdot (b',m') = (bb',bm' + b'm).$ There is the inclusion map $B \to B \ltimes M$ that serves as the identity and the map $$(B \ltimes M) \underset{B}{\times} (B \ltimes M) \to (B \ltimes M)$$ given by $$((b,m), (b,m')) \mapsto (b,m+m')$$ is a map of rings. The inverse map is given by 
$(b,m) \mapsto (b,-m).$ These all satisfy the associativity, commutativity, left-right identity, left-right inverse properties, so this indeed is an abelian group object.

We now prove that these are inverse equivalences. One direction is clear: given a $B$-module, the augmentation ideal of the associated abelian group object is isomorphic to the original module. In the reverse direction, suppose we have $R \to B$ an abelian group object of $(A \cat{-CAlg}/B).$ We want to show that $B \ltimes \ker(R \to B)$ is isomorphic to $R$ as elements of $(A \cat{-CAlg}/B)$. The map from left to right is addition $(b,m) \mapsto b+m$ where the map from $B$ to $R$ is suppressed in this notation. $R$ has a map $\sigma: R \underset{B}\times R \to R$. If $m_0,m_1 \in \ker(R \to B)$, let $n_0 = (m_0,0)$ and $n_1 = (0,m_1) \in R \underset{B}\times R$. We see that $\sigma(n_i) = m_i$ and because $\sigma$ is a map of rings, $m_0 m_1 = 0$. So the map $B \ltimes \ker(R \to B) \to R$ is a map of rings. The inverse map is $r \mapsto (\epsilon(r), r- \epsilon(r))$ where $\epsilon$ is the map from $R \to B$. So the above is isomorphic as commutative rings and it is easy to see that these maps preserve the maps from $A$ and to $B$ showing that these are isomorphic as elements in $(A \cat{-CAlg}/B)$. This was relayed in \cite{StricklandTambara}.

Let us now move into the equivariant setting, as in \cite{StricklandTambara}. For the exact same reason as in the non-equivariant setting, taking the augmentation ideal of an abelian group object of $\underline{A}\cat{-CAlg}/ \underline{B}$ is a $\underline{B}$-module. To generalize the square-zero extension, we let $\underline{B} \ltimes \underline{R} = \underline{B} \oplus \underline{R}$ as Mackey functors. The only additional element to describe is the norm maps. We assert that for any $\mathcal{O}$-admissible $G$ map $f:X\to Y$ such that $f^{-1}(y)$ has at least 2 elements (which we then call $f$ a 2-surjective map), $N_f((0,m)) = 0$ for all $m \in \underline{M}(X)$. This is the natural generalization of the square-zero extension in the following senses: it both includes the standard squares given by $f:X \coprod X \to X$ and fully determines the $\mathcal{O}$-Tambara structure of $\underline{B} \ltimes \underline{M}$: $N_f(b,m) = N_f((b,0)+(0,m))$, which if we rewrite in TNR-form we produce the following exponential diagram:
\[\xymatrixrowsep{.5in}
\xymatrixcolsep{1in}\xymatrix{
X \ar[d]^f& X \amalg X \ar[l]_{\nabla} & X \underset{Y}\times \Pi_f (X \amalg X) \ar[l]_{ev} \ar[d]^{\pi_2}\\
Y&& \Pi_f (X \amalg X) \ar[ll]^{q}
}\]
Now using the fact that $\underline{M}$ is killed by 2-surjective maps, we reduce to the case that 
\[\xymatrixrowsep{.5in}
\xymatrixcolsep{1in}\xymatrix{
X \ar[d]^f& X \amalg X \ar[l]_{\nabla} & X \amalg (X \underset{Y} \times X) \ar[l]_{ev} \ar[d]^{\pi_2}\\
Y&& Y \amalg (Y \underset{Y} \times X) \ar[ll]^{q}
}\]
where the first part of the coproduct represents the maps $X \to X \amalg X$ that map completely to the left side, and the second part of the coproduct represents the maps $X \to X \amalg X$ that map a single element to the right side. The term that picks up the $m$ in $X \underset{Y} \times X$ is the diagonal $\Delta(X)$. Let $f: X \to Y$ be an admissible map of $G$-sets and let $\pi_1,\pi_2 : X \underset{Y}\times X - \Delta(X) \to X$ be the projections. Thus we showed $$N_f(b,m) = (N_f(b), T_f(N_{\pi_1}R_{\pi_2}(b) \cdot m)).$$ \cite{StricklandTambara} 14.9 shows that this is an element of $(\underline{A}\cat{-CAlg}/\underline{B})_{ab}$ for the complete case, and the incomplete case is identical.

Any $\underline{B}$-module has a square-zero extension in $\underline{A}\cat{-CAlg}/\underline{B}$ and by the same argument as in the non-equivariant setting, it is an abelian group object. The issue is that being an abelian group object is insufficient to be isomorphic to the square-zero extension of an augmentation ideal. An extra condition is needed, as shown in the next proposition.

\begin{prop}
Let $\cat{D}$ be the full subcategory of $(\underline{A}\cat{-CAlg}/\underline{B})_{ab}$ such that the kernel of the map to $\underline{B}$ has vanishing norms for all 2-surjective admissible norms. Then $\cat{D}$ is equivalent to the category of $\underline{B}$-modules. The two maps are the equivariant square-zero extension and the augmentation ideal.
\end{prop}

Note that products are still trivial in the augmentation ideal for any element of $(\underline{A}\cat{-CAlg}/\underline{B})_{ab}$. So in the case of the trivial indexing system, we have the classical equivalence of abelian group objects and modules.

\begin{proof}
The augmentation ideal of a square-zero extension gives a module isomorphic to the original module. In the reverse direction, suppose we have $\underline{R} \to \underline{B}$ an abelian group object of $(\underline{A} \cat{-CAlg}/\underline{B})$ where the augmentation ideal has vanishing norms. We want to show that $\underline{B} \ltimes \ker(\underline{R} \to \underline{B})$ is isomorphic to 
$\underline{R}$ as elements of $(\underline{A} \cat{-CAlg}/\underline{B})$.  The map from left to right is addition $(b,m) \mapsto b+m$ at each $G$-set. This is evidently a map of Mackey functors. We now need to show that the map preserves norms. It is sufficient to prove that it preserves norms of the form $f: G/K \to G/H$ where $f$ is admissible, as the products case reduces to the classical statement.
We know that $N_f(b) = (N_f(b),0)$ and $N_f((0,m))= 0.$ We can rewrite $N_f(b+m)$ in TNR-form using the exponential diagram:
\[\xymatrixrowsep{.5in}
\xymatrixcolsep{1in}\xymatrix{
G/K \ar[d]^f& G/K \amalg G/K \ar[l]_{\nabla} & G/K \underset{G/H}\times \Pi_f (G/K \amalg G/K) \ar[l]_{ev} \ar[d]^{\pi_2}\\
G/H&& \Pi_f (G/K \amalg G/K) \ar[ll]^{q}
}\]
Again, we can restrict to the case where the sections have at most one element in the right set of the coproduct. So we get the following diagram:
\[\xymatrixrowsep{.5in}
\xymatrixcolsep{.5in}\xymatrix{
G/K \ar[d]^f& G/K \amalg G/K \ar[l]_{\nabla} & G/K \amalg [(G/K \underset{G/H}\times G/K - \Delta(G/K)) \amalg G/K] \ar[l]_(.65){ev} \ar[d]^{f \amalg (\pi_2 \coprod \id)}\\
G/H&& G/H \amalg G/K \ar[ll]^{q}
}\]
where $ev$ on the first term maps by the identity to the left term, and maps $(G/K \underset{G/H}\times G/K - \Delta(G/K))$ by $\pi_1$ to the left element of the coproduct and $G/K$ maps by the identity to the right most element. One can see the correspondences between these sets and the sets involving sections: the $G/K$ and $G/H$ in the left of the coproduct correspond to sections all mapping to the left most set of $G/K \amalg G/K$ and the $(G/K \underset{G/H}\times G/K - \Delta(G/K))\amalg G/K$ are the sections with exactly one element mapping to the $G/K$ on the right. $(G/K \underset{G/H}\times G/K - \Delta(G/K))$ corresponds to section-element pairs where the element does not map by the section to the right set, and the remaining $G/K$ corresponds to when the element does map to the section. We realized the equation $N_f(b+m) = N_f(b)+T_f(N_{\pi_1}R_{\pi_2}(b) \cdot m)$ from the above diagram. Therefore, the map $\underline{B} \ltimes \ker(\underline{R} \to \underline{B})$ is a map of $\mathcal{O}$-Tambara functors. The map also respects the maps from $\underline{A}$ and to $\underline{B}$, so this is a map in $\underline{A} \cat{-CAlg}/\underline{B}$. The inverse map is $r \mapsto (\epsilon(r), r - \epsilon(r))$ where $\epsilon: \underline{R} \to \underline{B}$ at every $G$-set $X$. This is readily seen to be a map of Mackey functors (as it is a Tambara functor or difference of Tambara functor maps on each summand). This is sufficient to show it is a $\mathcal{O}$-Tambara functor map, as it is a Mackey functor isomorphism where the inverse respects norms. We have shown that the two maps are essentially surjective. The maps being fully faithful are clear.

\end{proof}

\begin{prop}
In general, it is not true that $\cat{D}$ is equivalent to $(\underline{A}\cat{-CAlg}/\underline{B})_{ab}$. The following example was given in \cite{StricklandTambara}. Let $G = C_2$ and let $\underline{A} = \underline{B} = \underline{S}$ be the fixed point Tambara functor on $\Z$, and $\underline{T}$ be the following:
\[\xymatrixrowsep{.5in}
\xymatrixcolsep{.6in}\xymatrix{
&\underline{S}&&& \underline{T}\\
G/G: &\Z \ar@/^-1pc/[d]^{R^G_e}&&  &\Z[\alpha]/\alpha^2 \ar@/^-1pc/[d]^{R^G_e} \\
G/e: &\Z \ar@/^-1pc/[u]_{T^G_e} \ar@/^-2.5pc/[u]_{N^G_e} && & \Z[\beta,\gamma]/ (\beta^2,\beta\gamma,\gamma^2,2\gamma)   \ar@/^-2.5pc/[u]_{N^G_e}  \ar@/^-1pc/[u]_{T^G_e}\\
&\\
}\]
where $\underline{T}(G/G)$ has trivial $C_2$ action with the following on maps:
\begin{eqnarray*}
R^G_e (i+j \beta + k\gamma) &=& i+2 j \alpha\\
T^G_e (i + j \alpha) &=& 2i + j \beta \\
N^G_e (i + j \alpha) &=& i^2 +ij \beta + j^2 \gamma.
\end{eqnarray*}
There is an inclusion map $\underline{S} \to \underline{T}$ and an augmentation map $\underline{T} \to \underline{S}$ by just taking the constant terms. Then $\underline{T} \underset{\underline{S}}\times \underline{T}$ is 
\[\xymatrixrowsep{.5in}
\xymatrixcolsep{.1in}\xymatrix{
\underline{T} \underset{\underline{S}}\times \underline{T}\\
\Z \{ 1, \alpha_0, \alpha_1 \} / \sim \ar@/^-1pc/[d]^{R^G_e} \\
 \Z\{1, \beta_0, \beta_1\} \oplus (\Z/2)\{ \gamma_0,\gamma_1\} /\sim \ar@/^-2.5pc/[u]_{N^G_e}  \ar@/^-1pc/[u]_{T^G_e}\\
&\\
}\]
and there is a map of Tambara functors from $\underline{T} \underset{\underline{S}}\times \underline{T} \to \underline{T}$ which adds all the non-constant terms, just as in the non-equivariant setting. It is easy to check that this map is commutative, associative, and unital, so $\underline{T}$ is an abelian group object in $\underline{A}\cat{-CAlg}/\underline{B}.$ But the norms on the augmentation ideal are non-vanishing: the norm of $\alpha$ in $\underline{T}$ is $\gamma$. Therefore, this abelian group object is not isomorphic as elements of $\underline{A}\cat{-CAlg}/\underline{B}$ to the square-zero extension.
\end{prop}

In comparison to the above example, when the indexing system is trivial, the abelian objects is $\cat{D}$. In general $\cat{D}$ is the $\mathcal{O}$-commutative monoid objects. Much of this material is from \cite{HillAndre} and \cite{hill2016equivariant} in the complete setting. The reader is invited to read \cite{hill2016equivariant} for the definition and context of $\mathcal{O}$-commutative monoid objects.

Firstly, $\mathcal{O}$-Tambara functors form a symmetric monoidal coefficient system. The direct sum, the categorical product, of two $\mathcal{O}$-Tambara functors is an $\mathcal{O}$-Tambara functor. As described in \cite{incomplete} Section 6, induction gives a map between polynomials with restricted exponents
$$\uparrow^G_H: \mathit{P}_{i^*_H \mathcal{O}}^H \hookrightarrow \mathit{P}_{\mathcal{O}}^G$$ which gives a restriction map on Tambara functors $$i_H^* : \mathcal{O}-\cat{Tamb}_G \to (i_H^* \mathcal{O})-\cat{Tamb}_H.$$ Note that this restriction map is functorial in the sense that $i_K^* i_H^* \cong i_K^*$, and 
the restriction map commutes with symmetric monoidal product, so this forms a symmetric monoidal coefficient system. Given an $H$-set, we define the internal product $T \square \underline{M} := \underline{M}(T \times (-)),$ which by Corollary 6.7 of \cite{incomplete}, is a $\mathcal{O}$-Tambara functor. This gives $\mathcal{O}$-Tambara functors the structure of a $G$-symmetric monoidal coefficient system.

Now let us consider the category of $\underline{A}-\cat{CAlg}  = \underline{A} / (\mathcal{O}\cat{-Tamb})$. The direct sum still makes $\underline{A}-\cat{CAlg}$ a symmetric monoidal category. This is additionally true for the category $i^*_H \underline{A} / (i_H^* \mathcal{O}\cat{-Tamb})$ and $i_H^*$ respects the symmetric monoidal structure, making $G/H \mapsto i^*_H \underline{A} / (i_H^* \mathcal{O}\cat{-Tamb})$ a symmetric monoidal coefficient system. There is a natural map of $\mathcal{O}$-Tambara functors $\underline{R} \to \underline{R}_T := \underline{R}(T \times (-))$. Therefore, setting $T \square \underline{M} := \underline{M}_T$ makes $\underline{A}-\cat{CAlg}$ a $\mathcal{O}$-symmetric monoidal category.

Finally, let us consider the category at hand, $\underline{A}\cat{-CAlg}/\underline{B}$. Direct sum does not make this category a symmetric monoidal category: given $\underline{R} \to \underline{B}$ and $\underline{R'} \to \underline{B}$, then there is not a natural map of $\mathcal{O}$-Tambara functors $\underline{R} \oplus \underline{R'} \to \underline{B}.$ However the pullback
\[\xymatrixrowsep{.5in}
\xymatrixcolsep{.5in}\xymatrix{
\underline{R} \underset{\underline{B}}\oplus \underline{R'} \ar[r] \ar[d]& \underline{R} \ar[d] \\
\underline{R'} \ar[r]& \underline{B}
}\]
does make $\underline{A}\cat{-CAlg}/\underline{B}$ a symmetric monoidal category, with $\underline{B}$ being the unit. This also makes $i_H^* \underline{A}\cat{-CAlg}/ i_H^* \underline{B}$ a symmetric monoidal category and the restriction respects addition: $i_H^* (\underline{B} \underset{\underline{R}}\oplus \underline{B'}) = i_H^* \underline{B} \underset{i_H^* \underline{R}}\oplus i_H^* \underline{B'}.$ We see that $G/H \mapsto i_H^* \underline{A}\cat{-CAlg}/ i_H^* \underline{B}$ is a symmetric monoidal coefficient system. Note that one way of describing the exponentiation $\underline{R} \mapsto \underline{R}  \underset{\underline{B}}\oplus \underline{R}$ can be expressed in the following pullback 
\[\xymatrixrowsep{.5in}
\xymatrixcolsep{.5in}\xymatrix{
\underline{R} \underset{\underline{B}}\oplus \underline{R} \ar[r] \ar[d]& \underline{R} \oplus \underline{R} \ar[d] \\
\underline{B} \ar[r]& \underline{B} \oplus \underline{B}
}\]
so we can extend this to a $G$-symmetric monoidal coefficient system (even though $\underline{A}$ and $\underline{B}$ are only $\mathcal{O}$-Tambara functors, not necessarily complete Tambara functors) by defining $T \square \underline{R}$ to be the pullback
\[\xymatrixrowsep{.5in}
\xymatrixcolsep{.5in}\xymatrix{
T \square \underline{R} \ar[r] \ar[d]& \underline{R}(T \times (-)) \ar[d] \\
\underline{B} \ar[r]&  \underline{B}(T \times (-)) 
}\]
This definition generalizes the standard exponentiation and because 
\[\xymatrixrowsep{.5in}
\xymatrixcolsep{.5in}\xymatrix{
T \square \underline{R} (T' \times (-))\ar[r] \ar[d]& \underline{R}(T \times (T' \times-)) \ar[d] \\
\underline{B} (T' \times (-)) \ar[r]&  \underline{B}(T \times (T' \times -)) 
}\]
is also a pullback diagram then we have the diagram
\[\xymatrixrowsep{.5in}
\xymatrixcolsep{.5in}\xymatrix{
T' \square (T \square \underline{R}) \ar[r] \ar[d]& (T \square \underline{R})(T' \times (-)) \ar[d] \ar[r] &  \underline{R}(T \times (T' \times-)) \ar[d]  \\
\underline{B} \ar[r]&  \underline{B}(T' \times (-)) \ar[r]& \underline{B}(T \times (T' \times -)) 
}\]
where both boxes are pullbacks, implying that $T' \square (T \square \underline{R}) \cong (T' \times T) \square \underline{R}$. Now that we have a $G$-symmetric monoidal structure, we can ask what the $\mathcal{O}$-commutative monoids are.

Suppose $\underline{R}$ is a commutative monoid, which means that as a Mackey functor it is of the form $\underline{B} \oplus I(\underline{R})$. Now suppose all $\mathcal{O}$-admissible 2-surjective norms vanish on $I(\underline{R})$. Let $T \in \mathcal{O}(G/G)$. We have a map of Mackey functors given by the transfer: $\underline{B} \oplus I(\underline{R})(T \times (-)) \to \underline{B} \oplus I(\underline{R})$. Because $T \in \mathcal{O}(G/G)$, if $X \to Y$ is a map in $\mathcal{O}(G/G)$, then $T \times X \to T \times Y$ is a map in $\mathcal{O}(G/G)$. So because all $\mathcal{O}$-norms vanish on $I(\underline{R})$, the same is true for $I(\underline{R})(T \times (-)).$ Lastly, $I(\underline{R})(T \times (-))$ is a $\underline{R}$-module via the restriction map $\underline{R} \to \underline{R}(T \times (-))$. By \cite{StricklandTambara} 14.2, $I(\underline{R})(T \times (-)) \to I(\underline{R})$ is a map of $\underline{R}$-modules. From the above discussion, $\underline{B} \oplus I(\underline{R})(T \times (-)) \to \underline{B} \oplus I(\underline{R})$ is in fact a map of $\mathcal{O}$-Tambara functors, showing that $\underline{R}$ is in fact an $\mathcal{O}$-commutative monoid.

Now let $\underline{R}$ be an $\mathcal{O}$-commutative monoid, so $\underline{R}$ is necessarily an abelian group object, implying that as a Mackey functor $\underline{R} \cong \underline{B} \oplus I(\underline{R})$ where $I(\underline{R})$ is the kernel of the map $\underline{R} \to \underline{B}.$ We would like to show that all $\mathcal{O}$-admissible 2-surjective norms from $\mathcal{O}$ vanish on $I(\underline{R}).$ By the standard argument in commutative rings, products are trivial, so we only need to prove it for norms of the form $G/H \to G/K.$ Note that this norm is given by the norm of $i_K^* I(\underline{R})$ from $K/H \to K/K$. So it suffices to show that the norm $G/H \to G/G$ vanishes.

In this case we see that $$G/H \square \underline{R} \cong \underline{B} \oplus I(\underline{R})(G/H \times (-)).$$ By the definition, we have a map $$\underline{B} \oplus I(\underline{R})(G/H \times (-)) \to \underline{B} \oplus I(\underline{R})$$ of $\mathcal{O}$-Tambara functors. Note that because $\mathcal{O}$ is an indexing system and is closed under products, $\underline{B} \oplus I(\underline{R})(G/H \times (-))$ is also an $\mathcal{O}$-commutative monoid. In particular, products vanish on $I(\underline{R})(G/H \times (-))$ and the norm from $G/H \to G/G$ is isomorphic to a $|G|/|H|$ fold product, up to an action by $G$. So the norm vanishes on $I(\underline{R})(G/H \times (-)).$ We need only show that the map is surjective on the augmentation ideals. Then the norm on $I(\underline{R})$ will vanish because the map is a map of $\mathcal{O}$-Tambara functors.

Inspecting the map $I(\underline{R})(G/H \times (-)) \to I(\underline{R})(-)$, we see that this implies that $I(\underline{R})$ is a $\mathcal{O}$-commutative monoid in Mackey functors. This map is unique, which is the content of \cite{HillAndre} 3.20. We are interested in what the map $$I(\underline{R})(G \underset{H}\times (-)) \to I(\underline{R})$$ is. The point is for Mackey functors (not coefficient systems), $G \underset{H}\times (-)$ gives both left and right adjoints to the restriction map. Rewriting this, we are interested in the map $$\CoInd_H^G i_H^* I(\underline{R}) \overset{f^\sharp} \to \underline{R}$$ which has an adjoint $$i_H^* I(\underline{R}) \overset{f^\flat}\to i_H^* I(\underline{R})$$ which is equal to $$i_H^* I(\underline{R}) \overset{\eta_{i_H^*{I(\underline{R})}}} \longrightarrow i_H^* \CoInd_H^G i_H^* I(\underline{R}) \overset{i_H^* f^\sharp} \longrightarrow i_H^* I(\underline{R}).$$ First we note that $i_H^* I(\underline{R})$ is itself an abelian group object, and $i_H^* f^\sharp$ must be the abelian group structure (so a sum of the elements in $I(\underline{R})$), so 
\begin{eqnarray*}
i_H^* f^\sharp: i_H^* \CoInd_H^G i_H^* I(\underline{R}) &\cong& I(\underline{R}) (G \underset{H} \times (G \underset{H} \times (-)))\\ &\cong&  I(\underline{R}) (i_H^*(G/H) \times (G \underset{H} \times (-))) \to  I(\underline{R})(G \underset{H} \times(-)).
\end{eqnarray*}
 The unit of the adjunction comes from the inclusion of $H/H \to i_H^*(G/H)$ and is the inclusion of $I(\underline{R})(G \underset{H} (-)) \to I(\underline{R}) (i_H^*(G/H) \times (G \underset{H} \times (-))) .$ The composition is the identity. The adjoint of the identity map is the transfer. So the map $I(\underline{R})(G/H \times G/H) \to I(\underline{R})(G/H)$ is surjective. Thus we have proved the essential surjectivity of the following statement. The fully-faithfulness is trivial.

\begin{prop}
Let $G$ be a finite group, $\mathcal{O}$ be an indexing system and $\underline{A}$ and $\underline{B}$ be $\mathcal{O}$-Tambara functors. Then the full subcategory of abelian objects of $\underline{A}\cat{-CAlg}/\underline{B}$ with all 2-surjective $\mathcal{O}$-admissible norms vanishing on the augmentation ideal is equal to the full subcategory of $\mathcal{O}$-commutative monoid objects of $\underline{A}\cat{-CAlg}/\underline{B}.$
\end{prop}
An immediate corollary is
\begin{corollary}
There is an equivalence of categories between $\underline{B}$-modules and $\\ \mathcal{O}$-commutative monoids of $\underline{A}\cat{-CAlg}/\underline{B}.$ The map from left to right is $\underline{B} \ltimes (-)$ and the map from right to left is the augmentation ideal.
\end{corollary}

\section{Derivations}
The above equivalence leads us naturally to derivations. In the classical setting, a map from ${S} \to{B} \ltimes M$ in  ${A}\cat{-CAlg}/{B}$ is a map of abelian groups and decomposes to $s \mapsto (\epsilon(s),d(s)).$ The map $\epsilon$ must be a map of rings. We see that $A \to B \ltimes M$ maps to 0 in the $M$ component, so the map $A \to S \to B \ltimes M \to M$ is 0. Finally, if it is a map of rings then \begin{eqnarray*}
(\epsilon(s_1s_2),d(s_1s_2)) &=& (\epsilon(s_1),d(s_1)) \cdot (\epsilon(s_1),d(s_1)) \\
&=& (\epsilon(s_1)\epsilon(s_2), \epsilon(s_1)d(s_2) +\epsilon(s_2)d(s_1))
\end{eqnarray*}
That is, $d: S \to M$ is an $A$-derivation, where $M$ is considered an $S$-module via $S \to B$. If $d: S \to M$ is an $A$-derivation, then the map $s \mapsto (\epsilon(s),d(s))$ is a map of rings and therefore a map in ${A}\cat{-CAlg}/{B}.$

Following \cite{HillAndre}, we define $\mathcal{O}$-genuine derivations so as to maintain the same property. 
\begin{definition}
Let $\underline{A}$ be an $\mathcal{O}$-Tambara functor, $\underline{S}$ be an element of $\underline{A}\cat{-CAlg},$ and $\underline{M}$ be an $\underline{S}$-module. Then a map $d: \underline{S} \to \underline{M}$ is a \term{$\mathcal{O}$-genuine $\underline{A}$-derivation} if
\begin{enumerate}
\item It is a map of Mackey functors.
\item $\underline{A} \to \underline{S} \to \underline{M}$ is the zero map.
\item The map turns all admissible norms and products into transfers and sums in the following way: Let $f: X \to Y$ be an admissible norm map in $\mathcal{O}$, including products. Let $$\pi_1, \pi_2: X \underset{Y}\times X - \Delta(X) \to X$$ be the projections. Then for $s \in \underline{S}(X)$, we require $$d(N_f(s)) = T_f(N_{\pi_2}R_{\pi_1}(s) \cdot d(s)).$$
\end{enumerate}
We may sometimes refer to these as derivations, $\mathcal{O}$-derivations, or $\underline{A}$-derivations when the meaning is clear.
\end{definition}
The following propositions of \cite{HillAndre} Section 4 now go through with almost no change from the complete to incomplete case:

\begin{prop}\label{lmaderiv}
If $d: \underline{S} \to \underline{M}$ is a $\mathcal{O}$-genuine $\underline{A}$-derivation then:
\begin{itemize}
\item If $i: \underline{R} \to \underline{S}$ is a map of $\mathcal{O}$-Tambara functors, then $d \circ i$ is a $\mathcal{O}$-genuine $\underline{A}$-derivation.
\item If $f: \underline{M} \to \underline{M'}$ is a map of $\underline{S}$-modules, then $f \circ d$ is a $\mathcal{O}$-genuine $\underline{A}$-derivation.
\item The kernel $\ker(d)$ is a sub-$\mathcal{O}$-Tambara functor of $\underline{S}$.
\end{itemize}
If $\underline{B}$ is an $\mathcal{O}$-Tambara functor, $\underline{M}$ is an $\underline{B}$-module, $\underline{S} \overset{\epsilon}\to \underline{B}$ is an element of $\underline{A}\cat{-CAlg}/\underline{B}$ and $d: \underline{S} \to \underline{M}$ is a map of Mackey functors, then the map $$\epsilon \ltimes d: \underline{S} \to \underline{B} \ltimes \underline{M}$$ is a map in $\underline{A}\cat{-CAlg}/\underline{B}$ if and only if $d$ is a $\mathcal{O}$-genuine $\underline{A}$-derivation. In particular, $$\Hom_{\underline{A}\cat{-CAlg}/\underline{B}} (\underline{S},  \underline{B} \ltimes \underline{M}) \cong \Der_{\underline{A},\mathcal{O}}(\underline{S},\underline{M}).$$
\end{prop}
Note that the last line justifies our definition. By the same argument as in the classical case, we have the following equivalence:
\begin{prop}
A map $d: \underline{B} \to \underline{W}$ is an $\underline{A}$-derivation if and only if the map $d$ is a map of $\underline{A}$-modules and satisfies $$d(N_f(b)) = T_f(N_{\pi_1}R_{\pi_2}(b) \cdot d(b)).$$ 
\end{prop}

%% file: kahlerdiff2.tex
\section{Kahler Differentials}
Let $\underline{A}$ be an $\mathcal{O}$-Tambara functor and let $\underline{B}$ be an $\underline{A}$-algebra. Consider the natural homomorphism
$$\mu: \underline{B} \underset{\underline{A}} \boxtimes \underline{B} \to \underline{B}$$ given by the pushout of the two identity maps $\underline{B} \overset{\id}\to \underline{B}$. This maps $T_f(b \otimes b') \mapsto T_f(b \cdot b')$ for $f: X \to Y$ and $b \in \underline{B}(X)$. The kernel of this map is an $\mathcal{O}$-Tambara ideal $\underline{I}$. $\underline{I}$ gives rise to two other families of  ideals. First, we have $\underline{I}^{2}$, the image of $\underline{I} \boxtimes \underline{I} \to \underline{I}$, and likewise $\underline{I}^{3}, \underline{I}^4,\dots$ Second, we have the ideal $\underline{I}^{>1}$ which is the smallest ideal of $\underline{B}$ containing every element $N_f(i)$, where $i \in \underline{I}(X)$ and $f: X \to Y$ is a 2-surjective $\mathcal{O}$-admissible map. This is the same definition as appears in \cite{HillAndre} for the complete case. Similarly $\underline{I}^{>2}$ is defined likewise except $f$ must be 3-surjective, meaning $f^{-1}(y)$ has cardinality at least 3 for all $y \in Y$, and so on. We might also refer to these ideals as $\underline{I}^{\geq 2},\underline{I}^{\geq 3}, \ldots$ for easier to read formulas. We define $\Omega_{\underline{A}|\underline{B}} := \underline{I}/\underline{I}^{>1}$ and call it the module of K\"ahler differentials. 

There is a Mackey functor map $\underline{B} \to \underline{I}$ given by $1 \otimes b - b \otimes 1$. This gives rise to a natural $\underline{A}$-derivation $\delta: \underline{B} \to \Omega_{\underline{A}|\underline{B}}$ given by $b \mapsto \overline{1 \otimes b - b \otimes 1}$. Note that just in the classical case, while $\underline{I}$ has two $\underline{B}$-module structures, one on the left and one on the right, when passing to the quotient of any module including $\underline{I}^2$ they give the same structure.

\begin{lemma}\label{lemgenerate}
$\underline{I}$ is generated by the image of $\delta(\underline{B})$ as an ideal of $\underline{B} \underset{\underline{A}} \boxtimes \underline{B}$. Therefore, $\Omega_{\underline{A}|\underline{B}}$ is generated by $\delta(\underline{B})$ as a $\underline{B}$-module.
\end{lemma}

\begin{proof}
Let $T_f(b \otimes b')$ be an element of $\underline{B} \underset{\underline{A}} \boxtimes \underline{B}$ such that $T_f(b b') = 0$. Then 
\begin{eqnarray*}
T_f(b \otimes b') &=& T_f(b \otimes b')  - T_f(bb') \otimes 1 \\
&=&  T_f(b \otimes b')  - T_f(bb' \otimes R_f(1)) \\
&=&  T_f(b \otimes b' - bb' \otimes 1) \\
&=&  T_f((b \otimes 1) \cdot (1 \otimes b' - b' \otimes 1)). 
\end{eqnarray*}
\end{proof}

\begin{prop}\label{kahler}
If $\underline{W}$ is a $\underline{B}$ module, then given an element $h \in \Hom_{\underline{B}}(\Omega_{\underline{A}|\underline{B}}, \underline{W})$, the composition $h \circ \delta$ is a $\underline{A}$-derivation. Furthermore, all derivations are of this form: the map $$\Hom_{\underline{B}}(\Omega_{\underline{A}|\underline{B}}, \underline{W}) \to \Der_{\underline{A}}(\underline{B},\underline{W})$$ is an isomorphism.
\end{prop}
We will denote $\underline{\Der}_{\underline{A}}(\underline{B},\underline{W}) := \underline{\Hom}_{\underline{B}}(\Omega_{\underline{A}|\underline{B}}, \underline{W})$, the internal hom object. This is a $\underline{B}$-module that includes the genuine derivations $\underline{\Der}_{\underline{A}}(\underline{B},\underline{W})(G/G) = \Der_{\underline{A}}(\underline{B},\underline{W}) .$

\begin{proof}
~\cref{lmaderiv} shows that the map is well defined. Suppose $h \circ \delta$ is the zero map. Because $\delta(\underline{B})$ generates $\Omega_{\underline{B}|\underline{A}}$, $h$ must be 0. Therefore, the map is injective.

Suppose we have a derivation $d: \underline{B} \to \underline{W}$. We will show there is a map $h$ such that $d = h \circ \delta$. Firstly, we define $h: \underline{I} \to \underline{W}$ by $$h(T_q(b \otimes b')) =T_q( b \cdot d(b')).$$
We will show that this sends $\underline{I}^{>1}$ to 0. Once that is shown, the proposition is clear, as $$h \circ \delta(b) = h(1\otimes b- b\otimes 1) = d(b) - b d(1) = d(b).$$Suppose $g: X \to Y$ is 2-surjective and $Y$ decomposes into orbits $\coprod_i^n Y_i$. Then we get maps $g_i: X_i \to Y_i$ that are 2-surjective. Therefore, $N_g(1 \otimes (s_1,\ldots,s_n) - (s_1,\ldots,s_n) \otimes 1) = (N_{g_1} (1 \otimes s_1 - s_1 \otimes 1), \ldots , N_{g_n} (1 \otimes s_n - s_n \otimes 1)).$ In particular, this shows that $\underline{I}^{>1}$ is generated by  $N_g(1 \otimes s - s \otimes 1)$ where $g$ maps to an orbit. Additionally, every map onto an orbit $g: \coprod_j^n G/H_j \to G/H$ is the composition of $\coprod_j^n G/H_j \overset{h_j}\to \coprod_j^n G/H \overset{\nabla}{\to} G/H.$ If any of the $h_i$ are 2-surjective, then the fact that $\underline{I}^{>1}$ is an ideal means we can multiply by any element of $\underline{I}$ and stay in $\underline{I}^{>1}$. Thus $N_{h_j}(i_j)$ in the generating set implies $N_f(i)$ is in the ideal. If no $h_i$ are 2-surjective, then they are the identity, and our element must be $N_\nabla(i)$ for the map to be 2-surjective. In summary, $\underline{I}^{>1}$ is generated by $N_g(1 \otimes s - s \otimes 1)$ where $g$ is 2-surjective and is either an orbit mapping to an orbit $G/K \to G/H$ or a fold map.

With the fold map $\nabla: G/H \amalg G/H \to G/H$, the fact that $f(N_\nabla(1 \otimes s - s \otimes 1)) = 0$ is classical:
\begin{eqnarray*}
N_\nabla((1\otimes s_1 - s_1 \otimes 1, 1 \otimes s_2 - s_2 \otimes 1)) 
&=& (1\otimes s_1 - s_1 \otimes 1)( 1 \otimes s_2 - s_2 \otimes 1) \\
&=& 1\otimes s_1 s_2 - s_1 \otimes s_2 - s_2 \otimes s_1 + s_1 s_2 \otimes 1 \\
\end{eqnarray*}
giving
\begin{eqnarray*}
f(N_\nabla((1\otimes s_1 - s_1 \otimes 1, 1 \otimes s_2 - s_2 \otimes 1)))
&=& d(s_1 s_2) - s_1 d (s_2) - s_2 d( s_1)\\
&=& 0.
\end{eqnarray*}

Lastly, we need to show for $g: G/K \to G/H$, $f(N_g(1 \otimes s - s \otimes 1)) = 0$. Let us for a moment consider when $H = G$ and change our variables so that we are considering $K$-Tambara functors and the map $f: K/H \to K/K.$ We rewrite $N_f(1 \otimes s - s \otimes 1)$ in TNR-form again:
\[\xymatrix{
K/H \ar[d] &  K/H \coprod K/H \ar[l]& \{kH, s: K/H \to \{0,1\} \} \ar[d] \ar[l] \\
\star && \{s: K/H \to \{0,1\} \} \ar[ll]
	}\] 

We break up the 2 rightmost sets into the cardinality of $s^{-1}(1)$ indexed by the dummy variable $i$, and the top rightmost set into where $kH$ maps.

\[\xymatrixrowsep{.5in}
\xymatrixcolsep{.5in}\xymatrix{
K/H \ar[d] &
  K/H \coprod K/H \ar[l]& 
  \txt{ $\coprod_{i=0}^{|K/H|} ( \{kH, s | s^{-1}(1) = i, s(kH) =0\} $ \\  $\amalg \{kH, s | s^{-1}(1) = i, s(kH) =1\} )$ }  \ar[d] \ar[l] \\
\star && \coprod_{i=0}^{|K/H|}\{s: K/H \to \{0,1\} | s^{-1}(1) = i \} \ar[ll]
	}\] 

The key observation is that the map $$\{kH, s | s^{-1}(1) = i, s(kH) =1\} \to \{s | s^{-1}(1) = i \} $$ is precisely $i$-to-$1$. When we evaluate $N_f(s \otimes 1+(1 \otimes s - s\otimes 1))$ we get 
\begin{eqnarray*}
&&1 \otimes N_f(s)\\
&=& N_f[s \otimes 1 +(1 \otimes s - s\otimes 1)] \\
&=& c_0 + c_1 + c_2 + \ldots + c_{|K/H|-1} + c_{|K/H|}\\
&=& N_f(s) \otimes 1 + T_p(N_{\pi_2}^fR_{\pi_1}^f s \cdot (1\otimes s- s\otimes 1)) + c_2 + \ldots + c_{|K/H| - 1} + N_f(1 \otimes s- s\otimes 1) .
\end{eqnarray*}
Here we break up $N_f[1 \otimes s +(1 \otimes s - s\otimes 1)] $ by the $i$-index on the right two sets. $c_0, c_{|H/K|}$ can readily be seen to be the terms above $1 \otimes N_f(s) $ and $N_f(1 \otimes s- s\otimes 1)$ respectively. $c_1$ being $T_p(N_{\pi_2}^fR_{\pi_1}^f s \cdot (1\otimes s- s\otimes 1))$ is done in \cite{HillAndre} and not repeated here. Lastly, $c_i$ for all other $i$ is a multiple $N_g(R_\nabla(1 \otimes s - s\otimes 1))$ for where $R_\nabla$ is the restriction of some fold map, perhaps the identity, and $g$ an exactly $i$-to-$1$ map. 

One observation is because $G \underset{K}{\times} (-)$ respects exponential diagrams, we get the diagram

\[\xymatrixrowsep{.5in}
\xymatrixcolsep{.5in}\xymatrix{
G/H \ar[d]^f & 
 G/H \coprod G/H \ar[l]& 
 \txt{ $\coprod_{i=0}^{|K/H|} ( \{gH, s | s^{-1}(1) = i, s(gH) =0\} $ \\  $ \amalg \{gH, s | s^{-1}(1) = i, s(gH) =1\} )$ } \ar[d] \ar[l] \\
G/K && \coprod_{i=0}^{|K/H|}\{gK,s: f^{-1}(gK) \to \{0,1\} | s^{-1}(1) = i \} \ar[ll]
	}\] 

giving again
\begin{eqnarray*}
&&1 \otimes N_f(s)\\
&=& N_f[s \otimes 1 +(1 \otimes s - s\otimes 1)] \\
&=& c_0 + c_1 + c_2 + \ldots + c_{|K/H|-1} + c_{|K/H|}\\
&=& N_f(s) \otimes 1 + T_p(N_{\pi_2}^fR_{\pi_1}^f s \cdot (1\otimes s- s\otimes 1)) + c_2 + \ldots + c_{|K/H| - 1} + N_f(1 \otimes s- s\otimes 1) .
\end{eqnarray*}
with $c_i$ for all other $i$ is a multiple $N_g(R_\nabla(1 \otimes s - s\otimes 1))$.

We can now prove that $\underline{I}^{>1}$ is mapped to 0 by induction. The hypothesis, as a function of the $i$th step in the induction, is that $N_f(1 \otimes s - s \otimes 1)$ is mapped to 0 for all $f$ when $f$ is $j$-to-1 for $1 < j < i+1.$ For the first step, we have no hypothesis. The second step, we suppose the $N_f(1 \otimes s - s \otimes 1)$ maps to 0 for all $f$ 2-to-1 maps. The third step, all 2-to-1 and 3-to-1 maps, etc.

Suppose $h(N_f(1 \otimes s - s \otimes 1))=0$ for all $f$ when $f$ is $j$-to-1 for $1 < j < i+1.$ Now we consider $G/H \to G/K$ when $f$ is $(i+1)$-to-1. Then $c_2,\ldots,c_{|K/H|-1}$ are all 0 by the induction hypothesis. Therefore, 
$$h[1 \otimes N_f(s)] = h[N_f(s) \otimes 1 + T_p(N_{\pi_2}^fR_{\pi_1}^f s \cdot (1\otimes s- s\otimes 1))+ N_f(1 \otimes s- s\otimes 1) ].$$
Rearranging to 
\begin{eqnarray*}
h[N_f(1 \otimes s- s\otimes 1)] &=& h[1 \otimes N_f(s)] - h[N_f(s) \otimes 1] - h[T_p(N_{\pi_2}^fR_{\pi_1}^f s \cdot (1\otimes s- s\otimes 1))] \\
&=& 0
\end{eqnarray*} and then $h(N_f(1 \otimes s- s\otimes 1))=0$. Therefore, $N_f(1 \otimes s- s\otimes 1)$ is mapped to kernel for all $f$ $(i+1)$-to-1. By induction, $N_f(1 \otimes s- s\otimes 1)$ is mapped to 0 for all $f: G/H \to G/K$ and $f$ the fold map. This completes the proof that $\underline{I}^{>1}$ is mapped to 0.

\end{proof}

Suppose $M$ is a $\underline{B}$-module. Then we have the following isomorphisms:
\begin{eqnarray*}
\Hom_{\underline{A}\cat{-CAlg}/\underline{B}} (\underline{R},\underline{B} \ltimes \underline{M}) &\cong& \Der_{\underline{A},\mathcal{O}}(\underline{R},\underline{M}) \\ &\cong& \Hom_{\underline{R}-Mod}( \Omega_{\underline{R}/\underline{A}},\underline{M})\\ &\cong& \Hom_{\underline{B}-Mod}( \Omega_{\underline{R}/\underline{A}} \underset{\underline{R}} \boxtimes \underline{B},\underline{M})\\
\end{eqnarray*} 
Therefore there is a Mackeyization functor (this term is meant to emulate abelianization) between $$\underline{A}\cat{-CAlg}/\underline{B} \to (\underline{A}\cat{-CAlg}/\underline{B})_{\mathcal{O}-CommMon} \cong \cat{\underline{B}Mod}$$ where $(\underline{A}\cat{-CAlg}/\underline{B})_{\mathcal{O}-CommMon}$ are the $\mathcal{O}$-commutative monoid objects. The map takes  $$\underline{R} \mapsto  \Omega_{\underline{R}/\underline{A}}^{1,G,\mathcal{O}} \underset{\underline{R}} \boxtimes \underline{B}.$$ Using the equivalence of categories between $\underline{B}$-mod and the $\mathcal{O}$-commutative monoids of $\underline{A}\cat{-CAlg}/\underline{B} $, this is left adjoint to the inclusion functor of $$(\underline{A}\cat{-CAlg}/\underline{B})_{\mathcal{O}-CommMon}  \to \underline{A}\cat{-CAlg}/\underline{B}.$$

%% file: diffprop.tex
First we prove some functoriality properties of genuine derivations and K\"ahler differentials.

\begin{prop}\label{lemnat}
Suppose we have a commutative diagram of $\mathcal{O}$-Tambara functors 
\[\xymatrix{
\underline{A}' \ar[r] \ar[d]^{\alpha} &\underline{B}' \ar[d]^{\beta} \\
\underline{A} \ar[r] & \underline{B}
	}\] and a map $\omega: \underline{W'} \to \underline{W}$ of $\underline{B'}$-modules, with the $\underline{B'}$-module structure of $\underline{W}$ coming from a $\underline{B}$-module structure. Then there are natural $\underline{B'}$-module maps
$$ \underline{\Der}_{\underline{A}}(\underline{B}, \underline{W'}) \to  \underline{\Der}_{\underline{A'}}(\underline{B'}, \underline{W}) $$
$$ \Omega_{\underline{A'}|\underline{B'}} \underset{\underline{B'}} \boxtimes \underline{W'} \to   \Omega_{\underline{A}|\underline{B}} \underset{\underline{B}} \boxtimes \underline{W}.  $$
\end{prop}

\begin{proof}
There is the following commutative diagram
\[\xymatrix{ 
\underline{I'} \ar[d] \ar[r] & \underline{I} \ar[d]\\
\underline{B'} \underset{\underline{A'}} \boxtimes \underline{B'} \ar[r]^{\beta \otimes \beta} \ar[d]^{\mu'} &\underline{B} \underset{\underline{A}} \boxtimes \underline{B}  \ar[d]^{\mu} \\
\underline{B'} \ar[r]^{\beta} & \underline{B'}
	}\] and similarly $\underline{I'}^{>1} \to \underline{I}^{>1}$ inducing a map $\Omega_{\underline{A'}|\underline{B'}} \to \Omega_{\underline{A}|\underline{B}}$ which is a map of $\underline{B'}$-modules, as the map $\beta \otimes \beta$ is a map of $\underline{B'}$-modules (on the left). Therefore we get the series of maps 
\begin{eqnarray*}
\Omega_{\underline{A'}|\underline{B'}} \underset{\underline{B'}} \boxtimes \underline{W'} &\to& \Omega_{\underline{A}|\underline{B}} \underset{\underline{B'}} \boxtimes \underline{W'} \\
&\to& \Omega_{\underline{A}|\underline{B}} \underset{\underline{B'}} \boxtimes \underline{W} \\
&\cong& \Omega_{\underline{A}|\underline{B}} \underset{\underline{B}} \boxtimes \underline{W}.
\end{eqnarray*}
Let $f \in \underline{\Der}_{\underline{A}}(\underline{B}, \underline{W'})(X)$. Then $\omega_X \circ f \circ \beta \in  \underline{\Der}_{\underline{A'}}(\underline{B'}, \underline{W})(X).$
\end{proof}

\begin{prop}
Suppose $\{ \underline{W_k} \}$ is a set of $\underline{B}$-modules. Then there are natural isomorphisms $$\Omega_{\underline{A}|\underline{B}} \underset{\underline{B}} \boxtimes \left(\bigoplus_k \underline{W_k} \right) \cong \bigoplus_k \Omega_{\underline{A}|\underline{B}} \underset{\underline{B}} \boxtimes \underline{W_k} $$ $$\underline{\Der}_{\underline{A}}\left(\underline{B}, \prod_k \underline{W_k}\right) \cong \prod_k \underline{\Der}_{\underline{A}}(\underline{B},\underline{W_k})  $$
\end{prop}

\begin{proof}
The first is simply a property of the box product with respect to direct sums. The second follows from \cref{lmaderiv} and the fact that maps into a product is the product of maps.
\end{proof}

\begin{prop}
Suppose $\underline{C}$ an $\underline{A}$-algebra is a $\underline{W}$ an $\underline{B} \underset{\underline{A}} \boxtimes \underline{C}$-module. Then there are natural isomorphisms
 $$\Omega_{\underline{A}|\underline{B}} \underset{\underline{B}} \boxtimes \underline{W} \overset{\cong}\to \Omega_{\underline{C}|\underline{B} \underset{\underline{A}} \boxtimes \underline{C}} \underset{\underline{B} \underset{\underline{A}} \boxtimes \underline{C}} \boxtimes \underline{W} $$ $$ \underline{\Der}_{\underline{C}}(\underline{B} \underset{\underline{A}} \boxtimes \underline{C},\underline{W}) \overset{\cong}\to \underline{\Der}_{\underline{A}}(\underline{B}, \underline{W}) $$
\end{prop}

\begin{proof} Given the pushout diagram
\[\xymatrix{ 
\underline{A} \ar[d] \ar[r] & \underline{B} \ar[d] \\
\underline{C} \ar[r] & \underline{B} \underset{\underline{A}} \boxtimes \underline{C}
	}\] 
we have natural maps in the statement given by ~\cref{lemnat}. To show the second is an isomorphism, consider the isomorphism $$\underline{\Hom}_{\underline{C}-Mod}(\underline{B} \underset{\underline{A}}\boxtimes \underline{C}, \underline{W}) \cong \underline{\Hom}_{\underline{A}-Mod}(\underline{B},\underline{W}).$$
Suppose $f \in \underline{\Hom}_{\underline{C}-Mod}(\underline{B} \underset{\underline{A}}\boxtimes \underline{C}, \underline{W}) $ is additionally in $\underline{\Der}_{\underline{C}}(\underline{B} \underset{\underline{A}} \boxtimes \underline{C},\underline{W}) $. Then $f$ takes norms to the appropriate transfers. We show that the isomorphism extends to these subsets.

If the composition $\underline{B} \to \underline{B} \underset{\underline{A}}\boxtimes \underline{C} \to \underline{W}$ is 0, then the map $\underline{B} \underset{\underline{A}}\boxtimes \underline{C} \to \underline{W}$ is 0 because $T_q(c \otimes b) \mapsto T_q( c \cdot 0) = 0$. Thus the map on derivations is injective. Similarly if $g \in \underline{\Hom}_{\underline{A}-Mod}(\underline{B},\underline{W})$ takes norms to transfers, then the induced map in $\underline{\Hom}_{\underline{C}-Mod}(\underline{B} \underset{\underline{A}}\boxtimes \underline{C}, \underline{W})$, that is $T_q(c \otimes b) \mapsto T_q(c \cdot g(b)),$ takes norms to transfers by the laborious calculation in the next lemma. For the first isomorphism, we use the second isomorphism and the universal property of the module of K\"ahler differentials:
\begin{eqnarray*}
\Der_{\underline{C}}(\underline{B} \underset{\underline{A}}\boxtimes \underline{C}, \underline{W}) &\cong& \Der_{\underline{A}}(\underline{B},\underline{W}) \\
\Hom_{\underline{B} \underset{\underline{A}}\boxtimes \underline{C}}(\Omega_{\underline{C}|\underline{B} \underset{\underline{A}}\boxtimes \underline{C}}, \underline{W}) &\cong& \Hom_{\underline{B} }(\Omega_{\underline{A}|\underline{B}}, \underline{W}) \\
&\cong&\Hom_{\underline{B} \underset{\underline{A}}\boxtimes \underline{C}}(\Omega_{\underline{A}|\underline{B}} \underset{\underline{A}} \boxtimes \underline{C}, \underline{W}) \\
\Omega_{\underline{C}|\underline{B} \underset{\underline{A}}\boxtimes \underline{C}} &\cong& \Omega_{\underline{A}|\underline{B}} \underset{\underline{A}} \boxtimes \underline{C}\\
 \Omega_{\underline{A}|\underline{B}} \underset{\underline{B}} \boxtimes \underline{W}&\cong& \Omega_{\underline{C}|\underline{B} \underset{\underline{A}}\boxtimes \underline{C}} \underset{\underline{B} \underset{\underline{A}}\boxtimes \underline{C}} \boxtimes  \underline{W}
\end{eqnarray*}

\end{proof}

\begin{lemma}\label{lemlong}
Suppose $d \in \Der_{\underline{A}}(\underline{B},\underline{W})$ and let $\tilde{d}: \underline{C} \underset{\underline{A}}{\boxtimes} \underline{B} \to \underline{W}$ be the induced map from $d$, $$\tilde{d}: T_p(c \otimes b) \mapsto T_p(c \cdot d(b)).$$ Then $\tilde{d}$ is a genuine derivation for the $\mathcal{O}$-Tambara functor map $ \underline{C} {\to} \underline{C} \underset{\underline{A}}{\boxtimes} \underline{B}.$
\end{lemma}
\begin{proof}
This is a map of Mackey functors: it is clear it commutes with transfers and the fact that $d$ commutes with restrictions implies $\tilde{d}$ commutes with restrictions: Let $R_p T_q = T_a R_b,$
\begin{eqnarray*}
\tilde{d} \left(R_p T_q(c\otimes b)\right) &=& \tilde{d} \left(T_a R_b(c\otimes b)\right) \\
&=& T_a (R_b c \cdot d(R_b b)) \\
&=& T_a R_b (c \cdot d(b)) \\
&=& R_p T_q (c \cdot d(b)) \\
&=& R_p \tilde{d} [ T_q (c \otimes b) ].
\end{eqnarray*}
For the $\underline{C}$-action, compute
\begin{eqnarray*}
\tilde{d} \left(c' \cdot T_q(c\otimes b)\right) &=& \tilde{d} \left( T_q(R_q(c') \cdot c \otimes b)\right) \\
&=&T_q(R_q(c') \cdot c \cdot d(b)) \\
&=& c' \cdot T_q(c \cdot d(b)) \\
&=& c' \cdot \tilde{d}\left(T_q(c \otimes b)\right).
\end{eqnarray*}
Additionally, $\underline{C}$ maps to 0, since $d(1) = 0$. We simply need to show the multiplicative property.
We want to show that (equations we want to show and are not yet evident are labeled by `WTS'):
\begin{eqnarray*}
\tilde{d}(N_f(T_p(c \otimes b))) &\overset{WTS}{=}& T_f(N_{\pi_2}^fR_{\pi_1}^f(T_p(c\otimes b)) \cdot \tilde{d}(T_p(c\otimes b))) \\
&=&  T_f(N_{\pi_2}^fR_{\pi_1}^f(T_p(c\cdot b)) \cdot T_p(c\cdot d(b))). 
\end{eqnarray*}
To evaluate the left side, we need to put $N_f T_p$ in TNR form. Let $N_f T_p = T_aN_bR_c,$ then we want to show that \begin{eqnarray*}
\tilde{d}(N_f(T_p(c \otimes b))) &=& \tilde{d}(T_aN_bR_c(c \otimes b))  \\
&=&  \tilde{d}(T_a(N_bR_c c \otimes N_b R_c b))  \\
&=&  T_a(N_bR_cc \cdot d(N_b R_c b))  \\
&=& T_a(N_bR_cc \cdot T_b(N_{\pi_2}^bR_{\pi_1}^b R_c b \cdot d(R_c b))) \\ 
&=& T_aT_b(N_{\pi_2}^bR_{\pi_1}^b (R_c c \otimes R_c b) \cdot R_c c \cdot d(R_c b)) \\ 
&\overset{WTS}{=}&  T_f(N_{\pi_2}^fR_{\pi_1}^fT_p(c \otimes b) \cdot T_p(c\cdot d(b))). \\
\end{eqnarray*}
The author believes the only method is to convert $T_f(N_{\pi_2}^fR_{\pi_1}^f(T_p(c\cdot b)) \cdot T_p(c\cdot d(b)))$ into TNR-form. In the following diagrams, all the squares are pull backs, and all the diagrams with 2 arrows on top and 1 on bottom are exponential diagrams.

\begin{eqnarray*}
&& T_f(N_{\pi_2}^fR_{\pi_1}^fT_p(c\otimes b) \cdot T_p(c\cdot d(b)))\\ &=& T_fN_\nabla[N_{\pi_2}^fR_{\pi_1}^fT_pN_\nabla(c \otimes 1,1 \otimes b),T_pN_\nabla(c,d(b))] \\
&=& T_fN_\nabla[N_{\pi_2}^fT_{\tilde{p}}R_{\tilde{\pi}_1}(c \otimes b),T_p(c \cdot d(b))]  \textrm{ via Diagram 1}\\
&=& T_fN_\nabla[T_{a'}N_{b'}R_{c'}R_{\tilde{\pi}_1}(c \otimes b),T_p(c \cdot d(b))]  \textrm{ via Diagram 2}\\
&=& T_fT_{a''}N_{\nabla}[R_{c''_0}N_{b'}R_{c'}R_{\tilde{\pi}_1}(c \otimes b),R_{c}(c\cdot d(b))]  \textrm{ via Diagram 3}\\
&=& T_fT_{a''}N_{\nabla}[N_{b''}R_{c'''}R_{c'}R_{\tilde{\pi}_1}(c \otimes b),R_{c}(c\cdot d(b))]  \textrm{ via Diagram 4}\\
&=& T_fT_{a''}N_{\nabla}[N_{\pi_2}^bR_{\pi_1}^bR_{c}(c \otimes b),R_{c}(c \cdot d(b))]  \textrm{ via Diagram 5,6}\\
&=& T_aT_{b}[N_{\pi_2}^bR_{\pi_1}^bR_{c}(c \otimes b) \cdot R_{c}(c \cdot d(b))].\\
\end{eqnarray*}
Along with the complementary diagrams below, this completes the proof that this is a genuine derivation.

\begin{tcolorbox}[standard jigsaw,
    title=Exponential Diagram 0:,
    opacityback=0,]
\[\xymatrix{
 Y \ar[d]^f& X \ar[l]^p & \ar[l]^(.75)c \ar[d]^b\{(y,s: f^{-1}(f(y)) \to X )\} \\
Z && \{(z, f^{-1}(z) \to X) \} \ar[ll]^(.65)a
	}\] 
\end{tcolorbox}

\begin{tcolorbox}[standard jigsaw,
    title= Diagram 1:,
    opacityback=0,]
\[\xymatrix{
\txt{$X \underset{Z}{\times} Y - p(X)$ \\ $:= \{ (x,y) | f(p(x)) = f(y),p(x) \neq y\}$} \ar[r]^(.52){\tilde{p}} \ar[d]^(.6){\tilde{\pi}_1}& \txt{$Y \underset{Z}{\times} Y - Y $ \\  $:=\{ (y,y') | f(y) = f(y'), y \neq y'\} $}\ar[d]^(.6){\pi_1}\\
X \ar[r]^p & Y\\
	}\] 
\end{tcolorbox}

\begin{tcolorbox}[standard jigsaw,
    title= Diagram 2:,
    opacityback=0,]
\[\xymatrix{
Y \underset{Z}{\times} Y - Y \ar[d]^{\pi_2}& X \underset{Z}{\times} Y - p(X) \ar[l]^{\tilde{p}}& \txt{$\{(y,y',s: f^{-1}(f(y')) \setminus \{y'\} \to X) |$ \\ $ y \neq y', f(y) = f(y') \}$} \ar[l]^-{c'} \ar[d]^{b'} \\
Y && \{(y',s: f^{-1}(f(y')) \setminus \{y'\} \to X)\} \ar[ll]^(.65){a'} \\
	}\] 
\end{tcolorbox}	

\begin{tcolorbox}[standard jigsaw,
    title= Diagram 3:,
    opacityback=0,]
\[\xymatrixrowsep{.25in}
\xymatrixcolsep{.3in}
\xymatrix{
Y \coprod Y \ar[d]^{\nabla}& \{(y,s: f^{-1}(f(y)) - \{y\} \to X)\} \amalg X \ar[l]^(.78){(a',p)}&\txt{$\{(y,s: f^{-1}(f(y)) \to X)\}$ \\ $\amalg  \{(y,s: f^{-1}(f(y)) \to X)\}$} \ar[l]^-{(c''_0,c)} \ar[d]^(.6){\nabla} \\
Y &&  \ar[ll]^{a''} \{(y,s: f^{-1}(f(y)) \to X)\}\\
	}\] 
\end{tcolorbox}	
\begin{tcolorbox}[standard jigsaw,
    title= Diagram 4:,
    opacityback=0,]
\[\xymatrixrowsep{.25in}
\xymatrixcolsep{.2in}
\xymatrix{
\{(y,y',s: f^{-1}(f(y')) \to X) | y \neq y', f(y) = f(y')\} \ar[r]^(.64){b'' = \pi_2} \ar[d]^{c'''}& \{(y',s: f^{-1}(f(y')) \to X)\} \ar[d]^{c''_0}\\
\{(y,y',s: f^{-1}(f(y')) - y \to X) | y \neq y', f(y) = f(y')\} \ar[r]^>>{b'}& \{(y',s: f^{-1}(f(y')) - y' \to X)\} \\
	}\] 
\end{tcolorbox}	
\begin{tcolorbox}
[standard jigsaw,
    title= Diagram 5:,
    opacityback=0,]
\[\xymatrixrowsep{.25in}
\xymatrixcolsep{.2in}
\xymatrix{
\{ (y,y',s: f^{-1}(f(y')) \to X)  | y \neq y', f(y) = f(y')\} \ar[r]^(.63){\pi_1}&\{ (y',s: f^{-1}(f(y')) \to X) \}  \ar[r]^(.82){c} &X \\
	}\] 
\end{tcolorbox}	
		\begin{tcolorbox}[standard jigsaw,
    title= Diagram 6:,
    opacityback=0,]
\[\xymatrixrowsep{.2in}
\xymatrixcolsep{0in}
\xymatrix{
&\txt{$ \{(y,y',s:f^{-1}(f(y)) \to X) |$ \\ $y \neq y', f(y) = f(y')\}$} \ar[dr]_{\pi_1} \ar[dl]^{\pi_2 = b''}\\
\{(y,s:f^{-1}(f(y)) \to X) \} &&\{(y',s:f^{-1}(f(y')) \to X) \} \\
	}\] 
\end{tcolorbox}

\end{proof}

\begin{lemma}\label{lemprod}
There are natural isomorphisms
$$\Omega_{\underline{A}|\underline{B}} \oplus \Omega_{\underline{A}|\underline{C}} \cong \Omega_{\underline{A}|\underline{B} \underset{\underline{A}} \boxtimes \underline{C}}$$
$$\underline{\Der}_{\underline{A}}(\underline{B},\underline{W}) \oplus \underline{\Der}_{\underline{A}}(\underline{C},\underline{W}) \cong \underline{\Der}_{\underline{A}}(\underline{B} \underset{\underline{A}} \boxtimes \underline{C},\underline{W}).$$
\end{lemma}

\begin{proof}
We only show the property for $\underline{\Der}$, as the property for $\Omega$ follows by the universal property again. If $d$ is a derivation from $\underline{B} \underset{\underline{A}} \boxtimes \underline{C} \to \underline{W}$, then $$d(T_p(b \otimes c)) = T_p((b \otimes 1) \cdot d(1\otimes c) +(1 \otimes c) \cdot d(b\otimes 1) )$$ so the map of the lemma is injective. It is surjective by using ~\cref{lemlong} and by adding the two respective derivations to get a derivation from $\underline{B} \underset{\underline{A}} \boxtimes \underline{C}$, as the sum of two derivations is a derivation.
\end{proof}

%% file: KahlerOfFree2.tex
\subsection{Kahler differentials for Free Tambara functors} The canonical example to compute in the classical setting is the module of K\"ahler differentials for the inclusion of $R \to R[x]$. In this case, the module of K\"ahler differentials is the free $R[x]$-module $R[x] \cdot \{dx\}$.
We want to model $R[x]$. Maps from $\Z[x]$ to another ring $R$ are canonically isomorphic to elements of $R$, given by the image of $x$. In the equivariant setting, we have rings for every subgroup $H$ so we expect to have a family of free Tambara functors, indexed by $H$. Indeed, this definition is introduced in \cite{incomplete}.

\begin{defn}
For $H \subset G$, let $$\underline{A}^{\mathcal{O}}[x_H] = P_{\mathcal{O}}^G(G/H, -)$$ be the $\mathcal{O}$-Tambara functor represented by $G/H$. By the Yoneda lemma, $$\mathit{O-Tamb}(\underline{A}^{\mathcal{O}}[x_H], \underline{R}) \cong \underline{R}(G/H).$$
\end{defn}

Classically $R[x] = R \otimes \Z[x]$ and $R$-derivations from $R[x] \to M$ correspond with elements of $M$, the image of $dx$. This extends to the equivariant case.

\begin{prop} Let $\phi: \underline{R} \to \underline{R} \boxtimes \underline{A}^{\mathcal{O}}[x_H]$ be the inclusion $r \mapsto r \otimes 1$. Then there is an isomorphism
$$\Der_\phi(R \boxtimes \underline{A}^{\mathcal{O}}[x_H],\underline{M})  \overset{\cong}{\longrightarrow} \underline{M}(G/H).$$ The map is $d \mapsto d(1 \otimes \id)$, where $\id = T_{\id}N_{\id} R_{\id}$ is the identity polynomial.
\end{prop}

\begin{proof}
The image of $d(1 \otimes \id)$ completely determines the derivation $d$, because every element of $R \boxtimes \underline{A}^{\mathcal{O}}[x_H]$  is of the form $T_p(r \otimes N_bR_c)$ by \cite{StricklandTambara} Prop 3.14.  We necessarily have (using $d(r \otimes 1) = 0$):
\begin{eqnarray*}
d(T_p(r \otimes N_bR_c)) &=& T_p(d( (r\otimes 1) \cdot N_bR_c(1 \otimes \id))) \\
&=& T_p((r\otimes 1) \cdot d(N_bR_c(1 \otimes \id)))\\ &=& T_p((r\otimes 1) \cdot T_b(N_{\pi_2}^b R_{\pi_1}^b R_c(1 \otimes \id) \cdot R_c d(1\otimes \id))) .
\end{eqnarray*} This shows that the above map is injective. To show a bijection, we must show that the equation above does define a derivation for any element that we defined $d(1 \otimes \id)$ to be in $\underline{M}(G/H).$

From the definition, the $\underline{R}$ in $\underline{R} \boxtimes \underline{A}^{\mathcal{O}}[x_H]$ is mapped to 0:
\begin{eqnarray*}
d(r \otimes N_{\emptyset} R_{\emptyset} \Id) &=& (r \otimes 1) \cdot T_{\emptyset}(N^{\emptyset}_{\pi_2} R_{\pi_2}^\emptyset R_{\emptyset} (1 \otimes \Id) \cdot R_\empty d(1 \otimes \Id)) \\
&=& 0
\end{eqnarray*}
because $T_{\emptyset}$ maps every element to 0.
Transfers trivially commute with $d$. To show that restrictions commute with $d$, and that taking the norm of an element respects the multiplicative condition, requires we rearrange polynomials into TNR-form using exponential diagrams and pullbacks. Both these proofs are laborious but straightforward. 

First we show the commutativity with restrictions. We need to show that (maps defined by diagrams below) $$R_f d(T_a(r \otimes N_bR_c)) = d( R_f T_a(r \otimes N_bR_c)).$$
Let $R_f T_a = T_{\tilde{a}} R_{\tilde{f}}$ and $R_{\tilde{f}} N_b = N_{\tilde{b}} R_g.$ Then the right hand side is 
\begin{eqnarray*}
d( R_f T_a (r \otimes N_b R_c)) &=& d( T_{\tilde{a}}  R_{\tilde{f}} (r \otimes N_b R_c)) \\
&=& d( T_{\tilde{a}}   ( R_{\tilde{f}} r \otimes R_{\tilde{f}} N_b R_c)) \\
&=& d( T_{\tilde{a}}   ( R_{\tilde{f}} r \otimes N_{\tilde{b}}R_g R_c)) \\
&=& d( T_{\tilde{a}}   ( R_{\tilde{f}} r \otimes N_{\tilde{b}}R_g R_c)) \\
&=& T_{\tilde{a}} ((R_{\tilde{f}} r \otimes 1) \cdot T_{\tilde{b}} (N_{\pi_2}^{\tilde{b}} R_{\pi_1}^{\tilde{b}} R_g R_c \cdot R_g R_c d(\id)))\\
 R_f d(T_a (r \otimes N_b R_c)) &=& R_f T_a ((r \otimes 1) \cdot T_b (N_{\pi_2}^b R_{\pi_1}^b R_c \cdot R_c d(\id))) \\
  &=&  T_{\tilde{a}}  ( R_{\tilde{f}} (r \otimes 1) \cdot T_{\tilde{b}} R_g (N_{\pi_2}^b R_{\pi_1}^b R_c \cdot R_c d(\id))) \\
    &=&  T_{\tilde{a}}  ( R_{\tilde{f}} (r \otimes 1) \cdot T_{\tilde{b}} (R_g N_{\pi_2}^b R_{\pi_1}^b R_c \cdot R_c d(\id))) \\
\end{eqnarray*}
It remains to show that $ N_{\pi_2}^{\tilde{b}} R_{\pi_1}^{\tilde{b}} R_g= R_{g}  N_{\pi_2}^b R_{\pi_1}^b.$ We have the following diagrams:

\begin{tcolorbox}[
    standard jigsaw,
    title=Diagram 1,
    opacityback=0,  
]
\[\xymatrix{
X \underset{Z}{\times} Y  \ar[d]^g \ar[r]^{\tilde{b}}& Y \ar[d]^{\tilde{f}}\\
X \ar[r]^b& Z 
	}\] 
\end{tcolorbox}	

\begin{tcolorbox}[
    standard jigsaw,
    title=Diagram 2,
    opacityback=0,  
]
\[\xymatrix{
& &(X \underset{Z}{\times} Y) \underset{Y}{\times} (X \underset{Z}{\times} Y) - \Delta(X \underset{Z}{\times} Y )  \ar[dl]^{R_{\pi_1}^{\tilde{b}}} \ar[dr]_{N_{\pi_2}^{\tilde{b}}}& \\
&X \underset{Z}{\times} Y \ar[dl]^{R_g} \ar[dr] & &X \underset{Z}{\times} Y \ar[dl]\\
X& &X & \\
}\]
\end{tcolorbox}	

\begin{tcolorbox}[
    standard jigsaw,
    title=Diagram 3,
    opacityback=0,  
]
\[\xymatrix{
& &(X \underset{Z}{\times} Y) \underset{Y}{\times} (X \underset{Z}{\times} Y) - \Delta(X \underset{Z}{\times} Y )  \ar[dl] \ar[dr]& \\
&X \underset{Z}{\times} X - \Delta(X) \ar[dl]^{R_{\pi_1}^b} \ar[dr]^{N_{\pi_2}^b} & &X \underset{Z}{\times} Y \ar[dl]^{R_g}\\
X& &X & \\
}\]
\end{tcolorbox}	

In the third diagram, the square is a pullback. The pullback would be two distinct elements of $X$ and an element of $Y$, all mapping to the same element in $Z$, which is precisely what $(X \underset{Z}{\times} Y) \underset{Y}{\times} (X \underset{Z}{\times} Y) - \Delta(X \underset{Z}{\times} Y )$ is. Tracing through the maps, we see that $ N_{\pi_2}^{\tilde{b}} R_{\pi_1}^{\tilde{b}} R_g= R_{g}  N_{\pi_2}^b R_{\pi_1}^b.$ 

Now we have to show the derivation respects norms. We want to calculate (maps defined by diagrams below)
\begin{eqnarray*}
&& d(N_\gamma T_\beta (r \otimes N_\alpha R_k))\\ &=& d(T_{\alpha'} N_{\beta'} R_{\gamma'}  (r \otimes N_\alpha R_k)) \\
&=& d(T_{\alpha'} N_{\beta'}   ( R_{\gamma'} r \otimes  R_{\gamma'} N_\alpha R_k)) \\
&=& d(T_{\alpha'} N_{\beta'}   ( R_{\gamma'} r \otimes  N_{\tilde{\alpha}} R_{\tilde{\gamma'}} R_k)) \\
&=& d(T_{\alpha'}   ( N_{\beta'}  R_{\gamma'} r \otimes  N_{\beta'}  N_{\tilde{\alpha}} R_{\tilde{\gamma'}} R_k)) \\
&=& T_{\alpha'} ((N_{\beta'}  R_{\gamma'} r \otimes \id) \cdot T_{\beta' \circ \tilde{\alpha}} ((1 \otimes N_{\pi_2}^{\beta' \circ \tilde{\alpha}} R_{\pi_1}^{\beta' \circ \tilde{\alpha}} R_{\tilde{\gamma'}} R_k )\cdot R_{\tilde{\gamma'}} R_k d(1 \otimes \id)))\\
&=& T_{\alpha'}  T_{\beta' \circ \tilde{\alpha}} ((R_{\beta' \circ \tilde{\alpha}} N_{\beta'}  R_{\gamma'} r \otimes N_{\pi_2}^{\beta' \circ \tilde{\alpha}} R_{\pi_1}^{\beta' \circ \tilde{\alpha}} R_{\tilde{\gamma'}} R_k )\cdot R_{\tilde{\gamma'}} R_k d(1 \otimes \id))\\
&\overset{WTS}{=}& T_\gamma( N_{\pi_2}^\gamma R_{\pi_1}^\gamma T_\beta (r \otimes N_\alpha R_k) \cdot T_\beta ((r \otimes \id) \cdot T_\alpha ((1 \otimes N_{\pi_2}^\alpha R_{\pi_1}^\alpha R_k )\cdot R_k d(1 \otimes \id))))\\
\end{eqnarray*}
If the above equations are equal, then $d$ satisfies the multiplicative property. We rearrange the last expression in TNR-form:
\begin{eqnarray*}
&& T_\gamma( N_{\pi_2}^\gamma R_{\pi_1}^\gamma T_\beta (r \otimes N_\alpha R_k) \cdot T_\beta ((r \otimes \id) \cdot T_\alpha ((1 \otimes N_{\pi_2}^\alpha R_{\pi_1}^\alpha R_k )\cdot R_k d(1 \otimes \id)))) \\
&=& T_\gamma( N_{\pi_2}^\gamma T_{\tilde{\beta}} R_\lambda (r \otimes N_\alpha R_k) \cdot T_\beta ((r \otimes \id) \cdot T_\alpha ((1 \otimes N_{\pi_2}^\alpha R_{\pi_1}^\alpha R_k )\cdot R_k d(1 \otimes \id)))) \\
&=& T_\gamma( T_{\alpha''} N_{\beta''} R_{\gamma''} R_\lambda (r \otimes N_\alpha R_k) \cdot T_\beta ((r \otimes \id) \cdot T_\alpha ((1 \otimes N_{\pi_2}^\alpha R_{\pi_1}^\alpha R_k ) \cdot R_k d(1 \otimes \id)))) \\
&=& T_\gamma( T_{\alpha''} N_{\beta''} R_{\gamma''} R_\lambda (r \otimes N_\alpha R_k) \cdot T_\beta T_\alpha ((R_\alpha(r) \otimes N_{\pi_2}^\alpha R_{\pi_1}^\alpha R_k) \cdot R_k d(1 \otimes \id))) \\
&=& T_\gamma T_{\alpha'''} N_{\nabla} R_{\gamma'''} (  [N_{\beta''} R_{\gamma''} R_\lambda (r \otimes N_\alpha R_k)], [((R_\alpha(r) \otimes N_{\pi_2}^\alpha R_{\pi_1}^\alpha R_k) \cdot R_k d(1 \otimes \id))]) \\
&=& T_\gamma T_{\alpha'''} N_{\nabla}  (  [ R_{\gamma'''_0}N_{\beta''} R_{\gamma''} R_\lambda (r \otimes N_\alpha R_k)], [R_{\gamma'''_1}(R_\alpha(r) \otimes N_{\pi_2}^\alpha R_{\pi_1}^\alpha R_k) \cdot R_{\gamma'''_1}R_k d(1 \otimes \id)]) \\
&=& T_\gamma T_{\alpha'''} (    (R_{\gamma'''_0}N_{\beta''} R_{\gamma''} R_\lambda r \otimes R_{\gamma'''_0}N_{\beta''} R_{\gamma''} R_\lambda N_\alpha R_k) \\ && \cdot (R_{\gamma'''_1}R_\alpha(r) \otimes  R_{\gamma'''_1}N_{\pi_2}^\alpha R_{\pi_1}^\alpha R_k) \cdot R_{\gamma'''_1}R_k d(1 \otimes \id)) \\
&=& T_\gamma T_{\alpha'''} (    [(R_{\gamma'''_0}N_{\beta''} R_{\gamma''} R_\lambda (r) \cdot R_{\gamma'''_1}R_\alpha(r)) \otimes (R_{\gamma'''_0}N_{\beta''} R_{\gamma''} R_\lambda N_\alpha R_k \\ && \cdot R_{\gamma'''_1}N_{\pi_2}^\alpha R_{\pi_1}^\alpha R_k) ]  \cdot R_{\gamma'''_1}R_k d(1 \otimes \id)) \\
&=& T_\gamma T_{\alpha'''} (    [(R_{\gamma'''_0}N_{\beta''} R_{\gamma''} R_\lambda (r) \cdot R_{\gamma'''_1}R_\alpha(r)) \otimes (N_{\epsilon} R_\delta R_k \cdot N_\phi R_{\tilde{\gamma}'''_1} R_{\pi_1}^\alpha R_k) ]  \cdot R_{\gamma'''_1}R_k d(1 \otimes \id)) \\
\end{eqnarray*}

\begin{eqnarray*}
&&T_\gamma( N_{\pi_2}^\gamma R_{\pi_1}^\gamma T_\beta N_\alpha R_k \cdot T_\beta T_\alpha (N_{\pi_2}^\alpha R_{\pi_1}^\alpha R_k \cdot R_k d(\id)) )\\
&=& T_\gamma( N_{\pi_2}^\gamma T_{\tilde{\beta}} R_\lambda  N_\alpha R_k \cdot T_\beta T_\alpha (N_{\pi_2}^\alpha R_{\pi_1}^\alpha R_k \cdot R_k d(\id))) \\
&=& T_\gamma( T_{\alpha''} N_{\beta''} R_{\gamma''} R_\lambda  N_\alpha R_k \cdot T_\beta T_\alpha (N_{\pi_2}^\alpha R_{\pi_1}^\alpha R_k \cdot R_k d(\id))) \\
&=& T_\gamma T_{\alpha'''} N_{\nabla} R_{\gamma'''} (N_{\beta''} R_{\gamma''} R_\lambda  N_\alpha R_k, N_{\pi_2}^\alpha R_{\pi_1}^\alpha R_k \cdot R_k d(\id)) \\
&=& T_\gamma T_{\alpha'''}  N_{\nabla}( R_{\gamma'''_0} N_{\beta''} R_{\gamma''} R_\lambda  N_\alpha R_k, R_{\gamma'''_1}N_{\pi_2}^\alpha R_{\pi_1}^\alpha R_k \cdot R_{\gamma'''_1} R_k d(\id))\\
&=& T_\gamma T_{\alpha'''} ( (N_{\epsilon} R_\delta R_k \cdot N_\phi R_{\tilde{\gamma}'''_1} R_{\pi_1}^\alpha R_k) \cdot R_{\gamma'''_1} R_k d(\id))
\end{eqnarray*}
The following diagrams validate the above rearrangement:
\begin{tcolorbox}[
    standard jigsaw,
    title=Diagram 1,
    opacityback=0,  
]
\[\xymatrix{
A & B  \ar[l]_k \ar[r]^{\alpha}& C \ar[r]^\beta &D \ar[r]^\gamma& E \\
}\]
\end{tcolorbox}	

\begin{tcolorbox}[
    standard jigsaw,
    title=Diagram 2,
    opacityback=0,  
]
\[\xymatrix{
D \ar[d]^\gamma & C \ar[l]^\beta& \{(d,s:\gamma^{-1}(\gamma(d))) \to C)\} \ar[d]^{\beta'} \ar[l]^(.8){\gamma'}\\
E && \{(e,s:\gamma^{-1}(e) \to C)\} \ar[ll]^(.70){\alpha'}
}\]
\end{tcolorbox}	

\begin{tcolorbox}[
    standard jigsaw,
    title=Diagram 3,
    opacityback=0,  
]
\[\xymatrix{
&D\underset{E}{\times} D - \Delta(D) \ar[ld]^{\pi_1^\gamma} \ar[rd]_{\pi_2^\gamma}& \\
D& &D
}\]
\end{tcolorbox}	

\begin{tcolorbox}[
    standard jigsaw,
    title=Diagram 4,
    opacityback=0,  
]
\[\xymatrix{
&B\underset{C}{\times} B - \Delta(B) \ar[ld]^{\pi_1^\alpha} \ar[rd]_{\pi_2^\alpha}& \\
B& &B
}\]
\end{tcolorbox}	

\begin{tcolorbox}[
    standard jigsaw,
    title=Diagram 5,
    opacityback=0,  
]
\[\xymatrix{
\{(b,d,s:\gamma^{-1}(\gamma(d)) \to C) | \alpha(b) = s(d)\} \ar[d]^{\tilde{\gamma'}} \ar[r]^(.60){\tilde{\alpha}}& (d,s:\gamma^{-1}(\gamma(d)) \to C) \ar[d]^{\gamma'} \\
B \ar[r]^\alpha& C 
}\]
\end{tcolorbox}	

\begin{tcolorbox}[
    standard jigsaw,
    title=Diagram 6,
    opacityback=0,  
]
\[\xymatrix{
\{(b,d,s: \gamma^{-1}(\gamma(d)) \to C) | \alpha(b) = s(d)\} \ar[r]^(.6){\beta' \circ \tilde{\alpha}} & \{(e,s: \gamma^{-1}(e) \to C)\}
}\]
\end{tcolorbox}

\begin{tcolorbox}[
    standard jigsaw,
    title=Diagram 7,
    opacityback=0,  
]
\[\xymatrixrowsep{1in}
\xymatrixcolsep{1in}
\renewcommand{\labelstyle}{\textstyle}
\xymatrix@R=3pc@C=1pc{
 \txt{$\{(b,b',d,d',s)| \alpha(b) = s(d), \alpha(b') =s(d'), $ \\ $\gamma(d) =\gamma(d'), d\neq d' \textrm{ OR  }b \neq b'  \}$ }\ar[d]^{\pi_1^{\beta' \circ \tilde{\alpha}}} \ar[dr]^{\pi_2^{\beta' \circ \tilde{\alpha}}}&\\
\{(b,d,s: \gamma^{-1}(\gamma(d)) \to C) | \alpha(b) = s(d)\} &\{(b',d',s: \gamma^{-1}(\gamma(d)) \to C) | \alpha(b) = s(d)\}
}\]
\end{tcolorbox}	

\begin{tcolorbox}[
    standard jigsaw,
    title=Diagram 8,
    opacityback=0,  
]
\[\xymatrix{
\{(c,d')| \beta(c)\neq d', \gamma \beta(c) = \gamma(d') \} \ar[d]^\lambda \ar[r]^(.65){\tilde{\beta}}& D\underset{E}{\times} D - \Delta(D) \ar[d]^{\pi_1^\gamma}\\
C \ar[r]^\beta& D
}\]
\end{tcolorbox}

\begin{tcolorbox}[
    standard jigsaw,
    title=Diagram 9,
    opacityback=0,  
]
\[\xymatrixrowsep{.3in}
\xymatrixcolsep{.3in}\xymatrix{
D \underset{E}{\times} D - \Delta(D) \ar[d]^{\pi_2^\gamma} & \txt{$\{(c,d')| \beta(c)\neq d', $\\$  \gamma \beta(c) = \gamma(d') \} $}\ar[l]^(.55){\tilde{\beta}}& \txt{$\{(d,d',s:\gamma^{-1}(\gamma(d')) \setminus d' \to C)|$ \\ $d\neq d', \gamma(d) = \gamma(d')\} $}\ar[l]^(.6){\gamma''} \ar[d]^(.6){\beta''}\\
D && \{(d',s:\gamma^{-1}(\gamma(d')) \setminus d' \to C)\} \ar[ll]^{\alpha''}
}\]
\end{tcolorbox}	

\begin{tcolorbox}[
    standard jigsaw,
    title=Diagram 10,
    opacityback=0,  
]
\[\xymatrixrowsep{.2in}
\xymatrixcolsep{.4in}\xymatrix{
D \amalg D  \ar[d]^{\nabla} 
& \{(d',s:\gamma^{-1}(\gamma(d)) \setminus d' \to C)\} \amalg B \ar[l]_(.75){(\alpha'',\beta \circ \alpha)}
& \txt{$2 \cdot \{(b,d',s:\gamma^{-1}(\gamma(d')) \setminus d' \to C)|$ \\  $\beta \alpha(b) = d'\}$} \ar[l]^(.5){(\gamma'''_0 \gamma'''_1)} \ar[d]^(.6){\nabla}\\
D 
&& \txt{$\{(b,d',s:\gamma^{-1}(\gamma(d')) \setminus d' \to C)|$ \\ $ \beta \alpha(b) = d'\}$} \ar[ll]^{\alpha'''}
}\]
\end{tcolorbox}	

\begin{tcolorbox}[
    standard jigsaw,
    title=Diagram 11,
    opacityback=0,  
]
\[\xymatrixrowsep{.2in}
\xymatrixcolsep{.4in}\xymatrix{
\txt{$ \{(b,b',d', s:\gamma^{-1}(\gamma(d')) \setminus d' \to C) | \beta \alpha (b) = d',$ \\ $ b\neq b', \alpha(b) = \alpha(b')\}$} \ar[r]^(.7){\tilde{\gamma}_1'''} \ar[d]^(.55)\phi& B \underset{C}\times B - \Delta B  \ar[d]^{\pi_2^\alpha}\\
\{(b,d', s:\gamma^{-1}(\gamma(d')) \setminus d' \to C) | \beta \alpha (b) = d'\} \ar[r]^(.8){\gamma_1'''}& B
}\]
\end{tcolorbox}	

\begin{tcolorbox}[
    standard jigsaw,
    title=Diagram 12,
    opacityback=0,  
]
\[\xymatrixrowsep{.05in}
\xymatrixcolsep{.05in}
\renewcommand{\labelstyle}{\textstyle}
\xymatrix@R=3pc@C=1pc{
B  \ar[d]^{\alpha}& \txt{$\{(b,d,d',s) |$ \\ $d \neq d', \gamma (d) = \gamma(d'), \alpha(b) = s(d)\}$} \ar[l] \ar[d]
& \txt{$\{(b_0,b_1,d,d',s )| d \neq d', $ \\$\gamma (d) = \gamma(d'), \alpha(b_0) = s(d), \beta\alpha(b_1) = d' \}$} \ar[l] \ar[d] \ar@/^-4pc/[ll]^{\delta}  \ar@/^10pc/[dd]^{\epsilon}\\
C  &\txt{ $\{(d,d',s )|$ \\$ d \neq d', \gamma (d) = \gamma(d') \}$} \ar[l]^(.75){\lambda \circ \gamma''} \ar[d]^{\beta''}
&\txt{$\{(b,d,d',s)| $\\$d\neq d', \gamma(d) = \gamma(d'), \beta \alpha (b) = d'\}$}  \ar[l] \ar[d]&\\
   &\{(d',s:\gamma^{-1}(\gamma(d')) \setminus d' \to C) \} &\txt{$\{(b,d',s:\gamma^{-1}(\gamma(d')) \setminus d' \to C) |$ \\$ \beta \alpha (b) = d' \}$}\ar[l]^{\gamma_0'''}
}\]
\end{tcolorbox}

We can show the desired equality by showing the following:
\begin{eqnarray*}
T_{\alpha'} T_{\beta' \circ \tilde{\alpha}} &=& T_\gamma T_{\alpha'''} \\
R_{\tilde{\gamma}'} &=& R_{\gamma_1'''}\\
N_{\pi_2}^{\beta' \circ \tilde{\alpha}} R_{\pi_1}^{\beta' \circ \tilde{\alpha}} R_{\tilde{\gamma'}} R_k &=& N_{\epsilon} R_\delta R_k \cdot N_\phi R_{\tilde{\gamma}'''_1} R_{\pi_1}^\alpha R_k \\
R_{\gamma'''_0}N_{\beta''}R_{\gamma''}R_\lambda(r) \cdot R_{\gamma_1'''}R_\alpha(r) &=& R_{\beta' \circ \tilde{\alpha}} N_{\beta'}R_{\gamma'}(r) .
\end{eqnarray*}
The first one is true because the following sets are isomorphic:
$$\{(b,d,s:\gamma^{-1}(\gamma(d)) \to C) | \alpha(b) = s(d) \} \cong \{(b,d,s:\gamma^{-1}(\gamma(d)) \setminus d \to C) | \beta \alpha(b) = d \} $$ because the rightmost set has selected an element to map $s(d)$ to, the element $\alpha(b)$.

The second one is true for an identical reason. The third one is true because the set $$  \{ (b_0,b_1,d,d',s:\gamma^{-1}(\gamma(d')) \setminus d' \to C) | d \neq d', \gamma (d) = \gamma(d'), \alpha(b_0) = s(d), \beta\alpha(b_1) = d' \} $$ $$\coprod \{(b,b',d', s:\gamma^{-1}(\gamma(d')) \setminus d' \to C) | \beta \alpha (b) = d', b\neq b', \alpha(b) = \alpha(b')\}$$ is isomorphic to 
$$\{ (b,b',d,d',s: \gamma^{-1}(\gamma(d)) \to C) | \alpha(b) = s(d), \alpha(b') =s(d'),\\ \gamma(d) =\gamma(d'), d\neq d' \textrm{ OR  }b \neq b'  \}$$ by essentially the same argument as the first equation. The first set in the disjoint union is the case $d \neq d'$, then second is the case $d = d'.$ The fourth is true for an almost identical reason. Therefore, we have shown that the above defines a genuine derivation, and we have the isomorphism.
\end{proof}

\begin{cor}
The module of K\"ahler differentials for  $\underline{R} \to \underline{S} = \underline{R} \boxtimes \underline{A}^{\mathcal{O}}[x_H]$ is the free $\underline{S}$-module generated by a single element at $G/H$.
\end{cor}

\begin{cor}
The module of K\"ahler differentials for $\underline{R} \to \underline{S} = \underline{R} \boxtimes (\underset{i} \boxtimes \underline{A}^{\mathcal{O}}[x_{H_i}])$ is the free $\underline{S}$-module generated by an element at every $G/H_i$ when $i$ is finite.
\end{cor}
\begin{proof}
This follows from ~\cref{lemprod}.
\end{proof}

%% file: CotangentComplex.tex
\section{Cotangent Complex, (co)homology}

Using the Mackeyization functor, we will now describe the cotangent complex, the $q$th Andr\'e-Quillen (co)homology Mackey functors of the $\underline{A}$-Tambara functor $\underline{B}$ with coefficients in the $\underline{B}$ module $\underline{M}$, denoted $D^q(\underline{B}|\underline{A}; \underline{M})$ and $D_q(\underline{B}|\underline{A}; \underline{M}).$

Given any simplicial object in an abelian category (a contravariant functor from $\Delta^{op} \to \cat{C}$ our abelian category), we have the homology in the $q$th dimension $H^q(X)$, the $q$th homology of the chain complex with $d = \sum_{i} (-1)^i d_i$. Via the map in the Dold-Kan correspondence, this is the same as the homology of the normalized subcomplex, or phrased differently, the $q$th homotopy group. The homology and cohomology will be built via cofibrant replacement, so we briefly discuss the simplicial model structure on the category of simplicial $\underline{A}$-algebras over $\underline{B}$, following \cite{goerss2006model}, which is nearly identical to the classical case.

Chain complexes of modules over an incomplete Tambara functor $\underline{A}$ has a model category structure:
\begin{prop}\label{model}
The category $\cat{Ch_* \underline{A}}$ has the structure of a model category where $f:\underline{M}_* \to \underline{N}_*$ is 
\begin{itemize}
\item a weak equivalence if $H_* f$ is an isomorphism,
\item a fibration if $\underline{M}_n \to \underline{N}_n$ is surjective for $n \geq 1$, and
\item a cofibration if for $n \geq 0$, the map $\underline{M}_n \to \underline{N}_n$ is an injection with projective cokernel.
\end{itemize}
\end{prop}

\begin{proof}
This is similar to \cite{goerss2006model} Theorem 1.5.
We will show only the instances where the proofs differ from the classical case described in this source.
The axioms of limits and colimits, retracts, 2-out-of-3, and lifting to an acyclic fibration are all identical.

Now we show the acyclic-cofibration to a fibration axiom. The difference here is the presence of more than one $n$-disc. Let $\underline{D}(n,H)$ for $n\geq 1$ be the chain complex with $\underline{D}(n,H)_k = 0$ for $k \neq n, n -1$ and the free $\underline{A}$-module with one generator at $G/H$ otherwise. The only relevant differential is the identity. Then there is a natural isomorphism $\cat{Ch_* \underline{A}}(\underline{D}(n,H),\underline{N}_*) \cong \underline{N}_n(H).$ Therefore, $q : \underline{Q}_* \to \underline{N}_*$ is a fibration if and only if $q$ has the right lifting property with respected to $0 \to \underline{D}(n,H)$ for all $n>0, H < G.$

If $\underline{N}_*$ is any chain complex, define a new chain complex $P(\underline{N}_*)$ and an evaluation map $\epsilon: P(\underline{N}_*) \to \underline{N}_*$ by $$P(\underline{N}_*) = \bigoplus_{n>0, H < G} \bigoplus_{x \in \underline{N}_n(G/H)} \underline{D}(n,H) \to \underline{N}_*.$$ The map is clearly a fibration, but more so for any $\underline{M}_* \to \underline{N}_*$, we can factor this map by $$\underline{M}_* \to \underline{M}_* \oplus P(\underline{N}_*) \to \underline{N}_*$$ where the first map is an acyclic cofibration, completing one half of the factorization axioms. The rest of the proof is identical as in  \cite{goerss2006model}.

\end{proof}

Note from the proof above that chain complexes are cofibrantly generated. The set that generates the acyclic cofibrations are the maps $0 \to \underline{D}(n,H)$. The set that generates the cofibrations are maps of the form $0 \to \underline{D}(n,H)$, $\underline{S}(n-1,H) \to \underline{D}(n,H)$, and $0 \to \underline{S}(0,H)$, where $\underline{S}(n,H)$ is the chain complex that has a single non-zero module at dimension $n$, and is the free module generated by a single element at $G/H$.

\begin{prop}
$\cat{Ch_* \underline{A}}$ is a cofibrantly generated model category, with the generating sets for cofibrations and acyclic cofibrations as described above.
\end{prop}

\begin{proof}
Firstly, any bounded chain complex of finitely generated $\underline{A}$-modules is small for the set of all morphisms. The map $\colim C(\underline{A},\underline{X}_n) \to C(\underline{A},\colim \underline{X}_n)$ is always injective: if we have two maps in $\colim C(\underline{A},\underline{X}_n) $, then there exists some $n$ such that we can realize both maps in $C(\underline{A},\underline{X}_n) $. If they are unequal as maps, they give unequal maps in the colimit. The map is surjective because if the module is bounded and finitely generated, then for every map in $ C(\underline{A},\colim \underline{X}_n)$ we can realize this map as a map in $C(\underline{A},\underline{X}_n)$ for a specific $n$.

The property that a map is an acyclic fibration if has the RLP with respect to the generating cofibrations, and similarly for the fibrations and generating acyclic cofibrations, follows from \cref{model}.

\end{proof}

Via the Dold-Kan correspondence, which can be taken to be an inverse equivalence, we can lift a model structure to $\cat{sMod}_{\underline{A}}$, the simplicial modules over $\underline{A}$. The two items to check are (1) the normalization functor $N$ commutes with sequential colimits and (2) any cofibration in $\cat{sMod}_{\underline{A}}$ (maps which have the left lifting property with respect to acyclic fibrations) that has the left lifting property with respect to all fibrations is a weak equivalence. Any equivalence of categories satisfies both these conditions. So by \cite{goerss2006model} Theorem 3.6, we have 

\begin{prop}
The category $\cat{sMod}_{\underline{A}}$ has the structure of a model category where a morphism $f: \underline{X} \to \underline{Y}$ is 
\begin{itemize}
\item a weak equivalence if $\pi_* \underline{X} \to \pi_* \underline{Y}$ is an isomorphism; and
\item a fibration if $N \underline{X}_n \to N \underline{Y}_n$ is onto for $n \geq 1.$
\end{itemize}
The cofibrations are generated by $\Gamma [0 \to \underline{D}(n,H)]$ and $\Gamma [\underline{S}(n-1,H) \to \underline{D}(n,H)]$ where $\Gamma$ is the inverse to the normalization functor in the Dold-Kan correspondence.
\end{prop}

We now want to build a model category structure on  $s\underline{A}\cat{-CAlg}$. This can be done again by lifting a model category structure via the adjunction $$S_{\underline{A}} :  \cat{sMod}_{\underline{A}} \leftrightarrow s\underline{A}\cat{-CAlg}:U$$ where $S_{\underline{A}} $ is the symmetric algebra functor on each dimension, creating the free $\underline{A}$-algebra generated by the elements of the $\underline{A}$-module, subject to the $\underline{A}$-module relations. This is the standard free-forget adjunction. The proof follows as in \cite{goerss2006model} Theorem 4.17. It is clear the forgetful functor commutes with sequential colimits, as the colimit in both is given by the colimit as simplicial sets. Now we need to appeal to  \cite{goerss2006model} Theorem 3.8. All elements of $s\underline{A}\cat{-CAlg}$ are fibrant, and there is a natural path object for every element $\underline{B} \in s\underline{A}\cat{-CAlg}$ given by $\underline{B}^{\Delta^1}$, which exists because it is a simplicial category. Therefore, we get the following model category structure:

\begin{prop}
The category $s\underline{A}\cat{-CAlg}$ has the structure of a model category where a morphism $f: \underline{X} \to \underline{Y}$ is 
\begin{itemize}
\item a weak equivalence if $\pi_* \underline{X} \to \pi_* \underline{Y}$ is an isomorphism; and
\item a fibration if $N \underline{X}_n \to N \underline{Y}_n$ is onto for $n \geq 1$, where we consider $\underline{X}_n,\underline{Y}_n$ as simplicial modules.
\end{itemize}
The cofibrations are generated by $S \Gamma [0 \to \underline{D}(n,H)]$ and $S \Gamma [\underline{S}(n-1,H) \to \underline{D}(n,H)]$ where $\Gamma$ is the inverse to the normalization functor in the Dold-Kan correspondence and $S$ is the free symmetric algebra functor.
\end{prop}

Again we can classify cofibrations in a more direct way. We call a map $\underline{R}_\cdot \to \underline{S}_\cdot$ in $s\underline{A}\cat{-CAlg}$ free if there exists GIndexedSets $\underline{C}_*$, which when $\eta$ is a surjective map from $[p] \to [q]$ we have $\eta^* \underline{C}_q \subset \underline{C}_p$ such that $\underline{S}_q$ is a free $\underline{R}_q$-algebra generated by $\underline{C}_q$. Quillen shows (in the non-equivariant setting though the equivariant setting goes through just as formally) in \cite{quillen2006homotopical} that any free map is a cofibration, any map in $s\underline{A}\cat{-CAlg}$ can be factored in to a free map followed by an acyclic fibration, and that a map is a cofibration if and only if it is a retract of a free map.

This immediately gives us a model category structure on the over category $\\ s\underline{A}\cat{-CAlg}/ \underline{B}.$ Just as in the classical case, we define the Quillen homology to be the total left derived functor of Mackeyization. Namely, there are constant simplicial algebras $c\underline{A}$ and $c\underline{B}.$ Given an element $\underline{R} \in \underline{A}\cat{-CAlg}/ \underline{B}$, we can take the cofibrant replacement of the $c\underline{A} \to c\underline{B} $ to yield $c\underline{A} \to \underline{P} \to c\underline{B}.$ Applying the Mackeyization functor dimension-wise yields $$\underline{P}_n \mapsto \Omega_{\underline{P}_n/\underline{A}}^{1,G,\mathit{O}} \underset{\underline{R}}\boxtimes \underline{B}$$ which we call the \term{cotangent complex} $L_{\underline{B}/\underline{A}}$ or $L_\phi$ if $\phi: \underline{A} \to \underline{B}$. Note that by standard arguments of the cofibrant replacement, the cotangent complex is well defined up to chain homotopy equivalence. We can also always assume that $\underline{P}$ is built out of free $\mathcal{O}$-Tambara functors. Then we can define the \term{Andr\'e-Quillen homology and cohomology Mackey Functors} as:

$$D_q(\underline{B}|\underline{A};\underline{M}) = \pi_q(L_{\underline{B}|\underline{A}} \underset{\underline{B}} \boxtimes \underline{M}) \cong H_q(NL_{\underline{B}/\underline{A}} \boxtimes \underline{M})$$

$$D^q(\underline{B}|\underline{A};\underline{M}) = H^q(\underline{\Hom}_{\underline{B}}(NL_{\underline{B}|\underline{A}} , \underline{M})) .$$
If $\underline{M} = \underline{B}$, we write $D_q(\underline{B}|\underline{A}),D^q(\underline{B}|\underline{A}).$

%% file: CotangentComplexProp.tex
\subsection{Properties}
Let us now state some elementary properties of the homology and cohomology. A reader who is interested in categories similar to $\cat{\underline{R}Mod}$ where projective modules act more strangely (box products not preserving projectivity, projective not being flat, etc.) should take particular care in this section.
The properties are generalizations of properties from \cite{goerss2006model},\cite{iyengar2007andre}, and \cite{QuillenMimeo}.

To begin, when $\underline{B}$ is a free $\underline{A}$-algebra, then the identity map is a cofibrant replacement, because the map $c\underline{A} \to c\underline{B}$ is a free map. In this case, all the non-zero homology and cohomology vanishes, and at degree 0:
$$D^0 (\underline{B}|\underline{A};M) \cong \Hom_{\underline{R}-Mod}( \Omega_{\underline{B}/\underline{A}},\underline{M}) \cong \Der_{\underline{A},\mathcal{O}}(\underline{B},\underline{M})$$ 
$$D_0 (\underline{B}|\underline{A};M) \cong  \Omega_{\underline{B}/\underline{A}} \underset{\underline{B}}\boxtimes \underline{M}.$$

\begin{definition}
A simplicial module $\underline{X}_*$ over a simplicial Tambara-functor $\underline{A}_*$ is a \term{free simplicial module over $\underline{A}_*$} if there are $G$ indexed subsets $\underline{C}_q \subset \underline{X}_q$ such that $\eta^* \underline{C}_q \subset \underline{C}_p$ if $\eta:[p] \to [q]$ is a surjective monotone map, and $\underline{X}_q$ is the free $\underline{A}_q$ module generated by the elements $\underline{C}_q$. A \term{projective simplicial module} is a direct summand of a free simplicial module.
\end{definition}

Suppose $\underline{P}_*$ is a free $\underline{A}$ algebra resolution of $\underline{B}$ with generators $\underline{C}_*$ at each level. Then $ \Omega_{\underline{B}/\underline{A}} \underset{\underline{P}} \boxtimes \underline{B}$ is generated as a free simplicial module by $\{dx \otimes 1 | x \textrm{ an element of } C_q\}$. If $\underline{Q}_*$ is a projective resolution of $\underline{B}$, then it is a retract of a free simplicial resolution. Therefore, the cotangent complex is a summand of a free simplicial module, and therefore is a projective simplicial module.

From this fact, we can construct two K\"unneth spectral sequences. Suppose we have an $\underline{B}$-module $\underline{M}$. We can take a cofibrant replacement of $c\underline{M}$ in $s \cat{\underline{B}Mod}$, call it $\underline{Q}_*$. By \cite{goerss2006model} 4.4 and 4.5, $\underline{Q}_p$ is a projective $\underline{B}$-module for all $p \geq 0$ and has the same homotopy groups as $c\underline{M}$. Then $\underline{Q}_* \boxtimes L_{\underline{B}|\underline{A}}$ forms a bisimplicial object in the category of $\cat{\underline{B}Mod}$. Using the Eilenberg-Zilber theorem, we get two convergent first quadrant spectral sequences $$E_{pq}^1 = \pi_q^v(\underline{Q}_p \boxtimes (L_{\underline{B}|\underline{A}})_*), \quad E_{pq}^2 = \pi_p^h \pi_q^v(\underline{Q}_* \boxtimes (L_{\underline{B}|\underline{A}})_*) \Longrightarrow \pi_{p+q}(\underline{Q}_* \boxtimes (L_{\underline{B}|\underline{A}})_*)$$

$$E_{pq}^1 = \pi_q^h(\underline{Q}_* \boxtimes (L_{\underline{B}|\underline{A}})_p), \quad E_{pq}^2 = \pi_p^v \pi_q^h(\underline{Q}_* \boxtimes (L_{\underline{B}|\underline{A}})_*) \Longrightarrow \pi_{p+q}(\underline{Q}_* \boxtimes (L_{\underline{B}|\underline{A}})_*)$$

Because $L_{\underline{B}|\underline{A}}$ projective and therefore flat, the second spectral sequence has $E_{pq}^2 = \pi_p^v \pi_q^h(\underline{Q}_* \boxtimes (L_{\underline{B}|\underline{A}})_*) = \pi_p^v (\underline{M} \boxtimes (L_{\underline{B}|\underline{A}})_*)  $ if $q = 0$ and otherwise is zero. The spectral sequence collapses at the second page, and converges to $\pi_{p+q} (\underline{M} \boxtimes L_{\underline{B}|\underline{A}}) = D_{p+q} (\underline{B}|\underline{A},\underline{M}).$ 

Because the terms $\underline{Q}_*$ are flat, the $E_{pq}^2$ terms of the first spectral sequences are $E_{pq}^2 = \pi_p^h \pi_q^v(\underline{Q}_* \boxtimes (L_{\underline{B}|\underline{A}})_*) = \pi_p^h (\underline{Q}_* \boxtimes D_q({\underline{B}|\underline{A}})) = \underline{\Tor}_p^{\underline{B}}(D_q({\underline{B}|\underline{A}}), \underline{M}).$ This, and its dual cohomological version with Ext, gives the following proposition.

\begin{prop}\label{basicss}
We have two convergent spectral sequences 
$$E_{pq}^2 = \underline{\Tor}_p^{\underline{B}}(D_q(\underline{B}|\underline{A}),\underline{M}) \Longrightarrow D_{p+q}(\underline{B}|\underline{A},\underline{M}).$$

$$E^{pq}_2 = \underline{\Ext}^p_{\underline{B}}(D_q(\underline{B}|\underline{A}),\underline{M}) \Longrightarrow D^{p+q}(\underline{B}|\underline{A},\underline{M}).$$
\end{prop}

Note that $D_0(\underline{B}|\underline{A}) = \Omega_{\underline{B}|\underline{A}}$ since the Mackey functor objectivization is a left adjoint. By the two spectral sequences above in the 0th term, which is determined by the second page, we have isomorphisms $\Omega_{\underline{B}|\underline{A}} \underset{\underline{B}} \boxtimes \underline{M} \cong D_0(\underline{B}|\underline{A},\underline{M})$ and $D_0(\underline{B}|\underline{A};\underline{M}) \cong \underline{\Hom}_{\underline{B}}(\Omega_{\underline{B}|\underline{A}}, \underline{M}) \cong \underline{\Der}(\underline{B}|\underline{A},\underline{M}).$

We can define what an $\underline{A}$-algebra extension of $\underline{B}$ by $\underline{M}$. We mean a short exact sequence $$ 0 \to \underline{M}_Y \overset{i}\to \underline{X} \overset{u} \to \underline{B} \to 0$$ where $Y$ is a $G$-set, $u$ is a map of $\underline{A}$ algebras such that $\ker(u)$ is killed by admissible norms of all 2-surjective maps, and $\underline{M} \cong \ker u$ as $\underline{B}$-modules, with the kernel having a $\underline{B}$-module structure via $u(x) \cdot y:= x \cdot y$, which is well defined since $(\ker u)^2 = 0$. 
Let $\Exalcomm(\underline{B}|\underline{A};\underline{M})$ be the $G$ Indexed set of isomorphism classes of extensions, indexed by $Y$. 

\begin{prop}
$D^1(\underline{B}|\underline{A};\underline{M}) \cong \Exalcomm(\underline{B}|\underline{A};\underline{M})$
\end{prop}

\begin{proof}
Let $\underline{P}_*$ be a free $\underline{A}$-algebra  resolution of $\underline{B}$.
Given an extension $ 0 \to \underline{M}_Y \overset{i}\to \underline{X} \overset{u} \to \underline{B} \to 0$, we can construct a map in $\underline{A}\cat{-CAlg}/\underline{B}$ from $\theta: \underline{P}_0 \to \underline{X}$ by specifying where the generators map. This gives a map $\theta \circ (d_0 - d_1): \underline{P}_1 \to \underline{X}$, a Mackey functor map which sends $\underline{A}$ to 0 and maps to the kernel of $u$. Therefore, it induces a map $i^{-1} (\theta \circ (d_0 - d_1)): \underline{P}_1 \to \underline{M}_Y$ of Mackey functors such that $\underline{A}$ maps to 0. We now show that this is a derivation.

Suppose $f$ is a fold map, we show that the products turn to sums:
\begin{eqnarray*}
(\theta \circ (d_0 - d_1))(xy) &=& \theta d_0 x \cdot \theta d_0 y - \theta d_1 x \cdot \theta d_1 y \\
&=& (\theta d_0 y) (\theta d_0 x - \theta d_1 x ) + (\theta d_1 x)(\theta d_0 y - \theta d_1 y) \\
&=& y \cdot (\theta d_0 x - \theta d_1 x ) + x \cdot(\theta d_0 y - \theta d_1 y). \\
\end{eqnarray*}
with the last line using the $\underline{P}_1$-module structure on $\ker u.$

If it is an equivariant map $f: G/H \to G/K$, then we can use the same computation as in ~\cref{kahler}:

\begin{eqnarray*}
\theta d_0(N_f(x)) &=& N_f(\theta d_0(x)) \\
&=& N_f[\theta d_1(x) + (\theta d_0(x) - \theta d_1(x))] \\
&=& N_f[\theta d_1(x)] + T_f[ (N_{\pi_2}R_{\pi_1} \theta d_1 x) \cdot (\theta d_0(x)- \theta d_1(x))]  \\
&=& N_f[\theta d_1(x)] + T_f[ (N_{\pi_2}R_{\pi_1} x) \cdot (\theta d_0(x)-\theta d_1(x))].
\end{eqnarray*}
The third equation uses the fact that the kernel is 0 under any 2-surjective map, so there are no higher order terms involving the norm of an element of the kernel. The last line again uses the $\underline{P}_1$-module structure on the kernel. Therefore, we obtain:
$$\theta(d_0 - d_1)(N_f(x)) = T_f[ (N_{\pi_2}R_{\pi_1} x) \cdot (\theta (d_0-d_1) x)]$$ verifying that it is a derivation.
Note from the construction that the map from $\underline{P_1} \to \underline{M}_Y$ induces the 0 map $\underline{P_2} \to \underline{M}$ via the composition with $d_0 - d_1 + d_2.$ Additionally, if we change the map $\theta$ to $\theta'$, we change the derivation obtained, but only up to a derivation induced from $\underline{P}_0 \to \underline{M}_Y$. By the exact same computation as above, $\theta - \theta'$ gives a derivation from $\underline{P}_0 \to \underline{M}_Y$, which induces the map $(\theta - \theta')(d_0 - d_1): \underline{P}_1 \to \underline{M}_Y.$ To sum up, we have defined a map $\Phi: \Exalcomm(\underline{B}|\underline{A};\underline{M}) \to D^1(\underline{B}|\underline{A};\underline{M}).$

We now construct an inverse. Suppose we have a derivation $D: \underline{P}_1 \to \underline{M}.$ Let $$\underline{X} := \coker\{ \underline{P}_1 \overset{(d_0-d_1,D)} \longrightarrow \underline{P}_0 \ltimes \underline{M}_Y\}$$ taken in the category of $\underline{A}$-algebras, which comes with a canonical projection $p: \underline{P}_0 \ltimes \underline{M}_Y \to \underline{X}$. We naturally have a map $\underline{P}_0 \ltimes \underline{M}_Y \to \underline{B}$ from the map $\underline{P}_0 \to \underline{B}$ that surjects and induces a map $\underline{X} \to \underline{B}$. The kernel of this map is precisely $\underline{M}_Y \to \underline{P_0} \ltimes \underline{M}_Y$, which gives a sequence $0 \to \underline{M}_Y \to \underline{X} \to \underline{B} \to \underline{0}$, where all that remains to show is that $\underline{M}_Y \to \underline{X}$ is injective. If $(d_0 - d_1)(p_1) = 0$, then it is in the image of an element $(d_0 - d_1 + d_2)(p_2)$, so then $Dp_1 = 0$ via the cocycle condition of the derivation. Suppose we change $D$ to $D + (d_0 - d_1) \circ D'$ where $D': \underline{P}_0 \to \underline{M}$ is a derivation. Then we can construct isomorphisms between the cokernels 
\[\xymatrixrowsep{.5in}
\xymatrixcolsep{1.5in}\xymatrix{
\underline{X_0} :=  \coker\{ \underline{P}_1\ar[r]^(.55){(d_0-d_1,D)}& \underline{P}_0 \ltimes \underline{M}_Y\}\\
\underline{X_1} :=   \coker\{ \underline{P}_1 \ar[r]^(.55){(d_0-d_1,D + (d_0 - d_1) \circ D' )} &  \underline{P}_0 \ltimes \underline{M}_Y\}
}\]
via maps $\underline{P}_0 \ltimes \underline{M}_Y \to \underline{P}_0 \ltimes \underline{M}_Y$ such that $(p_0,m)\mapsto (p_0, m + D'(p_0))$ and $(p_0,m)\mapsto (p_0, m - D'(p_0))$. This is a well defined map $\Phi^{-1}: D^1(\underline{B}|\underline{A};\underline{M}) \to \Exalcomm(\underline{B}|\underline{A};\underline{M}).$

To show that these are mutual inverses, we see that if we have a derivation $D: \underline{P_1} \to \underline{M}_Y$ and construct $\underline{X}$ as the cokernel above, there is naturally a map $\underline{P_0} \overset{(\id,0)}\to \underline{P}_0 \ltimes \underline{M}_Y \to \underline{X}.$ The map induces $\underline{P_1} \overset{(d_0-d_1)} \to \underline{P_0} \overset{(\id,0)}\to \underline{P}_0 \ltimes \underline{M}_Y \to \underline{X}$ which by the property of the cokernel is the same as the map $\underline{P_1} \overset{(0, D)} \to \underline{P}_0 \ltimes \underline{M}_Y$, so we have that $\Phi \circ \Phi^{-1}$ is the identity on $D^1(\underline{B}|\underline{A};\underline{M}) $

Suppose we have an extension
\[\xymatrixrowsep{.5in}
\xymatrixcolsep{.5in}\xymatrix{
0 \ar[r]& \underline{M}_Y \ar[r] & \underline{X} \ar[r] & \underline{B} \ar[r] & 0
}\]
and we take $\theta: \underline{P}_0 \to \underline{X}$, the induced derivation $D: \underline{P}_1 \to \underline{M}_Y$, and the extension constructed by the cokernel $\underline{X'}$.
\[\xymatrixrowsep{.5in}
\xymatrixcolsep{.5in}\xymatrix{
0 \ar[r]& \underline{M}_Y \ar[r] & \underline{X'} \ar[r] & \underline{B} \ar[r] & 0\\
0 \ar[r]& \underline{M}_Y \ar[r] & \underline{X} \ar[r] & \underline{B} \ar[r] & 0
}\]

Let $\underline{P}_0 \ltimes \underline{M}_Y \to \underline{X}$ be the map $(p,m) \mapsto \theta(p) + i(m).$ This map is a map of Mackey functors and by the same computation as ~\cref{kahler}, one can show it commutes with norms. Therefore, it is a map of $\underline{A}$-algebras. Composing this map with the map $\underline{P}_1 \to \underline{P}_0 \ltimes \underline{M}_Y$ is $\theta \circ (d_0 - d_1) p_1 + i(D(p_1)) = 0$ by the construction of $D$. So it induces a map $\underline{X'} \to \underline{X}$. The map commutes with the identity maps between $\underline{M}_Y$ and $\underline{B}$ so we get the commutative diagram between extensions
\[\xymatrixrowsep{.5in}
\xymatrixcolsep{.5in}\xymatrix{
0 \ar[r] \ar[d]& \underline{M}_Y \ar[r] \ar[d]& \underline{X'} \ar[r] \ar[d]& \underline{B} \ar[r] \ar[d]& 0 \ar[d]\\
0 \ar[r]& \underline{M}_Y \ar[r] & \underline{X} \ar[r] & \underline{B} \ar[r] & 0
}\]
which is all that is necessary to show that $\underline{X'}$ and $\underline{X}$ are isomorphic as extensions. This shows that $\Phi^{-1} \circ \Phi$ is the identity on $\Exalcomm(\underline{B}|\underline{A};\underline{M})$ completing the proof.

\end{proof}

\begin{prop}
Suppose $\underline{A} \to \underline{B}$ is a surjection of $\mathcal{O}$-Tambara functors. The kernel $\underline{I}$ is an ideal of $\underline{A}$. Then $D_0(\underline{B}|\underline{A}) = 0$ and $D_1(\underline{B} | \underline{A}) = \underline{I}/\underline{I}^{>1}.$
\end{prop}

\begin{proof}
$D_0(\underline{B}|\underline{A}) = \Omega_{\underline{B}|\underline{A}} = 0$.  From the second spectral sequence of ~\cref{basicss}, and the fact that $D_0(\underline{B}|\underline{A}) = 0$, the five term exact sequence gives an isomorphism between $$\underline{\Hom}_{\underline{B}} (D_1(\underline{B}|\underline{A}), \underline{M}) \cong D^1(\underline{B}|\underline{A},\underline{M})$$ which is by the the above proposition $\Exalcomm(\underline{B}|\underline{A};\underline{M}).$

The isomorphism class $0 \to \underline{I}/\underline{I}^{>1} \to \underline{A}/\underline{I}^{>1} \to \underline{B} \to 0$ defines an element of $\xi \in \Exalcomm(\underline{B}|\underline{A}, \underline{I}/\underline{I}^{>1}).$ Functoriality of Exalcomm in the module coordinate gives a map $\xi_*: \underline{\Hom}_{\underline{B}}(\underline{I}/\underline{I}^{>1},\underline{M}) \to \Exalcomm(\underline{B}|\underline{A};\underline{M}).$ If we have an extension $0 \to \underline{M}_Y \to \underline{X} \to \underline{B} \to 0$, then we have a map $\underline{A} \to \underline{X}$ since $\underline{X}$ is an $\underline{A}$-algebra. This induces a map $\underline{I}/\underline{I}^{>1} \to \underline{M}_Y$, which is an inverse to $\xi_*.$ So $\underline{\Hom}_{\underline{B}}(D_1(\underline{B}|\underline{A}),\underline{M}) \cong \underline{\Hom}_{\underline{B}}(\underline{I}/\underline{I}^{>1}, \underline{M})$ which by the enriched Yoneda Lemma gives $D_1(\underline{B}|\underline{A}) \cong \underline{I}/\underline{I}^{>1}.$
\end{proof}

We have a functoriality property in the coefficients. In particular, if we have a map of $\underline{S}$-modules $\underline{M} \to \underline{M'}$, then for $\phi: \underline{R} \to \underline{S}$ we have an induced map of $L_\phi \underset{\underline{S}} \boxtimes \underline{M} \to L_\phi \underset{\underline{S}} \boxtimes \underline{M'}$ which, when taking the homology and cohomology, is functorial. Additionally, because $L_\phi \underset{\underline{S}} \boxtimes (-)$ is exact, since $L_\phi$ is projective as a simplicial $\underline{S}$-module, given an exact sequence of
\[\xymatrixrowsep{.5in}
\xymatrixcolsep{.5in}\xymatrix{
0 \ar[r] & \underline{M'} \ar[r] & \underline{M} \ar[r]& \underline{M''} \ar[r] & 0
}\]
 we have an exact sequence of simplicial modules 
 \[\xymatrixrowsep{.5in}
\xymatrixcolsep{.5in}\xymatrix{
0 \ar[r] & L_\phi \underset{\underline{S}} \boxtimes \underline{M'} \ar[r] &  L_\phi \underset{\underline{S}} \boxtimes\underline{M} \ar[r]& L_\phi \underset{\underline{S}} \boxtimes \underline{M''} \ar[r] & 0
}\]
which taking the (co)homology gives two exact sequences:
 \[\xymatrixrowsep{.5in}
\xymatrixcolsep{.25in}\xymatrix{
0 \ar[r] & D^0(\underline{S}|\underline{R};\underline{M'}) \ar[r] & D^0(\underline{S}|\underline{R};\underline{M}) \ar[r] & D^0(\underline{S}|\underline{R};\underline{M''}) \ar[r] & D^1(\underline{S}|\underline{R};\underline{M'}) \ar[r] & \cdots \\
\cdots \ar[r] & D_1(\underline{S}|\underline{R};\underline{M''}) \ar[r] & D_0(\underline{S}|\underline{R};\underline{M'}) \ar[r] & D_0 (\underline{S}|\underline{R};\underline{M}) \ar[r] & D_0(\underline{S}|\underline{R};\underline{M''}) \ar[r] & 0
}\]

Now suppose that we have two maps of $\mathcal{O}$-Tambara functors that fit together into the following commutative diagram

\[\xymatrixrowsep{.5in}
\xymatrixcolsep{.5in}\xymatrix{
\underline{R'} \ar[r] \ar[d]& \underline{S'} \ar[d] \\
\underline{R} \ar[r]& \underline{S}
}\]
that give us the cotangent complexes $L_{\underline{S}|\underline{R}}$ and $L_{\underline{S'}|\underline{R'}}$ which are respectively $\underline{S}$ and $\underline{S'}$-modules. We might suspect that a map on on the relative sets of Tambara functors induces a map on their homology, which it does:

\begin{prop}
Suppose we have a commutative diagram as above. Then we get a morphism of the form $\underline{R} \underset{\underline{R'}} \boxtimes L_{\underline{S'}|\underline{R'}} \to L_{\underline{S}|\underline{R}}$ which is a map of $\underline{R'} \underset{\underline{R}} \boxtimes \underline{S}$-modules.
\end{prop}

\begin{proof}
Suppose we take the cofibrant replacement $\underline{A'}$ of $\underline{S'}$ as a simplicial $\underline{R'}$-algebra. We can assume that $\underline{A'}_*$ is free as an $\underline{R'}$-algebra. Since $\underline{A'} \underset{\underline{R'}} \boxtimes \underline{R}$ is the pushout of the diagram made by $\underline{A'},\underline{R'},\underline{R}$, we get a map of $\underline{R}$-algebras $$\underline{A'} \underset{\underline{R'}} \boxtimes \underline{R} \to \underline{S'} \underset{\underline{R'}} \boxtimes \underline{R}.$$ Since $\underline{A'}$ is a free resolution of $\underline{R'}$, $\underline{A'} \underset{\underline{R'}} \boxtimes \underline{R}$ is free over $\underline{R}$ and so is a cofibration in $\underline{R}$-algebras. Let $\underline{A}$ be the cofibrant replacement of $\underline{S}$ as an $\underline{R}$-algebra. $\underline{A} \to \underline{S}$ is an acyclic fibration in $\underline{R}$-algebras so the commutative diagram 
\[\xymatrixrowsep{.5in}
\xymatrixcolsep{.5in}\xymatrix{
\underline{R} \ar[r] \ar[d]& \underline{A} \ar[d] \\
\underline{R} \underset{\underline{R'}} \boxtimes \underline{A'} \ar[r]& \underline{S}
}\]
gives a lift, unique up to homotopy, of $\underline{R} \underset{\underline{R'}} \boxtimes \underline{A'} \to \underline{A}.$ There is a natural map $\underline{R} \to \underline{R} \underset{\underline{R'}} \boxtimes \underline{A'} \to \underline{A}$, inducing the map $$\Omega_{ \underline{R} \underset{\underline{R'}} \boxtimes \underline{A'}  | \underline{R}} \underset{ \underline{R} \underset{\underline{R'}} \boxtimes \underline{A'} }  \boxtimes \underline{A} \to \Omega_{\underline{A} | \underline{R}}.$$ Applying $(-) \underset{\underline{A}} \boxtimes \underline{S}$ gives a map $$\Omega_{ \underline{R} \underset{\underline{R'}} \boxtimes \underline{A'}  | \underline{R}} \underset{ \underline{R} \underset{\underline{R'}} \boxtimes \underline{A'} }  \boxtimes \underline{S} \to \Omega_{\underline{A} | \underline{R}} \underset{\underline{A}} \boxtimes \underline{S} \cong L_{\underline{S}|\underline{R}}.$$
Now 
\begin{eqnarray*}
L_{\underline{S'}|\underline{R'}} \underset{\underline{R'}} \boxtimes \underline{R} &=& (\Omega_{\underline{A'}|\underline{R'}} \underset{\underline{A'}} \boxtimes \underline{S'}) \underset{\underline{R'}} \boxtimes \underline{R} \\
&\cong&  (\Omega_{\underline{A'}|\underline{R'}} \underset{\underline{A'}} \boxtimes \underline{R}) \underset{\underline{R'} \underset{\underline{R'}} \boxtimes \underline{R}} \boxtimes (\underline{S'} \underset{\underline{R'}} \boxtimes \underline{R}) \\
&\cong&  [\Omega_{\underline{A'}|\underline{R'}} \underset{\underline{A'}} \boxtimes (\underline{A'} \underset{\underline{R'}} \boxtimes \underline{R}) ]\underset{\underline{R'} \underset{\underline{R'}} \boxtimes \underline{R}} \boxtimes (\underline{S'} \underset{\underline{R'}} \boxtimes \underline{R}) \\
&\cong& (\Omega_{\underline{A'} \underset{\underline{R'}} \boxtimes \underline{R}|\underline{R}}) \underset{\underline{R'} \underset{\underline{R'}} \boxtimes \underline{R}} \boxtimes (\underline{S'} \underset{\underline{R'}} \boxtimes \underline{R}). 
\end{eqnarray*}
The second line and third lines are formal. The last line, or that $$  [\Omega_{\underline{A'}|\underline{R'}} \underset{\underline{A'}} \boxtimes (\underline{A'} \underset{\underline{R'}} \boxtimes \underline{R}) ] \cong \Omega_{\underline{A'} \underset{\underline{R'}} \boxtimes \underline{R}|\underline{R}}$$
follows because $\underline{A'}$ is a free $\underline{R'}$-algebra, so derivations are determined by the generators in both cases. Therefore the modules are isomorphic, modulo a change of units from $\underline{A'}$ to $\underline{A'} \underset{\underline{R'}} \boxtimes \underline{R}.$ The canonical map $\underline{S'} \underset{\underline{R'}} \boxtimes \underline{R} \to \underline{S}$ gives us the map $L_{\underline{S'}|\underline{R'}} \underset{\underline{R'}} \boxtimes \underline{R}  \to  L_{\underline{S}|\underline{R}}.$

\end{proof}

A related property is the base change. In the previous proof, if $\underline{S} = \underline{S'} \underset{\underline{R'}} \boxtimes \underline{R}$, then the map $$ L_{\underline{S'}|\underline{R'}} \underset{\underline{R'}} \boxtimes \underline{R} \cong (\Omega_{\underline{A'} \underset{\underline{R'}} \boxtimes \underline{R}|\underline{R}}) \underset{\underline{R'} \underset{\underline{R'}} \boxtimes \underline{R}} \boxtimes (\underline{S'} \underset{\underline{R'}} \boxtimes \underline{R}) \to \Omega_{ \underline{R} \underset{\underline{R'}} \boxtimes \underline{A'}  | \underline{R}} \underset{ \underline{R} \underset{\underline{R'}} \boxtimes \underline{A'} }  \boxtimes \underline{S} $$ is an isomorphism.

Suppose that $\underline{\Tor}_{n}^{\underline{R'}}(\underline{S'},\underline{R}) =0 $ for $n \geq 1$. The calculation for $\Tor$ is the homotopy groups of $\underline{A'} \underset{\underline{R'}} \boxtimes \underline{R}$. If those are all 0 for $n \geq 1$, then the map $\underline{A'} \underset{\underline{R'}} \boxtimes \underline{R} \to c\underline{S}$ is a weak equivalence. So $\underline{A'} \underset{\underline{R'}} \boxtimes \underline{R}$ is a cofibrant replacement for $\underline{R} \to \underline{S}$. By the uniqueness of the cofibrant replacement up to homotopy, $\underline{A'} \underset{\underline{R'}} \boxtimes \underline{R}$ is homotopy equivalent to $\underline{A}$, implying that $L_{\underline{S'}|\underline{R'}} \underset{\underline{R'}} \boxtimes \underline{R}$ is homotopy equivalent to $L_{\underline{S}|\underline{R}}$ via the natural map above.

Next we state a property analogous to the property of (co)homology in topology, that the homology takes finite coproducts to their sums under certain conditions.

\begin{prop}
If $\underline{S'} = \underline{S} \underset{\underline{R}} \boxtimes \underline{R'}$ and $\underline{\Tor}_{n}^{\underline{R'}}(\underline{S'},\underline{R}) =0 $ for $n \geq 1$. Then $$L_{\underline{S'}|\underline{R}} \cong L_{\underline{S}|\underline{R}} \underset{\underline{R}} \boxtimes \underline{R'} \oplus  L_{\underline{R'}|\underline{R}} \underset{\underline{R}} \boxtimes \underline{S}.$$
\end{prop}

\begin{proof}
By the same argument we found in the base change proposition, if $\underline{R} \to \underline{Q} \to \underline{R'}$ and $\underline{R} \to \underline{P} \to \underline{S}$ are cofibrant replacements, then $\underline{R} \to \underline{P} \underset{\underline{R}} \boxtimes \underline{Q} \to \underline{S'}$ is a cofibrant replacement. Therefore,
\begin{eqnarray*}
L_{\underline{S'}|\underline{R}}& \cong& \Omega_{ \underline{P} \underset{\underline{R}} \boxtimes \underline{Q}  | \underline{R}} \underset{ \underline{P} \underset{\underline{R}} \boxtimes \underline{Q} } \boxtimes \underline{S'} \\
&\cong& (\Omega_{\underline{P}|\underline{R}} \underset{\underline{R}} \boxtimes \underline{Q} \oplus \Omega_{\underline{Q}|\underline{R}} \underset{\underline{R}} \boxtimes \underline{P} ) \underset{ \underline{P} \underset{\underline{R}} \boxtimes \underline{Q} } \boxtimes \underline{S'} \\
&\cong& (\Omega_{\underline{P}|\underline{R}} \underset{\underline{P}} \boxtimes \underline{S}) \underset{\underline{R}} \boxtimes \underline{R'} \oplus(\Omega_{\underline{Q}|\underline{R}} \underset{\underline{Q}} \boxtimes \underline{R'}) \underset{\underline{R}} \boxtimes \underline{S}\\
&\cong& L_{\underline{S}|\underline{R}} \underset{\underline{R}} \boxtimes \underline{R'} \oplus  L_{\underline{R'}|\underline{R}} \underset{\underline{R}} \boxtimes \underline{S}.
\end{eqnarray*}
\end{proof}
The property given on (co)-homology is as follows:
\begin{prop}
If $\underline{S'} = \underline{S} \underset{\underline{R}} \boxtimes \underline{R'}$ and $\underline{\Tor}_{n}^{\underline{R'}}(\underline{S'},\underline{R}) =0 $ for $n \geq 1$. Then for an $\underline{S'}$-module $\underline{N}$: $$D^q(\underline{S'}|\underline{S};\underline{N}) \cong D^q(\underline{R}|\underline{R'};\underline{N})$$ $$D^q(\underline{S'}|\underline{R};\underline{N}) \cong D^q(\underline{S}|\underline{R};\underline{N}) \oplus D^q(\underline{R'}|\underline{R};\underline{N})$$ and the identical is true for the homology.
\end{prop}

The next result shows that the homology extends the Jacobi-Zariski sequence. Suppose that we have a short exact $0 \to X' \to X \to X'' \to 0$ of simplicial modules. Formally, this induces a long exact sequence on the homotopy groups of $X',X,X''$, which will be an extension to the left of the Jacobi-Zariski.

\color{red}{Here we offer a correction from the publication of this thesis in 2019. The following proof, modeled after classical proofs of Quillen \cite{quillen1970co} and Andr\'e \cite{andre1974homologie}, contains an inaccuracy. Namely, the property of pushouts of weak equivalences along cofibrations being weak equivalent requires the category to be left proper, which the category of rings is, but the category of $\mathcal{O}$-Tambara functors is not. See the thesis of Stahlhauer \cite{stahlhauer} which points out this issue.

In order to resolve this, we give the limitation that the map $\underline{R} \to \underline{S}$ be a cofibration, which allows us to make the argument about weak equivalence. The proof can be made much simpler but we leave it as unchanged as possible from the 2019 publication.

As a side note through private communication with Mike Hill: this argument would hold for $E^\infty$ algebras rather than strict commutative monoids as that category is left proper.
}
\color{black}

\begin{prop}
Let $\underline{R} \overset{u} \to \underline{S} \overset{v} \to \underline{T}$ be a map of $\mathcal{O}$-Tambara functors \color{red} where  $\underline{R} \to \underline{S}$ is a cofibration.\color{black} Then $$ L_{\underline{S}|\underline{R}} \boxtimes \underline{T} \to L_{\underline{T}|\underline{R}} \to L_{\underline{T}|\underline{S}} \to $$ is an exact triangle in the derived category of $\underline{T}$-modules. Tensoring with an $\underline{T}$-module $\underline{M}$ and taking the homotopy, we get a long exact sequence 
 \[\xymatrixrowsep{.5in}
\xymatrixcolsep{.12in}\xymatrix{
\cdots \ar[r] & D_1(\underline{T}|\underline{R};\underline{M}) \ar[r] & D_1(\underline{T}|\underline{S};\underline{M}) \ar[r]   & D_0(\underline{S}|\underline{R};\underline{M}) \ar[r] & D_0(\underline{T}|\underline{R};\underline{M}) \ar[r] & D_0(\underline{T}|\underline{S};\underline{M}) \ar[r] & 0
}\]
or
 \[\xymatrixrowsep{.5in}
\xymatrixcolsep{.12in}\xymatrix{
\cdots \ar[r] & D_1(\underline{T}|\underline{R};\underline{M}) \ar[r] & D_1(\underline{T}|\underline{S};\underline{M}) \ar[r]   & \underline{M} \underset{\underline{S}} \boxtimes \Omega_{\underline{S}|\underline{R}} \ar[r] & \underline{M} \underset{\underline{T}} \boxtimes \Omega_{\underline{T}|\underline{R}}  \ar[r] & \underline{M} \underset{\underline{T}} \boxtimes \Omega_{\underline{T}|\underline{S}} \ar[r] & 0.
}\]
while taking $\underline{\Hom}$ to a module and taking the homotopy gives the following 
 \[\xymatrixrowsep{.5in}
\xymatrixcolsep{.15in}\xymatrix{
0 \ar[r] & D^0(\underline{T}|\underline{S};\underline{M}) \ar[r] & D^0(\underline{T}|\underline{R};\underline{M}) \ar[r] & D^0(\underline{S}|\underline{R};\underline{M}) \ar[r]   & D^1(\underline{T}|\underline{S};\underline{M}) \ar[r] & D^1(\underline{T}|\underline{R};\underline{M}) \ar[r] & \cdots
}\]
or
 \[\xymatrixrowsep{.5in}
\xymatrixcolsep{.15in}\xymatrix{
0 \ar[r] & \underline{\Der}(\underline{T}|\underline{S};\underline{M}) \ar[r] & \underline{\Der}(\underline{T}|\underline{R};\underline{M}) \ar[r] & \underline{\Der}(\underline{S}|\underline{R};\underline{M}) \ar[r]   & D^1(\underline{T}|\underline{S};\underline{M}) \ar[r] & D^1(\underline{T}|\underline{R};\underline{M}) \ar[r] &  \cdots
}\]
\end{prop}
A clean resulting statement is that $L_{\underline{S}|\underline{R}} \boxtimes \underline{T} \to L_{\underline{T}|\underline{R}} \to L_{\underline{T}|\underline{S}}  \to \Sigma(L_{\underline{S}|\underline{R}} \boxtimes \underline{T} )$ is a cofiber in the homotopy category of $\underline{T}$-modules.

\begin{proof}
Let $\underline{X}$ be the cofibrant replacement of $\underline{R} \to \underline{S}$, and let $\underline{Y}$ be the cofibrant replacement of $\underline{X} \to \underline{T}$, and take $\underline{X}$ and $\underline{Y}$ to be free.
So we have a map of rings $\underline{R} \to \underline{X} \to \underline{Y}$, giving a map (the Jacobi-Zariski map)

 \[\xymatrixrowsep{.5in}
\xymatrixcolsep{.25in}\xymatrix{
\underline{Y} \underset{\underline{X}} \boxtimes \Omega_{\underline{X}|\underline{R}} \ar[r] & \Omega_{\underline{Y}|\underline{R}} \ar[r] &  \Omega_{\underline{Y}|\underline{X}} \ar[r] & 0
}\]
which, because $\underline{T} \underset{\underline{Y}} \boxtimes (-)$ is a right exact functor, gives the exact sequence
 \[\xymatrixrowsep{.5in}
\xymatrixcolsep{.25in}\xymatrix{
\underline{T} \underset{\underline{X}} \boxtimes \Omega_{\underline{X}|\underline{R}} \ar[r] & \underline{T} \underset{\underline{Y}} \boxtimes \Omega_{\underline{Y}|\underline{R}} \ar[r] &  \underline{T} \underset{\underline{Y}} \boxtimes  \Omega_{\underline{Y}|\underline{X}} \ar[r] & 0
}\] which is in fact a short exact sequence
 \[\xymatrixrowsep{.5in}
\xymatrixcolsep{.25in}\xymatrix{
0 \ar[r] & \underline{T} \underset{\underline{X}} \boxtimes \Omega_{\underline{X}|\underline{R}} \ar[r] & \underline{T} \underset{\underline{Y}} \boxtimes \Omega_{\underline{Y}|\underline{R}} \ar[r] &  \underline{T} \underset{\underline{Y}} \boxtimes  \Omega_{\underline{Y}|\underline{X}} \ar[r] & 0.
}\] 
This is because the induced map $\Der(\underline{Y}|\underline{R}, \underline{M}) \to \Der(\underline{X}|\underline{R},\underline{M})$ is surjective, because $\underline{Y}$ is free over $\underline{X}$. Now we need to relate the last term to $L_{\underline{T}|\underline{S}}$.

A pushout of a cofibration is a cofibration, so $\underline{S} \to \underline{S} \underset{\underline{X}} \boxtimes \underline{Y}$ is a cofibration. Additionally, pushout of a weak equivalence between cofibrant objects (which $\underline{X}$ and \color{red}$\underline{S}$ \color{black} are) along a cofibration is a weak equivalence, see \cite{hirschhorn2009model} chapter 13. So $\underline{Y} \to \underline{S} \underset{\underline{X}} \boxtimes \underline{Y}$ is a weak equivalence, so by the 2-out-of-3 axiom $\underline{S} \underset{\underline{X}} \boxtimes \underline{Y} \to \underline{T}$ is a weak equivalence. It is surjective, because $\underline{Y} \to \underline{T}$ is, so it is a fibration. Taken together, we have that $\underline{S} \to \underline{S} \underset{\underline{X}} \boxtimes \underline{Y} \to \underline{T}$ is a cofibrant replacement, an $L_{\underline{T}|\underline{S}} = \underline{T} \underset{ \underline{S} \underset{\underline{X}} \boxtimes \underline{Y}} \boxtimes \Omega_{ \underline{S} \underset{\underline{X}} \boxtimes \underline{Y} | \underline{S}}.$ Finally, $\underline{T} \underset{ \underline{S} \underset{\underline{X}} \boxtimes \underline{Y}} \boxtimes \Omega_{ \underline{S} \underset{\underline{X}} \boxtimes \underline{Y} | \underline{S}} \cong \underline{T} \underset{\underline{Y}} \boxtimes  \Omega_{\underline{Y}|\underline{X}}  $ because $\underline{Y}$ is a free $\underline{X}$-algebra, so derivations are given by where the generators map, giving the isomorphism modulo change of units. Therefore we have the exact sequence
 \[\xymatrixrowsep{.5in}
\xymatrixcolsep{.25in}\xymatrix{
0 \ar[r] & \underline{T} \underset{\underline{X}} \boxtimes \Omega_{\underline{X}|\underline{R}} \ar[r] & \underline{T} \underset{\underline{Y}} \boxtimes \Omega_{\underline{Y}|\underline{R}} \ar[r] & \underline{T} \underset{ \underline{S} \underset{\underline{X}} \boxtimes \underline{Y}} \boxtimes \Omega_{ \underline{S} \underset{\underline{X}} \boxtimes \underline{Y} | \underline{S}} \ar[r] & 0
}\] 
which gives the desired cofiber sequence.

\end{proof}

%% file: SpectralSeq.tex
\chapter{Fundamental Spectral Sequence}
\section{Creating the Spectral Sequence}

\color{red} Throughout this chapter, we assume that the map $\underline{A} \to \underline{P}$ is a cofibration in order to use the transitivity long exact sequence. See the correction at the end of the last chapter. \color{black}

Quillen \cite{QuillenMimeo} created a spectral sequence for computing the Andr\'e-Quillen cohomology groups in the non-equivariant case. The $E^2$ page involves homology of algebras generated by the cotangent complex, and it converges something more easily computable, Tor which does not need free resolutions of algebras but rather just resolutions of modules. The spectral sequence is therefore useful for trying to compute the cohomology groups by going backwards. First we make a reduction to the case where $\underline{B}$ is a quotient of $\underline{A}$, just as Quillen did in the non-equivariant setting.

Note that for any $\underline{A} \to \underline{B}$, we can express $\underline{B}$ as a quotient of a free algebra $\underline{P}$ over $\underline{A}$, so we have maps $\underline{A} \to \underline{P} \to \underline{B}$ and $\underline{B} \cong \underline{P}/\underline{I}$, for some ideal $\underline{I}$. Since $\underline{P}$ is free over $\underline{A}, D^q(\underline{P}|\underline{A}; \underline{M}) = 0$, $D_q(\underline{P}|\underline{A}; \underline{M}) = 0$ for $q > 0.$ Therefore, given the transitivity exact sequence applied to $\underline{A} \to \underline{P} \to \underline{B}$, we get $D^q(\underline{B}|\underline{P};\underline{M}) \cong D^q(\underline{B}|\underline{A};\underline{M})$ when $q>1$ (and identically for cohomology). In terms of calculating the higher cohomology groups, we can consider the case where $\underline{B}$ is a quotient, or $\underline{B} \cong \underline{A}/\underline{I}.$ In addition to the above isomorphisms, in the low degree case we have the following sections of a long exact sequence
 \[\xymatrixrowsep{.5in}
\xymatrixcolsep{.15in}\xymatrix{
 0 \ar[r] & D_1(\underline{B}|\underline{A};\underline{M}) \ar[r] & D_1(\underline{B}|\underline{P};\underline{M})  \ar[r]   & \underline{M} \underset{\underline{P}} \boxtimes \Omega_{\underline{P}|\underline{A}} \ar[r] & \underline{M} \underset{\underline{B}} \boxtimes \Omega_{\underline{B}|\underline{A}}  \ar[r] & \underline{M} \underset{\underline{B}} \boxtimes \Omega_{\underline{B}|\underline{P}} = 0 .
}\]

 \[\xymatrixrowsep{.5in}
\xymatrixcolsep{.15in}\xymatrix{
\underline{\Der}(\underline{B}|\underline{P};\underline{M}) = 0 \ar[r] & \underline{\Der}(\underline{B}|\underline{A};\underline{M}) \ar[r] & \underline{\Der}(\underline{P}|\underline{A};\underline{M}) \ar[lld]   & \hspace{.05in} \\  D^1(\underline{B}|\underline{P};\underline{M}) \ar[r] & D^1(\underline{B}|\underline{A};\underline{M}) \ar[r] & D^1(\underline{P}|\underline{A};\underline{M})=0.
}\]
Also, because $D_1(\underline{B}|\underline{P}) = \underline{I}/\underline{I}^{>1}$, the K\"unneth spectral sequence gives that $D_1(\underline{B}|\underline{P};\underline{M}) = \underline{I}/\underline{I}^{>1} \underset{\underline{B}} \boxtimes \underline{M}.$

If $\underline{P}$ is a free simplicial $\underline{A}$-algebra resolution of $\underline{B}$, then applying $(-) \underset{\underline{A}} \boxtimes \underline{B}$ to $\underline{A} \to \underline{P} \to \underline{B}$ yields $$\underline{B} \to \underline{P} \underset{\underline{A}} \boxtimes \underline{B} \to \underline{B} \underset{\underline{A}} \boxtimes \underline{B} \to \underline{B}$$ with the last arrow being the $\underline{A}$-algebra structure of $\underline{B}.$ If $\underline{P}$ is free over $\underline{A}$, then $ \underline{P} \underset{\underline{A}} \boxtimes \underline{B}$ is free over $\underline{B}$. Considering the sequence $\underline{B} \to \underline{P} \underset{\underline{A}} \boxtimes \underline{B} \to \underline{B}$ and applying the above argument dimension-wise, realizing that $D_i(\underline{B}|\underline{B}) = 0$ and similarly for the cohomology, we get the isomorphisms 
$$D_1(\underline{B}|\underline{P}_* \underset{\underline{A}} \boxtimes \underline{B} ) \cong \underline{B} \underset{\underline{P}_* \underset{\underline{A}} \boxtimes \underline{B} } \boxtimes \Omega_{\underline{P}_* \underset{\underline{A}} \boxtimes \underline{B}|\underline{B}}. $$ Which, again by the above, if $\underline{J}$ is the kernel of $\underline{P} \underset{\underline{A}} \boxtimes \underline{B} \to \underline{B}$ then $$(\underline{J}/\underline{J}^{>1})_* \cong \underline{B} \underset{\underline{P}_* \underset{\underline{A}} \boxtimes \underline{B} } \boxtimes \Omega_{\underline{P}_* \underset{\underline{A}} \boxtimes \underline{B}|\underline{B}}$$ in each dimension, so $$\underline{J}/\underline{J}^{>1} \cong \underline{B} \underset{\underline{P} \underset{\underline{A}} \boxtimes \underline{B} } \boxtimes \Omega_{\underline{P} \underset{\underline{A}} \boxtimes \underline{B}|\underline{B}} \cong \underline{B} \underset{\underline{P}} \boxtimes \Omega_{\underline{P} |\underline{A}} \cong L_{\underline{B}|\underline{A}}$$ which is another way to calculate the cotangent complex.

Now, we have the filtration $$\underline{Q} =  \underline{P} \underset{\underline{A}} \boxtimes \underline{B} \supset \underline{J} \supset \underline{J}^{>1} \supset \underline{J}^{>2} \supset \cdots$$
of simplicial ideals of $\underline{P} \underset{\underline{A}} \boxtimes \underline{B}$. 
The filtration, after taking the normalization, naturally gives a spectral sequence with $$E_{pq}^0 = (N\underline{J}^{>-p-1})_{p+q} / (N\underline{J}^{>-p})_{p+q}, \quad E_{pq}^1 = H_{p+q}(\underline{J}^{>-p-1}/\underline{J}^{>-p})$$ where $\underline{J}^{>0} = \underline{J}$ and $\underline{J}^{>-1} = \underline{Q}$ and is 0 otherwise. Note that the filtration is clearly exhaustive and bounded above.

Now we must relate the other terms to the cotangent complex. $\underline{J}$ is the ideal generated by the indeterminants in $\underline{Q}.$ Another way of constructing this would be to take the Mackey functor freely generated over $\underline{B}$ by the same elements, and then consider taking the symmetric algebra functor over $\underline{B}$. This functor assigns norms to all elements and agrees with the norm maps on $\underline{B}$. The Mackey functor generated by the same elements is $\underline{J}/\underline{J}^{>1}$. So we have an identification $\underline{J} \cong \Sigma^{\underline{B}} (\underline{J}/\underline{J}^{>1})$, and we see that $$(\Sigma^{\underline{B}}_{q} (\underline{J}/\underline{J}^{>1})) := (\Sigma^{\underline{B}}_{>q} (\underline{J}/\underline{J}^{>1})) / (\Sigma^{\underline{B}}_{>q+1} (\underline{J}/\underline{J}^{>1})) \cong \underline{J}^{>q} /\underline{J}^{>q+1}$$ where $\Sigma^{\underline{B}}_{>q} \subset \Sigma^{\underline{B}}$ is the sub-algebra containing all the elements coming from a norm $>q.$

Combining all of the above, we have the following theorem:
\begin{prop}
If $\underline{B} \underset{\underline{A}} \boxtimes \underline{B} \cong \underline{B}$ (for example $\underline{A}$ surjects onto $\underline{B}$), then there is a  spectral sequence $$E_{pq}^1 = H_{p+q}(\Sigma_{-p-1}^{\underline{B}}(L_{\underline{B}|\underline{A}})) \Longrightarrow \underline{\Tor}_{p+q}^{\underline{A}}(\underline{B},\underline{B}).$$
\end{prop}

Note that we can call the $E_1$ page the $E_2$ page, and reindex everything to make the differentials correct, so that the terms that 0 given by the lemma are in the second quadrant:
\begin{prop}
If $\underline{B} \underset{\underline{A}} \boxtimes \underline{B} \cong \underline{B}$, then there is a spectral sequence $$E_{pq}^2 = H_{p+q}(\Sigma_q^{\underline{B}}(L_{\underline{B}|\underline{A}})) \Longrightarrow \underline{\Tor}_{p+q}^{\underline{A}}(\underline{B},\underline{B}). $$
\end{prop}
In the classical case, the spectral sequence above was a first quadrant spectral sequence, proving the convergence. This is no longer necessarily true, as we will see non-trivial terms in the second quadrant.

%% file: Convergence.tex
\section{Convergence}

Much of the argument here about convergence is similar to that given by Andr\'e in \cite{andre1974homologie} for the classical case. The overall structure is the same, but the presence of non-trivial norms provides non-trivial hurdles. In particular, for the classical case $I \otimes I \to I^2$ is a surjection, whereas in the equivariant case the image of $\underline{I} \boxtimes \underline{I} \to \underline{I}^{>1}$ does not necessarily contain monomials with just one nontrivial norm of order 2. To fix this issue, we must consider the numbers that exhibit sharing for $|G|$.

A much larger issue is that free $\mathcal{O}$-Tambara functors are not necessarily free, projective, or even flat as Mackey functors. In order to deal with flatness issues, we consider only finite groups $G$ and indexing system for which the free Tambara functor $\underline{A}^{\mathcal{O}}[x_H]$ are flat as Mackey functors. We call this condition on $\mathcal{O}$ being \term{free-flat}. The condition of being free-flat is very restrictive, and we do not have an example outside of $G$ being the trivial group for which this holds. Even for $C_p$, the computations of \cite{blumberg2019right} section 3 show that for both the trivial indexing system or complete indexing system, the free-flat condition does not hold by inspection.

In personal communication with Mike Hill, it was relayed to the author that the free-flat condition could be sidestepped by considering operadic algebras in simplicial Mackey functors. The author hopes in the future to resolve this issue and show convergence of these spectral sequences completely in a later joint paper:

\begin{conj}
For any finite group $G$ and indexing system $\mathcal{O}$, the fundamental spectral sequence of the pervious section converges. In particular, Corollary 4.2.10 holds without the flat-free condition on $\mathcal{O}$.
\end{conj}

Here, we include this proof with this very restrictive hypothesis to show the difficulty in showing convergence and as a basis for future work.

\begin{lemma}
Suppose $\mathcal{O}$ is free-flat. Then for any $\mathcal{O}$-Tambara functor $\underline{B}$, any free $\underline{B}$-algebra is flat as a $\underline{B}$-module.
\end{lemma}

\begin{proof}
Box products of flat-modules are free, so any free  $\underline{A}^{\mathcal{O}}$-algebra on finitely many generators is free. For infinitely many generators, the direct limit of flat modules are flat, so any free $\underline{A}^{\mathcal{O}}$-algebra is flat as an $\underline{A}^{\mathcal{O}}$-module. Lastly, $(\underline{B} \boxtimes \underline{A}^{\mathcal{O}}[x_{G/Hi}]) \underset{\underline{B}} \boxtimes (-) \cong \underline{A}^{\mathcal{O}}[x_{G/Hi}] \boxtimes (-)$ and the right side is exact, so we have the lemma.
\end{proof}

Next we have to resolve the $I \otimes I \to I^2$ discrepancy, which we do with the following notion.
\begin{definition}
Let $n$ be a positive integer. Suppose there is a positive integer $m$ such that for every positive integer $k$, every partition of $k \cdot m$ into divisors of $n$ is the union of $k$ many partitions of $m$. Then we say that \textbf{$m$ exhibits sharing for $n$}. If any such $m$ exists for $n$ then we say that $n$ is a \textbf{sharing number}.
\end{definition}

To get a handle on the definition, the following are true, most of them trivially:
\begin{itemize}
\item 1 is a sharing number and any $m$ exhibits sharing for 1.
\item 2 is a sharing number and any even $m$ exhibits sharing for 2.
\item Moreover, $p^n$ is a sharing number for any prime $p$, with $p^n$ exhibiting sharing for $p^n$.
\item $p^n \cdot q$ is a sharing number for any primes $p$ and $q$, with $p^n \cdot q$ exhibiting sharing for $p^n \cdot q.$
\item 36 exhibits sharing for 36.
\item 30 does not exhibit sharing for 30. The partition $(15,10,10,6,6,6,6,1)$ is a partition of 60 that is not the union of two partitions of 30. It is the first number which does not exhibit sharing for itself. 60 exhibits sharing for 30.
\end{itemize}
 
\begin{lemma}
Let $n$ be a positive integer. If $d$ is a proper divisor of $n$, let $f(d)$ be the smallest multiple of $d$ that is also a divisor of $n$. Let $m$ be the first multiple of $n$ that is bigger than $$\sum_{d | n, d \neq n} f(d) - d.$$ Then $m$ exhibits sharing for $n$. In particular, all positive integers are sharing numbers.
\end{lemma}

\begin{proof}
Suppose we are given a partition $p$ of $k \cdot m$ for the above $m$ into divisors of $n$. Suppose that for the partition $p$ of $k \cdot m$, that a particular divisor $d$ appears $f(d)/d$ times. Let $p'$ be the partition which is the same as $p$, except $f(d)/d$ many instances of $d$ are replaced by a single instance of $f(d)$. If $p'$ can be split into $k$ equally valued partitions, then trivially so can $p$. Doing this replacement arbitrarily many times, we can assume our partition has, for a proper divisor $d$ of $n$, at most $f(d)/d - 1$ instances of $d$. Because of the inequality $$m > \sum_{d | n, d \neq n} f(d) - d$$ then the number of instances of $n$ must sum up to more than $(k-1)\cdot m$, as the sum $(\sum_{d | n, d \neq n} f(d) - d)$ with the instance of $n$ must equal $k \cdot m$. But if in our partition, the instances of $n$ add up to more than $(k-1)\cdot m$, then the partition is trivially a union of $k$ many such partitions of $m$.
\end{proof}

For all of the following we will assume that $m$ is an integer that exhibits sharing for $|G|.$ Before getting to the convergence, we must show various lemmas. The convergence proof involves computing the homology of various box products that relate to the homology of interest by various exact sequences, so the following lemma will be critical.

\begin{lemma} \label{lem8.18}
Suppose we have 2 simplicial $\underline{A}$-modules $\underline{E}_*'$ and $\underline{E}_*''$ such that homology for $\underline{E}_*'$ vanishes up to and including degree $m$ and likewise for $\underline{E}_*''$ up to and including degree $n$. Then the homology of $\underline{E}_*' \underset{\underline{A}} \boxtimes \underline{E}_*''$ vanishes up to degree $m+n+1$ if one of the modules $\underline{E}_*'$ or $\underline{E}_*''$ is projective.
\end{lemma}

\begin{proof}
First we must construct a few additionally simplicial modules: consider first
$$0 \to \underline{F}'_* \to  \underline{\Gamma}_*' \to  \underline{E}_*' \to 0$$
where $ \underline{\Gamma}_n' =  \underline{E}_{n+1}'$ and for an element $x \in  \underline{\Gamma}_n'(X)$, we let $\epsilon_n^i(x) = \epsilon_{n+1}^{i+1}(x)$ for $x \in  \underline{E}_{n+1}'(X)$ and $\sigma_{n}^i(x) = \sigma_{n+1}^{i+1}(x)$, where $\epsilon$ and $\sigma$ represent face and degeneracy maps respectively. The map $\gamma':  \underline{\Gamma}_*' \to  \underline{E}_*'$ sends $x$ to $\epsilon^0_{n+1}(x).$ It is an easy check that this map $\gamma$ is a map of simplicial modules and that $\gamma$ is surjective as $\sigma_{n}^0(x)$ provides a pre-image. We now show various properties of this sequence.

Firstly, the non-zero homotopy Mackey functors of $\underline{\Gamma_*}$ are all zero: consider $x \in \underline{E}_{n+1}(X)$ with $\epsilon_{n+1}^i(x) = 0$ and $1 \leq i \leq n+1.$ So one would need to find an element $y \in \underline{E}_{n+1}(X)$ such that $\epsilon_{n+2}^i(y) = 0$ for $1 \leq i \leq n+1$ and $\epsilon_{n+2}^{n+2}(y) = x.$ This necessarily exists because of the Kan condition on $\underline{E}_*$ because $\underline{E}_*$ is an abelian group object.

Secondly, because $\underline{\Gamma}_n = \underline{E}_n$, if $\underline{E}_n$ is projective then so is $\underline{\Gamma}_n.$ $\underline{F}$, being the quotient of a map between projective objects is also projective.

Next, when $H_0[\underline{E}_*'] \cong 0$, we consider $$0 \to  \underline{G}_*' \to  \underline{\Delta}_*' \to  \underline{E}_*' \to 0.$$
$\underline{\Delta}_*$ is a submodule of $\underline{\Gamma}_*$. The modules $\underline{\Delta}_n$ are formed by elements of $\underline{E}_{n+1}$ such that $$(\epsilon_1^1 \circ \epsilon_2^2 \circ \cdots \circ \epsilon_{n+1}^{n+1})(x) = 0.$$ One can check, using the simplicial identities, that $\underline{\Delta}_*$ is a simplicial module. Additionally, because $H_0[\underline{E}_*']$ is 0, for every element $x \in \underline{E}_n'(X)$, there is an element $\tau \in \underline{E}_1'(X)$ such that $\epsilon_1^0(\tau) =0$ and $\epsilon_1^1(\tau) = (\epsilon_1^1 \circ \cdots \circ \epsilon_{n+1}^{n+1} \circ \sigma_n^0)(x)$. Consider $y \in \underline{E}_{n+1}'$ defined by $y = \sigma_n^0(x)-(\sigma_n^n \circ \cdots \circ_1^1)(\tau)$. Using the simplicial identities, we see that $\gamma(y)=x$ and that $y$ is an element of $\underline{\Delta}_n(X),$ so $\delta_n:\underline{\Delta}_n \to \underline{E}_n'$ is a surjection.

The positive homotopy Mackey functors of $\underline{\Delta}_n'$ also vanish by the same reason as $\underline{\Gamma}_*'$. The zeroth homotopy Mackey functor is also 0: consider $x \in \underline{E}_1'$ with $\epsilon_1^1(x) = 0$. Then again by the Kan condition there exists $y \in \underline{E}_2'$ such that $\epsilon_2^1(y) = 0$ and $\epsilon_2^2(y) = x$, so in particular $(\epsilon_1^1 \circ \epsilon_2^2)(y) = 0$. For every element $x \in \underline{\Delta}_0$, there is an element $y \in \underline{\Delta}_1$ such that $\epsilon_1^0(y) = 0$ and $\epsilon_1^1(y) = x$. 

Lastly, because of the exact sequence $$0 \to \underline{\Delta}_n' \to \underline{\Gamma}_n' \to \underline{E}_0' \to 0,$$ if $\underline{E}_*'$ is projective, then so are $\underline{\Delta}_n'$ and $\underline{G}_*'$.

On to the proof of the lemma. We go by induction on $s = m+n.$ We can assume that one of the entries, say $m \geq 0.$ For every element $x'$ of $\underline{E}_0'(X')$, there exists an element $y'$ in $\underline{E}_1'(X)$ with $\epsilon_1^0(y') = 0$ and $\epsilon_1^1(y') = x'$ and for every element $x''$ of $\underline{E}_0''(X)$, there exists an element $Y''$ of $\underline{E}_1''(X)$ such that $\epsilon_1^0(y'') = x''$ and $\epsilon_1^1(y'') = x''.$ The equations $$\epsilon_1^0(y' \otimes y'') = 0$$ and $$\epsilon_1^1(y' \otimes y'') = x' \otimes x''$$ shows that the $H_0[\underline{E}_*' \underset{\underline{A}} \boxtimes \underline{E}_*''] \cong 0$, so the lemma is shown for $s= -1.$ It now remains to show that given $H_{m+n+1}[\underline{E}_*' \underset{\underline{A}} \boxtimes \underline{E}_*'']$ using the induction hypothesis.

From the inductive hypothesis on $H_p[\underline{\Delta}_*'] \cong 0$ for $p \leq s$ and $H_q[\underline{F}_*''] \cong 0$ for $q \leq -1$, we have that $H_s[\Delta_*' \underset{\underline{A}} \boxtimes \underline{F}_*''] \cong 0.$ Note that this is true because at least one of them is flat.
Similarly, since $H_p[\underline{G}_*'] \cong 0$ for $p \leq m-1$ and $H_q[\underline{E}''_*] \cong 0$ for $q \leq n$ gives $H_s[\underline{G}_*' \underset{\underline{A}} \boxtimes \underline{E}''_*] \cong 0.$ 

We have the following short exact sequence of simplicial $\underline{A}$-modules
$$0 \to \underline{\Delta}_*' \underset{\underline{A}} \boxtimes \underline{F}''_* \to  \underline{\Delta}_*' \underset{\underline{A}} \boxtimes \underline{\Gamma}''_* \to  \underline{\Delta}_*' \underset{\underline{A}} \boxtimes \underline{E}''_*\to 0 $$
which gives the following part of the long exact sequence:
$$ H_{s+1}[\underline{\Delta}_*' \underset{\underline{A}} \boxtimes \underline{\Gamma}''_*] \to  H_{s+1}[\underline{\Delta}_*' \underset{\underline{A}} \boxtimes \underline{E}''_*] \to H_s[\underline{\Delta}_*' \underset{\underline{A}} \boxtimes \underline{F}''_* ].$$ The projectivity of either $\underline{\Delta}_*'$ or $\underline{E}_*''$ insures the exactness. Because the third term is 0, the first homomorphism is surjective. In the same vein, we have 
$$0 \to \underline{G}_*' \underset{\underline{A}} \boxtimes \underline{E}''_* \to  \underline{\Delta}_*' \underset{\underline{A}} \boxtimes \underline{E}''_* \to  \underline{E}_*' \underset{\underline{A}} \boxtimes \underline{E}''_*\to 0 $$
which gives the following part of the long exact sequence:
$$ H_{s+1}[\underline{\Delta}_*' \underset{\underline{A}} \boxtimes \underline{E}''_*] \to  H_{s+1}[\underline{E}_*' \underset{\underline{A}} \boxtimes \underline{E}''_*] \to H_s[\underline{G}_*' \underset{\underline{A}} \boxtimes \underline{F}''_* ]$$ with the first map being surjective. Therefore, the composites $$H_{s+1}[\underline{\Delta}_*' \underset{\underline{A}} \boxtimes \underline{\Gamma}_*''] \to  H_{s+1}[\underline{E}_*' \underset{\underline{A}} \boxtimes \underline{E}_*'']$$ and $$H_{s+1}[\underline{\Gamma}_*' \underset{\underline{A}} \boxtimes \underline{\Gamma}_*''] \to  H_{s+1}[\underline{E}_*' \underset{\underline{A}} \boxtimes \underline{E}_*'']$$ are surjective. But we can consider the exact sequence
 $$0 \to \underline{F}_* \to \underline{\Gamma}_* = \underline{\Gamma}_*' \underset{\underline{A}} \boxtimes \underline{\Gamma}_*'' \to \underline{E}_* = \underline{E}_*' \underset{\underline{A}} \boxtimes \underline{E}''_* \to 0.  $$ The surjective homomorphism $$H_{s+1}[\underline{\Gamma}_*' \underset{\underline{A}} \boxtimes \underline{\Gamma}_*''] \to H_{s+1}[\underline{E}_*' \underset{\underline{A}} \boxtimes \underline{E}_*'']  $$ has domain zero since $s+1$ is non-zero, which shows that the homology is 0 as desired.

\end{proof}

\begin{lemma}\label{lem12.15}
Let $\underline{A}$ be a free $\underline{B}$-algebra, with the projection map $\underline{A} \to \underline{B}$ having kernel $\underline{K}$. Then for $p \neq 0$, the map $$\underline{\Tor}_{p}^{\underline{A}}(\underline{A}/ \underline{K}^{\geq m}, \underline{A}/ (\underline{K}^{\geq m})^n) \to \underline{\Tor}_{p}^{\underline{A}}(\underline{A}/ \underline{K}^{\geq m}, \underline{A}/ (\underline{K}^{\geq m})^{n-1})$$ induced from the short exact sequence $$0 \to (\underline{K}^{\geq m})^{n-1} / (\underline{K}^{\geq m})^n \to  \underline{A}/ (\underline{K}^{\geq m})^n \to \underline{A}/ (\underline{K}^{\geq m})^{n-1} \to 0 $$ is the zero map.

\end{lemma}

\begin{proof}
Let the projective resolution $$\cdots \to \underline{P}_2 \to \underline{P}_1 \to \underline{P}_0 \to \underline{B}$$ be the resolution given by the free-forgetful adjunction where $\underline{P}_0 = \underline{A}$. That is $\underline{P}_1 = \underline{A}[K_0]$ the free $\underline{A}$-module generated by $K$, $\underline{K}$ treated as a $G$-indexed set, and $\underline{P}_n = \underline{A}[K_{n-1}]$ where $K_{n-1}$ is the kernel of $\underline{P}_{n-1} \to \underline{P}_{n-2}$ again treated as a $G$ indexed set. 
An element of $\underline{P}_1(X)$ is of the form $T_{f_1}(a_1 \cdot R_{g_1}[k])$ where $f_1: Y \to X, g_1: Y \to Z, a_1 \in \underline{A}(Y)$ and $k$ is a formal element of the $G$-indexed set $\underline{K}$ (if $X,Y,Z$ are not orbits, then the formal elements are always tuples of formal elements from orbits). 
An element of $\underline{P}_n$ is of the form $$T_{f_n}(a_n \cdot R_{g_n}[T_{f_{n-1}}(a_{n-1} \cdot R_{g_{n-1}}[\cdots T_{f_1}(a_1 \cdot R_{g_1}[k]) \cdots ])])$$ with the condition that if we remove the first set of formal brackets in $$T_{f_k}(a_k \cdot R_{g_k}[T_{f_{k-1}}(a_{k-1} \cdot R_{g_{k-1}}[\cdots T_{f_1}(a_1 \cdot R_{g_1}[k]) \cdots ])])$$ it yields  $$T_{f_k}(a_k \cdot R_{g_k}T_{f_{k-1}}(a_{k-1} \cdot R_{g_{k-1}}[\cdots T_{f_1}(a_1 \cdot R_{g_1}[k]) \cdots ])) = 0$$ the zero element of $\underline{P}_{k-1}$ for all $k \neq n.$

Suppose we have an element of $\underline{\Tor}_{p}^{\underline{A}}(\underline{A}/ \underline{K}^{\geq m}, \underline{A}/ (\underline{K}^{\geq m})^n).$ This can be represented by an element of $\underline{P}_n \underset{A} \boxtimes (\underline{A}/ \underline{K}^{\geq m})^n$ which appears of the form $$T_{f_n}(a_n \cdot R_{g_n}[T_{f_{n-1}}(a_{n-1} \cdot R_{g_{n-1}}[\cdots T_{f_1}(a_1 \cdot R_{g_1}[k]) \cdots ])])$$ such that $a_n$ (and $a_n$ alone, not all $a_i$) is an element of $ \underline{A}/ (\underline{K}^{\geq m})^n$ with the property that $$T_{f_n}(a_n \cdot R_{g_n}T_{f_{n-1}}(a_{n-1} \cdot R_{g_{n-1}}[\cdots T_{f_1}(a_1 \cdot R_{g_1}[k]) \cdots ])) =0 $$ is 0 in $ \underline{P}_{n-1} \underset{A} \boxtimes \underline{A}/ (\underline{K}^{\geq m})^n.$ Therefore, $$T_{f_n}(a_n \cdot R_{g_n}T_{f_{n-1}}(a_{n-1} \cdot R_{g_{n-1}}[\cdots T_{f_1}(a_1 \cdot R_{g_1}[k]) \cdots ])) = T_f( c \cdot R_g ([k_{n-1}]))$$ as elements in $\underline{P}_{n-1} $ for $p_{n-1} \in \underline{K}_{n-1}$ and $c \in (\underline{K}^{\geq m})^n.$ The element $c$ can be written as $T_h(c_1 c_{n-1})$ by the sharing property, where $c_1 \in (\underline{K}^{\geq m})$ and $c_{n-1} \in (\underline{K}^{\geq m})^n.$ So the elements can be rearranged into TNR form for $$T_f( c \cdot R_g ([k_{n-1}])) = T_fT_h( (c_1 c_{n-1}) \cdot R_h R_g ([k_{n-1}])) $$Therefore, the element $$[T_{f_n}(a_n \cdot R_{g_n}[T_{f_{n-1}}(a_{n-1} \cdot R_{g_{n-1}}[\cdots T_{f_1}(a_1 \cdot R_{g_1}[k]) \cdots ])]) - T_fT_h( c_{n-1} \cdot [c_1\cdot R_h R_g ([k_{n-1}]))]  ]$$ is in $\underline{P}_{n+1} \underset{A} \boxtimes (\underline{A}/ \underline{K}^{\geq m})^n$ as $$T_{f_n}(a_n \cdot R_{g_n}[T_{f_{n-1}}(a_{n-1} \cdot R_{g_{n-1}}[\cdots T_{f_1}(a_1 \cdot R_{g_1}[k]) \cdots ])]) - T_fT_h( c_{n-1} \cdot [c_1\cdot R_h R_g ([k_{n-1}]))]  $$ is an element of $\underline{K}_n.$ Additionally, in $\underline{P}_{n} \underset{A} \boxtimes (\underline{A}/ \underline{K}^{\geq m})^{n-1}$, the image is $$T_{f_n}(a_n \cdot R_{g_n}[T_{f_{n-1}}(a_{n-1} \cdot R_{g_{n-1}}[\cdots T_{f_1}(a_1 \cdot R_{g_1}[k]) \cdots ])])    $$ because the second term is zero in $\underline{P}_{n} \underset{A} \boxtimes (\underline{A}/ \underline{K}^{\geq m})^{n-1}$. Being the image of an element of $\underline{P}_{n+1} \underset{A} \boxtimes (\underline{A}/ \underline{K}^{\geq m})^{n-1}$, the homology class in  $\underline{P}_{n} \underset{A} \boxtimes (\underline{A}/ \underline{K}^{\geq m})^{n-1}$ is the 0 class in $\newline \underline{\Tor}_{p}^{\underline{A}}(\underline{A}/ \underline{K}^{\geq m}, \underline{A}/ (\underline{K}^{\geq m})^{n-1}).$ Therefore the map on Tor is 0.
\end{proof}

\begin{corollary}\label{lem12.16}
Let $\underline{A}$ be a free $\underline{B}$-algebra, with the projection map $\underline{A} \to \underline{B}$ having kernel $\underline{K}$. Then for $p \neq 0$, there exists a short exact sequence 
\[\xymatrixrowsep{.5in}
\xymatrixcolsep{.3in}\xymatrix{
0 \ar[r]& \underline{\Tor}_{p+1}^{\underline{A}}(\underline{A}/ \underline{K}^{\geq m}, \underline{A}/ (\underline{K}^{\geq m})^{n-1}) \ar[r]& \\ \underline{\Tor}_{p}^{\underline{A}}(\underline{A}/ \underline{K}^{\geq m},(\underline{K}^{\geq m})^{n-1}/ (\underline{K}^{\geq m})^n) \ar[r]& \underline{\Tor}_{p}^{\underline{A}}(\underline{A}/ \underline{K}^{\geq m}, \underline{A}/ (\underline{K}^{\geq m})^n) \ar[r]& 0
}\]
coming from the derived long exact sequence for Tor.
\end{corollary}

\begin{lemma}\label{lem13.1}
Let $\underline{A}_* \to \underline{A}$ be a surjective map of simplicial $\underline{A}$-algebras, where $\underline{A}_n$ is a free $\underline{A}$-algebra for all $n$ and $\underline{A}$ is considered as the trivial simplicial $\underline{A}$-algebra. Let $\underline{K}_*$ be the kernel of this map, and suppose $\underline{K}_0 = 0$. Assume $\mathcal{O}$ has the flat-free condition. Then $$H_k[\underline{\Tor}_q^{\underline{A}_*}(\underline{A}_*/ \underline{K}_*^{\geq m}, \underline{A}_*/ (\underline{K}_*^{\geq m})^n)] \cong 0$$ for $k < n$ and $q \neq 0$.
\end{lemma}

First, we give a warning/definition of $\underline{\Tor}_q^{\underline{A}_*}(\underline{M}_*,\underline{W}_*)$ where $\underline{A}_*$ is a simplicial $\mathcal{O}$-Tambara functor and $\underline{M}_*,\underline{W}_*$ are simplicial $\underline{A}_*$ modules. It is incorrect to resolve $\underline{M}_*$ by some sequence of modules $\underline{P}_*$ that looks like the following diagram:
 \[\xymatrixrowsep{.5in}
\xymatrixcolsep{.25in}\xymatrix{
\cdots \ar[r]&\underline{P}_2  \ar[r]  \ar[d] &\underline{P}_1 \ar[r]  \ar[d] & \underline{P}_0 \ar[d] \\
\cdots \ar[r] &\underline{M}_2 \ar[r] &\underline{M}_1 \ar[r]  & \underline{M}_0  \\
}\]
where $\underline{P}_0$ is suitable $\underline{A}_0$-module and so on. This is only one part of the construction. The correct method is to get a resolution of $\underline{M}_*$ by projective $\underline{A}_*$-modules:
 \[\xymatrixrowsep{.5in}
\xymatrixcolsep{.25in}\xymatrix{
\cdots \ar[r]&\underline{P}_*^2  \ar[r] \ar[d] &\underline{P}_*^1 \ar[r]\ar[d] & \underline{P}_*^0  \ar[d] \\
\cdots \ar[r] &\underline{M}_* \ar[r] &\underline{M}_* \ar[r] & \underline{M}_*  \\
}\]
where every arrow is a simplicial map before the usual tensoring with $\underline{W}_*$
 \[\xymatrixrowsep{.5in}
\xymatrixcolsep{.25in}\xymatrix{
\cdots \ar[r]&\underline{P}_*^2 \underset{\underline{A}_*} \boxtimes \underline{W}_* \ar[r] \ar[d] &\underline{P}_*^1  \underset{\underline{A}_*} \boxtimes \underline{W}_* \ar[r] \ar[d] & \underline{P}_*^0   \underset{\underline{A}_*} \boxtimes \underline{W}_* \ar[d] \\
\cdots \ar[r] &\underline{M}_* \ar[r] &\underline{M}_* \ar[r]  & \underline{M}_*  \\
}\]
expanding out in the bottom coordinate gives the complicated but still somewhat illuminating diagram:
 \[\xymatrixrowsep{.1in}
\xymatrixcolsep{-.1in}\xymatrix{
&& \underline{P}_2^0   \underset{\underline{A}_2} \boxtimes \underline{W}_2  \ar[ld]   \ar[ddd] &&& \underline{P}_2^1   \underset{\underline{A}_2} \boxtimes \underline{W}_2  \ar[ld]  \ar[lll]   \ar[ddd]&&&  \underline{P}_2^2   \underset{\underline{A}_2} \boxtimes \underline{W}_2  \ar[ld] \ar[lll]   \ar[ddd] \\
&\underline{P}_1^0   \underset{\underline{A}_1} \boxtimes \underline{W}_1  \ar[ld]  \ar[ur]  \ar[ddd] &&&\underline{P}_1^1   \underset{\underline{A}_1} \boxtimes \underline{W}_1  \ar[ld]  \ar[lll] \ar[ur]  \ar[ddd]&&& \underline{P}_1^2   \underset{\underline{A}_1} \boxtimes \underline{W}_1   \ar[ld]  \ar[lll] \ar[ur]  \ar[ddd] \\
\underline{P}_0^0   \underset{\underline{A}_0} \boxtimes \underline{W}_0  \ar[ur]  \ar[ddd]&&& \underline{P}_0^1   \underset{\underline{A}_0} \boxtimes \underline{W}_0 \ar[lll] \ar[ur]  \ar[ddd]&&& \underline{P}_0^2   \underset{\underline{A}_0} \boxtimes \underline{W}_0 \ar[lll] \ar[ur]  \ar[ddd]\\
&&\underline{M}_2 \ar[ld]   &&& \underline{M}_2 \ar[ld]  \ar[lll]  &&& \underline{M}_2 \ar[ld] \ar[lll]   \\
&\underline{M}_1 \ar[ld]   \ar[ur] &&& \underline{M}_1 \ar[ld] \ar[lll] \ar[ur] &&& \underline{M}_1 \ar[ld]  \ar[lll] \ar[ur] \\
\underline{M}_0   \ar[ur] &&& \underline{M}_0 \ar[lll] \ar[ur] &&& \underline{M}_0  \ar[lll] \ar[ur] 
}\]
Now one would take the homology along the top coordinate in order to get $\underline{\Tor}$:
 \[\xymatrixrowsep{.1in}
\xymatrixcolsep{.1in}\xymatrix{
&&\underline{\Tor}_0^{\underline{A}_2}(\underline{M}_2,\underline{W}_2) \ar[ld]   &&&  \\
&\underline{\Tor}_0^{\underline{A}_1}(\underline{M}_1,\underline{W}_1) \ar[ld]  \ar[ur] &&& \underline{\Tor}_1^{\underline{A}_1}(\underline{M}_1,\underline{W}_1) \ar[ld]  \\
\underline{\Tor}_0^{\underline{A}_0}(\underline{M}_0,\underline{W}_0)   \ar[ur] &&&\underline{\Tor}_1^{\underline{A}_0}(\underline{M}_0,\underline{W}_0) \ar[ur] 
}\]
to realize $\underline{\Tor}$ as $\underline{A}_*$-module. We can take the homology again in the final direction
 \[\xymatrixrowsep{.1in}
\xymatrixcolsep{.01in}\xymatrix{
&&H_2[\underline{\Tor}_0^{\underline{A}_*}(\underline{M}_*,\underline{W}_*)]   &&   \\
&H_1[\underline{\Tor}_0^{\underline{A}_*}(\underline{M}_*,\underline{W}_*)]  && H_1[\underline{\Tor}_1^{\underline{A}_*}(\underline{M}_*,\underline{W}_*)]     \\
H_0[\underline{\Tor}_0^{\underline{A}_*}(\underline{M}_*,\underline{W}_*)]    && H_0[\underline{\Tor}_1^{\underline{A}_*}(\underline{M}_*,\underline{W}_*)]  }\]
to get the objects of the lemma.

\begin{proof}
The proof goes by induction on $k$. As a base case, when $k = 0$, because $$\underline{\Tor}_q^{\underline{A}_0}(\underline{A}_0/ \underline{K}_0^{\geq m}, \underline{A}_0/ (\underline{K}_0^{\geq m})^n) \cong \underline{\Tor}_q^{\underline{A}_0}(\underline{A}_0, \underline{A}_0) \cong 0 $$ implying the zeroth homology must be zero. Now assume the statement is true for $k-1$.

The short exact sequence of \cref{lem12.16} shows that we have a short exact sequence of simplicial modules $$ 0 \to
 \underline{\Tor}_{p+1}^{\underline{A}_*}(\underline{A}_*/ \underline{K}_*^{\geq m}, \underline{A}_*/ (\underline{K}_*^{\geq m})^{n-1}) \to \underline{\Tor}_{p}^{\underline{A}_*}(\underline{A}_*/ \underline{K}_*^{\geq m},(\underline{K}_*^{\geq m})^{n-1}/ (\underline{K}_*^{\geq m})^n)$$ $$ \to 
 \underline{\Tor}_{p}^{\underline{A}_*}(\underline{A}_*/ \underline{K}_*^{\geq m}, \underline{A}_*/ (\underline{K}_*^{\geq m})^n) \to 0$$ which gives rise to the long exact sequence in homology
\begin{eqnarray*}
 \cdots &\to&  
H_k[\underline{\Tor}_{p}^{\underline{A}_*}(\underline{A}_*/ \underline{K}_*^{\geq m}, (\underline{K}_*^{\geq m})^{n-1}/ (\underline{K}_*^{\geq m})^n)]\\ &\to& 
H_k[\underline{\Tor}_{p}^{\underline{A}_*}(\underline{A}_*/ \underline{K}_*^{\geq m}, \underline{A}_*/ (\underline{K}_*^{\geq m})^{n})]\\ &\to& 
H_{k-1}[\underline{\Tor}_{p+1}^{\underline{A}_*}(\underline{A}_*/ \underline{K}_*^{\geq m}, \underline{A}_*/ (\underline{K}_*^{\geq m})^{n-1})] \to \cdots
\end{eqnarray*}
the third term is zero by the induction hypothesis assuming $k < n$. So it remains to show that
$$H_i[\underline{\Tor}_q^{\underline{A}_*}(\underline{A}_*/\underline{K}_*^{\geq m}, (\underline{K}_*^{\geq m})^{n-1} / (\underline{K}_*^{\geq m})^{n} )] \cong 0.$$
$(\underline{K}_*^{\geq m})^{n-1} / (\underline{K}_*^{\geq m})^{n}$ is a flat $\underline{A}_*/\underline{K}_*$-module by the flat-free condition.  Then 
\begin{eqnarray*}
&&\underline{\Tor}_q^{\underline{A}_*}(\underline{A}_*/\underline{K}_*^{\geq m}, (\underline{K}_*^{\geq m})^{n-1} / (\underline{K}_*^{\geq m})^{n} ) \\&\cong& 
\underline{\Tor}_q^{\underline{A}_*}[\underline{A}_*/\underline{K}_*^{\geq m}, (\underline{A}_*/\underline{K}_* ) \underset{\underline{A}_*} \boxtimes (\underline{A}_* \underset{\underline{A}_*/\underline{K}_*} \boxtimes (\underline{K}_*^{\geq m})^{n-1} / (\underline{K}_*^{\geq m})^{n}) ]\\
&\cong& \underline{\Tor}_q^{\underline{A}_*}[\underline{A}_*/\underline{K}_*^{\geq m}, (\underline{A}_*/\underline{K}_* ) ] \underset{\underline{A}_*} \boxtimes (\underline{A}_* \underset{\underline{A}_*/\underline{K}_*} \boxtimes (\underline{K}_*^{\geq m})^{n-1} / (\underline{K}_*^{\geq m})^{n})\\
&\cong& \underline{\Tor}_q^{\underline{A}_*}[\underline{A}_*/\underline{K}_*^{\geq m}, (\underline{A}_*/\underline{K}_* ) ] \underset{\underline{A}_*/\underline{K}_*} \boxtimes (\underline{K}_*^{\geq m})^{n-1} / (\underline{K}_*^{\geq m})^{n}.
\end{eqnarray*}
From the short exact sequence
$$0 \to (\underline{K}_*^{\geq m}) \to \underline{A}_* \to \underline{A}_* / (\underline{K}_*^{\geq m}) \to 0$$ 
and the functor $(-) \underset{\underline{A}_*} \boxtimes \underline{A}_*/ (\underline{K}_*^{\geq m})^{n-1}$,
we obtain the isomorphism $$H_j[ (\underline{K}_*^{\geq m})^{n-1} / (\underline{K}_*^{\geq m})^{n}] \cong H_j[\underline{\Tor}_1^{\underline{A}_*}(\underline{A}_* / (\underline{K}_*^{\geq m}) , \underline{A}_*/ (\underline{K}_*^{\geq m})^{n})   ] \cong 0$$ for $j\leq  k-1$ by the induction hypothesis. Additionally, $$H_0[\underline{\Tor}_q^{\underline{A}_*}[\underline{A}_*/\underline{K}_*^{\geq m}, (\underline{A}_*/\underline{K}_* ) ] ] \cong 0 $$ for the same reason as the base case: $$\underline{\Tor}_q^{\underline{A}_0}[\underline{A}_0/\underline{K}_0^{\geq m}, (\underline{A}_0/\underline{K}_0 ) ]  \cong \underline{\Tor}_q^{\underline{A}_0}[\underline{A}_0, \underline{A}_0 ] \cong 0.$$
Then  \cref{lem8.18} shows that 
$$H_j[\underline{\Tor}_q^{\underline{A}_*}(\underline{A}_*/\underline{K}_*^{\geq m}, (\underline{K}_*^{\geq m})^{n-1} / (\underline{K}_*^{\geq m})^{n} )] \cong 0$$ for $j \leq k$ which completes the proof.
\end{proof}

We now give the main convergence theorem.
\begin{theorem}
Let $\underline{A}_*$ a free $\underline{A}$-algebra, with a surjective map onto the trivial simplicial $\underline{A}$-algebra $\underline{A}$. Let $\underline{K}_*$ be the kernel. Assume that $\underline{K}_0 = 0$ and $\mathcal{O}$ has the flat-free condition. Then $$H_k(\underline{K}_*^{\geq nm}) \cong 0$$ if $k < n.$
\end{theorem}

\begin{proof}
The proof goes by induction on $n$. Suppose $n = 1$, then the zeroth homology vanishes because $\underline{K}_0 = 0.$ Now we assume the statement is true for $n$ to prove for $n+1$.
Given the short exact sequence of simplicial $\underline{A}$-modules $$0 \to \underline{K}^{\geq m}_* \to \underline{A}_* \to \underline{A}_*/\underline{K}^{\geq m}_* \to 0 $$ consider the functor $(-) \underset{\underline{A}_*} \boxtimes (\underline{K}_*^{\geq m})^n$. Because $m$ is a sharing number for $n$, the image of $\underline{K}_*^{\geq m} \underset{\underline{A}_*} \boxtimes (\underline{K}_*^{\geq m})^n$ is $(\underline{K}_*^{\geq m(n+1)})$ so we get the short exact sequence which is part of the long exact sequence for Tor:
$$0 \to \underline{\Tor}_1^{\underline{A}_*}(\underline{A}_*/\underline{K}_*^{\geq m}, (\underline{K}_*^{\geq m})^n) \to (\underline{K}_*^{\geq m}) \underset{\underline{A}_*} \boxtimes (\underline{K}_*^{\geq m})^n  \to (\underline{K}_*^{\geq m(n+1)}) \to 0.$$ This gives rise to the following part of the long exact sequence in homology:
$$H_k( (\underline{K}_*^{\geq m}) \underset{\underline{A}_*} \boxtimes (\underline{K}_*^{\geq m})^n ) \to H_k((\underline{K}_*^{\geq m})^{n+1}) \to H_{k-1}(\underline{\Tor}_1^{\underline{A}_*}(\underline{A}_*/\underline{K}_*^{\geq m}, (\underline{K}_*^{\geq m})^n) ).$$
The third element is zero because from the short exact sequence $$0 \to (\underline{K}_*^{\geq m})^n) \to \underline{A}_* \to \underline{A}_* /(\underline{K}_*^{\geq m})^n  \to 0 $$ we get an isomorphism $$\underline{\Tor}_1^{\underline{A}_*}(\underline{A}_*/\underline{K}_*^{\geq m}, (\underline{K}_*^{\geq m})^n) \cong \underline{\Tor}_2^{\underline{A}_*}(\underline{A}_*/\underline{K}_*^{\geq m}, \underline{A}_*/(\underline{K}_*^{\geq m})^n)$$ and $$H_{k-1}(\underline{\Tor}_2^{\underline{A}_*}(\underline{A}_*/\underline{K}_*^{\geq m}, \underline{A}_*/(\underline{K}_*^{\geq m})^n)) \cong 0$$ for $k  \leq n$ by \cref{lem13.1}. All that remains to show is that $H_k( (\underline{K}_*^{\geq m}) \underset{\underline{A}_*} \boxtimes (\underline{K}_*^{\geq m})^n ) \cong 0$ for $k \leq n.$

Consider the homomorphisms of simplicial $\underline{A}_*/\underline{K}_*$ modules:
$$d_j:  \underline{K}_*^{\geq m } \underset{\underline{A}_*/\underline{K}_*}  \boxtimes (\underset{\underline{A}_*/\underline{K}_*}  {\overset{j} \boxtimes} \underline{K}_*)\underset{\underline{A}_*/\underline{K}_*} \boxtimes (\underline{K}_*^{\geq m})^n \to  \underline{K}_*^{\geq m } \underset{\underline{A}_*/\underline{K}_*} \boxtimes (\underset{\underline{A}_*/\underline{K}_*}  {\overset{j-1} \boxtimes}  \underline{K}_*)\underset{\underline{A}_*/\underline{K}_*} \boxtimes (\underline{K}_*^{\geq m})^n $$
with $$d_j( T_p(k_1 \otimes \cdots \otimes k_{j+1})) = \sum_{i=0}^j T_p(k_1 \otimes \cdots \otimes k_i k_{i+1} \otimes \cdots \otimes k_{j+1})$$ on the first component and the identity on the second component. This makes a chain complex $L$. Let $\underline{Q}_*^j$ be the image of this map and $\underline{P}_*^j$ be the kernel. Then the above describes a chain complex of simplicial modules.
\begin{lemma}
The homology of this complex computes $\underline{\Tor}_q^{\underline{A}_*}(\underline{K}_*^{\geq m}, (\underline{K}_*^{\geq m})^n).$
\end{lemma}
\begin{proof}
The first homology $H_0[L]$ is the quotient of $ \underline{K}_*^{\geq m} \underset{\underline{A}_*/\underline{K}_*} \boxtimes (\underline{K}_*^{\geq m})^n$ by the elements of the form $x_0 \otimes x_1 x_2 - x_0 x_1 \otimes x_2$ where $x_1$ is an element of $\underline{K}_*.$ This is precisely $ \underline{K}_*^{\geq m} \underset{\underline{A}_*} \boxtimes (\underline{K}_*^{\geq m})^n.$

We now need to show that every other homology Mackey functor vanishes. First let us consider a similar chain complex, $S_*$ 
$$d_j:  \underline{K}_*^{\geq m } \underset{\underline{A}_*/\underline{K}_*}  \boxtimes (\underset{\underline{A}_*/\underline{K}_*}  {\overset{j} \boxtimes}  \underline{K}_*)\underset{\underline{A}_*/\underline{K}_*} \boxtimes \underline{A}_* \to  \underline{K}_*^{\geq m } \underset{\underline{A}_*/\underline{K}_*} \boxtimes (\underset{\underline{A}_*/\underline{K}_*}  {\overset{j-1} \boxtimes}  \underline{K}_*)\underset{\underline{A}_*/\underline{K}_*} \boxtimes \underline{A}_* $$
define $s_n(x_0 \otimes \cdots \otimes x_n \otimes x_{n+1}) \mapsto x_0 \otimes \cdots \otimes x_{n+1} \otimes 1$ if $x_{n+1}$ is in $\underline{K}_*$ and is 0 if $x_{n+1}$ is in $\underline{A}_*/\underline{K}_*$. This defines a map from $S_k \to S_{k+1}$ and has the property that $s_{n-1} \circ d_n - d_{n+1} \circ s_n = (-1)^n \Id.$ Now applying the functor $(-) \underset{\underline{A}_*} \boxtimes (\underline{K}_*^{\geq m})^n$ to this projective complex gives the desired lemma.

\end{proof}
We have the short exact sequence $$0 \to \underline{Q}_*^1 \to \underline{K}_*^{\geq m} \underset{\underline{A}_*/\underline{K}_*} \boxtimes (\underline{K}_*^{\geq m})^n \to \underline{K}_*^{\geq m} \underset{\underline{A}_*} \boxtimes (\underline{K}_*^{\geq m})^n \to 0 $$
So the homology of the second term is 0 using  \cref{lem8.18} up to $n+1$ and the inductive hypothesis. So now it remains to show that $H_{i-1}(\underline{Q}_*^1) \cong 0$ for $i \leq n.$

We have the short exact sequence: 
$$0 \to \underline{Q}_*^{q+1} \to  \underline{P}_*^{q} \to \underline{\Tor}_q^{\underline{A}_*}(\underline{K}_*^{\geq m}, (\underline{K}_*^{\geq m})^n) \to 0 $$
which gives rise to the homology long exact sequence $$H_k[\underline{Q}_*^{q+1}] \to  H_k[\underline{P}_*^{q}] \to H_k[\underline{\Tor}_q^{\underline{A}_*}(\underline{K}_*^{\geq m}, (\underline{K}_*^{\geq m})^n)] $$
Now again by the short exact sequence  
$$0 \to (\underline{K}_*^{\geq m})\to \underline{A}_* \to \underline{A}_* /(\underline{K}_*^{\geq m})  \to 0 $$ we see that $$\underline{\Tor}_q^{\underline{A}_*}(\underline{K}_*^{\geq m}, (\underline{K}_*^{\geq m})^n) \cong \underline{\Tor}_{q+1}^{\underline{A}_*}(\underline{A}_*/\underline{K}_*^{\geq m}, (\underline{K}_*^{\geq m})^n) \cong  \underline{\Tor}_{q+2}^{\underline{A}_*}(\underline{A}_*/\underline{K}_*^{\geq m}, \underline{A}_*/(\underline{K}_*^{\geq m})^n)$$ so the third term in the homology long exact sequence is 0 for $k < n$ by \cref{lem13.1}, which shows that $H_k[\underline{Q}_*^{q+1}] \to  H_k[\underline{P}_*^{q}]$ is surjective.

We have another short exact sequence:
$$0 \to \underline{P}_*^{q} \to  \underline{K}_*^{\geq m } \underset{\underline{A}_*/\underline{K}_*}  \boxtimes (\underset{\underline{A}_*/\underline{K}_*} {\overset{q} \boxtimes} \underline{K}_*)\underset{\underline{A}_*/\underline{K}_*} \boxtimes (\underline{K}_*^{\geq m})^n  \to \underline{Q}_*^{q} \to 0 $$
which gives rise to the homology long exact sequence $$ H_{k+1}[ \underline{K}_*^{\geq m } \underset{\underline{A}_*/\underline{K}_*}  \boxtimes (\underset{\underline{A}_*/\underline{K}_*}  {\overset{q} \boxtimes} \underline{K}_*)\underset{\underline{A}_*/\underline{K}_*} \boxtimes (\underline{K}_*^{\geq m})^n ] \to H_{k+1}[\underline{Q}_*^{q}] \to H_k[\underline{P}_*^{q}] $$ the first term is 0 for $k<n+1$, again by repeated use \cref{lem8.18}, implying that $H_{k+1}[\underline{Q}_*^{q}] \to H_k[\underline{P}_*^{q}] $ is injective. 

Putting the two statements together, that $H_k[\underline{Q}_*^{q+1}] \to  H_k[\underline{P}_*^{q}]$ is surjective and $H_{k+1}[\underline{Q}_*^{q}] \to H_k[\underline{P}_*^{q}] $ is injective, if $H_k[\underline{Q}_*^{q+1}] \cong 0$ then $H_{k+1}[\underline{Q}_*^{q}] \cong 0$. Since $H_0[\underline{Q}_*^{i}] \cong 0$ for all $i$, we get the desired result.
\end{proof}

\begin{corollary}
Suppose we have a simplicial resolution $\underline{B}_*$ of the $\underline{A}$-algebra $\underline{B}$ with $\underline{B}_0 \cong \underline{A}$ and $\mathcal{O}$ has the flat-free condition. Then the kernel $\underline{J}_*$ of the map $$\underline{B}_* \underset{\underline{A}} \boxtimes \underline{B} \to \underline{B} $$ has the property $H_k[\underline{J}_*^{\geq mn} \underset{\underline{B}} \boxtimes \underline{W}] \cong 0$ for $k < m$ and any $\underline{B}$-module $\underline{W}.$
\end{corollary}

\begin{proof}
$\underline{B}_n \underset{\underline{A}} \boxtimes \underline{B}$ is a free $\underline{B}$-algebra, and $\underline{J}_0 =0$, so $H_k[\underline{J}_*^{\geq mn}] \cong 0$ for $k \leq n.$ We proceed by induction. The statement is clearly true for $k = 0$ because $\underline{J}_0 = 0.$ Now $H_{k-1}[\underline{J}_*^{\geq nm} \underset{\underline{B}} \boxtimes (-)]$ is exact on the left (because it is always 0). Because $\underline{J}^{\geq mn}$ is a free simplicial $\underline{B}$-algebra, we have the corollary by the next lemma.
\end{proof}

\begin{lemma}
Let $\underline{L}_*$ be a complex of projective $\underline{A}$-modules, and let $\underline{W}$ be a $\underline{A}$-module. Then there exists an isomorphism of $\underline{A}$-modules $$H_n[\underline{L}_*] \underset{\underline{A}} \boxtimes \underline{W} \cong H_n[\underline{L}_* \underset{\underline{A}} \boxtimes \underline{W}]$$ if the functor $H_{n-1}[\underline{L}_* \underset{\underline{A}} \boxtimes (-)]$ is exact on the left.
\end{lemma}

\begin{proof}
Given a short exact sequence of $\underline{A}$-modules:
$$0 \to \underline{W}' \to \underline{W} \to \underline{W}'' \to 0$$ we get a exact sequence of simplicial $\underline{A}$-modules, using $\underline{L}_*$ is projective
$$0 \to \underline{L}_* \underset{\underline{A}} \boxtimes \underline{W}' \to \underline{L}_* \underset{\underline{A}} \boxtimes \underline{W} \to \underline{L}_* \underset{\underline{A}} \boxtimes \underline{W}'' \to 0$$
which gives the long exact sequence in homology
$$H_n[\underline{L}_* \underset{\underline{A}} \boxtimes \underline{W}] \to H_n[\underline{L}_* \underset{\underline{A}} \boxtimes \underline{W}''] \to H_{n-1}[\underline{L}_* \underset{\underline{A}} \boxtimes \underline{W}'] \to H_{n-1}[\underline{L}_* \underset{\underline{A}} \boxtimes \underline{W}]  $$
so $H_{n-1}[\underline{L}_* \underset{\underline{A}} \boxtimes (-)]$ being exact on the left is the same as $H_{n}[\underline{L}_* \underset{\underline{A}} \boxtimes (-)]$ being exact on the right.

There is always a map $$H_n[\underline{L}_*] \underset{\underline{A}} \boxtimes \underline{W} \to H_n[\underline{L}_* \underset{\underline{A}} \boxtimes \underline{W}].$$ This is an isomorphism if $\underline{W}$ is a free $\underline{A}$-module.  But every $\underline{W}$ is a quotient of 2 free modules: $$\underline{F}' \to \underline{F} \to \underline{W} \to 0 $$
which gives the following commutative diagram, and upon using the 5-lemma realizes the desired isomorphism.
 \[\xymatrixrowsep{.4in}
\xymatrixcolsep{.4in}\xymatrix{
H_n[\underline{L}_*] \underset{\underline{A}} \boxtimes \underline{F}'  \ar[r] \ar[d]^\cong& H_n[\underline{L}_*] \underset{\underline{A}} \boxtimes \underline{F}  \ar[r] \ar[d]^\cong& H_n[\underline{L}_*] \underset{\underline{A}} \boxtimes \underline{W}   \ar[r] \ar[d]& 0 \ar[d]^\cong \\
H_n[\underline{L}_* \underset{\underline{A}} \boxtimes \underline{F}']  \ar[r] & H_n[\underline{L}_* \underset{\underline{A}} \boxtimes \underline{F}] \ar[r] & H_n[\underline{L}_* \underset{\underline{A}} \boxtimes \underline{W}]  \ar[r] & 0 
}
\]
\end{proof}

\section{Computations}
As an example of the spectral sequence in action, we compute
$\underline{\Tor}^{\underline{A}[x_G]}(\underline{A},\underline{A})$ in two ways.

First we do $\underline{\Tor}^{\underline{A}[x_G]}(\underline{A},\underline{A})$ via projective resolution. The free resolution of modules terminates. Take 
 \[\xymatrixrowsep{.48in}
\xymatrixcolsep{.95in}\xymatrix{
0 \to \underline{A}[x_G]\{b_G \} \ar[r]^{b_G \mapsto t_q} & \underline{A}[x_G]\{p_G,q_e \} \ar[r]^{p_G \mapsto t_m}_{q_e \mapsto n_e - \overline{n_e}} & \underline{A}[x_G]\{m_e, n_e\} \ar[r]_(.6){m_e \mapsto w_e-\overline{w_e}}^(.6){n_e \mapsto rw_e - R(w_G)} & \cdots \\  
\underline{A}[x_G]\{w_G, w_e\} \ar[r]_(.5){w_e \mapsto rR(z_G)-R(u_G)}^(.5){w_G \mapsto nz_G - xu_G} &\underline{A}[x_G]\{z_G, u_G\} \ar[r]_{u_G \mapsto ny_G}^{z_G\mapsto xy_G} & \underline{A}[x_G]\{y_G\} \ar[r]^{y_G \mapsto1} & \underline{A}
}
\]
changing the units to $\underline{A}$ gives 
 \[\xymatrixrowsep{.5in}
\xymatrixcolsep{1in}\xymatrix{
0 \to \underline{A}\{b_G \} \ar[r]^{b_G \mapsto t_q} & \underline{A}\{p_G,q_e \} \ar[r]^{p_G \mapsto t_m}_{q_e \mapsto n_e - \overline{n_e}} & \underline{A}\{m_e, n_e\} \ar[r]_(.6){m_e \mapsto w_e-\overline{w_e}}^(.6){n_e \mapsto - R(w_G)} & \cdots \\  
\underline{A}\{w_G, w_e\} \ar[r]_(.5){w_e \mapsto -R(u_G)}^(.5){w_G \mapsto 0} &\underline{A}\{z_G, u_G\} \ar[r]_{u_G \mapsto 0}^{z_G\mapsto 0} & \underline{A}\{y_G\} \ar[r]^{y_G \mapsto1} & \underline{A}
}
\]
which upon taking the homology gives 
 \[\xymatrixrowsep{.5in}
\xymatrixcolsep{.5in}\xymatrix{
\txt{deg $\geq$ 3} & \txt{deg 2} & \txt{deg 1} & \txt{deg 0} \\
0  \ar@/^-1pc/[d]_{R^G_e}  &\Z\{w_G\}  \ar@/^-1pc/[d]_{R^G_e}     &\Z[t]/(t^2-2t)\{z_G\} \oplus \Z\{u_G\}   \ar@/^-1pc/[d]_{R^G_e}  &\Z[t]/(t^2-2t)\{y_G\}   \ar@/^-1pc/[d]_{R^G_e} \\
0  \ar@/^-1pc/[u]_{T^G_e}  &  0  \ar@/^-1pc/[u]_{T^G_e} & \Z\{R(z_G)\} \oplus 0 \ar@/^-1pc/[u]_{T^G_e}  & \Z\{R(y_G)\}  \ar@/^-1pc/[u]_{T^G_e}  
}
\]
Now we do the resolution in algebras. The free resolution is in fact `generated' by a single element, in that all the other generators are degenerate:
 \[\xymatrixrowsep{.5in}
\xymatrixcolsep{.5in}\xymatrix{
\cdots \ar@<-.5ex>[r] \ar@<.5ex>[r] \ar@<-1.5ex>[r] \ar@<1.5ex>[r]&  (\underline{A}[x_G])[z^0_G,z^1_G] \ar[r] \ar@<-.5ex>[r] \ar@<.5ex>[r] &  (\underline{A}[x_G])[z_G] \ar@<-.5ex>[r]_{z_G \mapsto 0} \ar@<.5ex>[r]^{z_G \mapsto x} & \underline{A}[x_G] \ar[r]^{x_G \mapsto 0} & \underline{A}
}
\]
restricting to the different degrees and taking the homology fills out the following second page of the spectral sequence:

 \[\xymatrixrowsep{.5in}
\xymatrixcolsep{.5in}\xymatrix{
0 & 0 & \underline{K} &0 & 0&\cdots \\
& 0 & \underline{K} &0  &0& \cdots \\
& & 0 &\underline{A} & 0  & \cdots\\
& & & \underline{A} & 0 & \cdots
}
\]
where $\underline{A}$ is the burnside Mackey functor and $\underline{K}$ is the Mackey functor 
 \[\xymatrixrowsep{.5in}
\xymatrixcolsep{.5in}\xymatrix{
\Z \ar@/^-1pc/[d]_{R^G_e}  \\
0  \ar@/^-1pc/[u]_{T^G_e} }
\]